\documentclass{amsart}
\usepackage[OT1]{fontenc}
\usepackage{amsmath,amssymb,amsthm,mathtools}
\usepackage[expansion=false]{microtype}
\usepackage[colorlinks=true, linkcolor=blue, citecolor=blue, urlcolor=blue]{hyperref}
\usepackage{xcolor} 

\numberwithin{equation}{section}
\newtheorem{theorem}{Theorem}[section]
\newtheorem{proposition}[theorem]{Proposition}
\newtheorem{lemma}[theorem]{Lemma}
\newtheorem{corollary}[theorem]{Corollary}
\newtheorem{knowntheorem}{Theorem}

\newtheorem{openproblem}{Open Problem}

\newtheorem{conjecture}{Conjecture}

\theoremstyle{remark}
\newtheorem{remark}[theorem]{Remark}

\newcommand{\R}{\mathbb R}
\newcommand{\Z}{\mathbb Z}
\newcommand{\T}{\mathbb R/\mathbb Z}
\newcommand{\dd}{\,d}

\DeclareMathOperator{\ind}{ind}

\title[Classification of critical points of heat kernels]
{Classification of critical points of heat kernels on flat tori 
and Riemann theta functions}

\author{Zhijie Chen}
\address{Department of Mathematical Sciences, Yau Mathematical Sciences Center,
Tsinghua University, Beijing, 100084, China }
\email{zjchen2016@tsinghua.edu.cn}

\author{Chen-Gang Long}
\address{Yau Mathematical Sciences Center,
Tsinghua University, Beijing, 100084, China}
\email{lcg21@tsinghua.org.cn}

\author{Shihong Zhang}
\address{Yau Mathematical Sciences Center,
Tsinghua University, Beijing, 100084, China}
\email{shihong-zhang@tsinghua.edu.cn}

\begin{document}

\begin{abstract}
For any two-dimensional flat torus $E_\tau:=\mathbb{C}/(\Z+\Z\tau)$
and any $t>0$, 
we prove that the heat kernel $p_{\tau}(z,t)$ has either four or six critical points on $z\in E_{\tau}$ (i.e., including four trivial critical
points and at most one pair of nontrivial critical points).
This is interestingly
 consistent with Lin-Wang's remarkable result (Ann. Math. 2010) for the Green function.
The key idea is to show that
 any nontrivial critical point (if it exists) is a nondegenerate global minimum point.
Furthermore,
  we classify the types of all critical points, as well as the splitting 
  phenomena that occur when the number of critical points changes. Applications to Riemann theta functions are also given. In particular, we resolve a conjecture of Luo and Wei (Arch. Ration. Mech. Anal. 2022).

\end{abstract}

\maketitle


\section{Introduction}
\label{sec:introduction}
Let $\tau=a+ib\in\mathbb H:=\{\tau\in\mathbb{C}:\operatorname{Im}\tau>0\}$,
$\Lambda_\tau:=\Z+\Z\tau$ and  $E_\tau:=\mathbb{C}/\Lambda_\tau$
be a flat torus with the Euclidean metric. 
The heat kernel $p_{\tau}(z,t)$ on the torus $E_{\tau}$
 is the unique function satisfying
\begin{equation}\label{eq1-1}
	\left\{\begin{aligned}
		\partial_t p_{\tau}(z,t)&=\triangle_z p_{\tau}(z,t),\quad (z,t)\in E_{\tau}\times (0,+\infty),\\
		p_{\tau}(z,0)&=\delta_0(z),\quad z\in E_{\tau}.
	\end{aligned}\right.
 \end{equation}
Here, $\delta_0$ denotes the Dirac measure at $0$, and the initial condition is understood in the sense of distributions.
The heat kernel $p_{\tau}(z,t)$ has an
explicit expression 
\begin{equation}\label{eq:heat}
 p_\tau(z,t)=\frac1{4\pi t}\sum_{\omega\in\Lambda_{\tau}} e^{-\frac{|z+\omega|^2}{4t}}=\frac1{4\pi t}
 \sum_{m,n\in\Z}
 e^{-\frac{|z+m+n\tau|^2}{4t}},\quad z\in E_{\tau},\;t>0.
\end{equation}
See \cite[Lemma 2.1]{LWZ} for a proof.
For an introduction to
heat kernels on manifolds, we refer to the textbook of Grigor'yan \cite{Gr09}.

In this paper, 
we study the critical points of $p_{\tau}(z,t)$
with respect to $z$ for fixed  $\tau\in\mathbb{H}$ and $t>0$, i.e.,
\begin{equation*}
       \nabla_z p_{\tau}(z,t)=0.
\end{equation*}
Since the heat kernel $p_{\tau}(\cdot,t)$ is even and doubly periodic,
the half-period points $0,\frac{1}{2},\frac{\tau}{2},\frac{1+\tau}{2}$
are always critical points of $p_{\tau}(\cdot,t)$, and
 $-z$ is also a critical point of $p_{\tau}(\cdot,t)$ if $z$ 
is.
We call the half-period points $0,\frac{1}{2},\frac{\tau}{2},\frac{1+\tau}{2}$
\emph{trivial} critical points, and call
 non-half-period critical points (i.e., $z\neq-z$ in $E_{\tau}$) \emph{nontrivial} critical points.
Therefore, nontrivial critical points
must appear in pairs if they exist.
A basic question is:

\vskip 0.1in
\noindent{\bf Problem:} For fixed \(\tau\in\mathbb H\) and \(t>0\), does 
\(p_\tau(\cdot,t)\) have
nontrivial critical points? If so, how many are there?
\vskip 0.1in

One of our motivations comes from Lin and Wang \cite{Li10} on the study of critical points of the Green function $G(z)=G(z;\tau)$
 on $E_\tau$, which
is defined by 
\begin{equation}
       -\triangle G(z)=\delta_0(z)-\frac{1}{|E_{\tau}|},\quad\int_{E_{\tau}} G(z)dz=0,
\end{equation}
where $|E_{\tau}|$ is the area of the torus $E_{\tau}$. Since $G(z)$ is even and has a singularity at $0$, $G(z)$ has always three trivial critical points $\frac{1}{2},\frac{\tau}{2},\frac{1+\tau}{2}$, and nontrivial critical points also appear in pairs.
Lin and Wang \cite{Li10} developed deep connections among the critical points of $G(z)$,
the bubbling phenomenon of mean field equations, and elliptic functions. Consequently,  they proved a remarkable result that $G(z)$ has at most one pair of nontrivial critical points. See also Bergweiler and Eremenko \cite{Bergweiler2016} for a new proof. Further properties of the unique pair of nontrivial critical points were studied in \cite{Chen2025,Li17}.
We summarize these results as follows.
\begin{knowntheorem}[\cite{Bergweiler2016,Chen2025,Li10,Li17}]\label{thm:Green function}
       $G(z;\tau)$ has at most one pair of nontrivial critical points, i.e., has either three or five critical points (depending on the choice of $\tau$). 
       
       Furthermore, if the unique pair of nontrivial  critical points $\pm q$ exists, then
       $\pm q$ are both nondegenerate global minimum points of $G(z;\tau)$.
\end{knowntheorem}

A complete characterization of 
$\tau$ such that $G(z;\tau)$ has nontrivial critical points was given in \cite{Bergweiler2016,Chen2018}.
Later, subsequent works \cite{Chai2015,Chen-Fu-Lin2026,Chen2021,Lin2026,Lin2017,Lin2020} 
 studied the generalized multiple Green function 
$$G_n(z_1,\cdots,z_n;\tau):=\sum_{j<k}G(z_j-z_k;\tau)-n\sum_{j=1}^n G(z_j;\tau)$$
due to its applications in mean field equations and classical Lam\'{e} equations. 
Among them, Chai, Lin and Wang \cite{Chai2015} proved that $G_n(\cdot;\tau)$ has at most $2n+1$ trivial critical points, but how many nontrivial critical points $G_n(\cdot;\tau)$ might have remained completely open there.
Recently, the first author and Lin \cite{Chen2026} proved that
 $G_n(\cdot;\tau)$ has at most $n$ pairs of nontrivial critical points, and each nontrivial critical point is always nondegenerate and contributes the same degree $(-1)^{n+1}$.

 The heat kernel $p_{\tau}(z,t)$ and the Green function $G(z;\tau)$ satisfy the following important formula
 \begin{equation}\label{eq1-4}
        G(z;\tau)=\int_0^{\infty} \left(p_{\tau}(z,t)-\frac{1}{|E_{\tau}|}\right)dt,\quad z\in E_{\tau}\setminus\{0\}.
 \end{equation}
 See \cite{LWZ} for a proof, where
 the second author, Wei, and Zou proved that the heat kernel $p_{\tau}(\cdot,t)$
 also has at most one pair of nontrivial critical points on 
 rectangular, hexagonal, and rhombic tori like the Green function. 
 Specifically,

\begin{knowntheorem}[\cite{LWZ}]
\label{thm:known-heat}
For any $t>0$, the following statements hold.
\begin{enumerate}
\item If $\tau=ib$ with $b>0$, 
the heat kernel $p_{\tau}(\cdot,t)$ has exactly four critical
points.
\item If $\tau=\frac12+ib$ with $b>0$, there exist
$\frac{\sqrt{3}}{6}<b_1(t)<\frac{1}{2}<b_2(t)<\frac{\sqrt{3}}{2}$ such that 
$p_{\tau}(\cdot,t)$ has exactly four critical
points for $b\in[b_1(t),b_2(t)]$ and six
for $b\notin[b_1(t),b_2(t)]$.
\end{enumerate}
\end{knowntheorem}

It is worth mentioning that previous studies on critical points of 
the Green function $G$ relied crucially on the formula
\begin{equation}\label{eq:Green function expression}
       -4\pi \frac{\partial G }{\partial z}(z;\tau)=\zeta(z;\tau)-\alpha z-\beta\bar{z},
\end{equation}
where $z=x+iy\in\mathbb C$ is considered as the complex variable,  $\zeta(z;\tau)$ is the Weierstrass zeta function, and $\alpha=2\zeta(\frac{1}{2};\tau)-\frac{\pi}{b}$,
$\beta=\frac{\pi}{b}$.
However, although the heat kernel $p_{\tau}$ has a connection with the Green function via formula \eqref{eq1-4},
the derivative $\frac{\partial}{\partial z}p_{\tau}(z,t)$ cannot satisfy a formula similar to \eqref{eq:Green function expression},
and the methods of studying critical points of the Green function in \cite{Bergweiler2016,Chen2025,Chen2018,Li10,Li17} cannot work for the heat kernel, so 
studying the critical points of $p_{\tau}(\cdot,t)$ requires new approaches. The approach of proving Theorem \ref{thm:known-heat} in \cite{LWZ} used certain symmetries of  rectangular, hexagonal, and rhombic tori, and cannot work for general tori.
An open problem about general tori was proposed in \cite{LWZ}.

\begin{openproblem}[\cite{LWZ}, Open problem~7.1 (1)-(2)]\label{Open Problem}
       Fix any $\tau\in\mathbb{H}$ and $t>0$.
       \begin{enumerate}
       \item Does \(p_\tau(\cdot,t)\) have at most one pair of nontrivial critical points?
       \item Are all nontrivial critical points of \(p_\tau(\cdot,t)\) minimum points?
       \end{enumerate}
\end{openproblem}

Our second motivation comes from periodic particle energies.
The heat kernel $p_{\tau}$ can be viewed as the potential energy
on the lattice $\Lambda_{\tau}$ with respect to the Gaussian potential
$f(|\cdot|^2)=e^{-\frac{1}{4t} |\cdot|^2}$:
 \begin{equation}
      p_\tau(z,t)= \frac1{4\pi t} \sum_{\omega\in\Lambda_{\tau}} f(|z+\omega|^2).
 \end{equation}
This potential energy expresses the interaction energy between particles
on $\Lambda_{\tau}$ under the Gaussian potential and
appears frequently in the structure theory of matter
at low temperatures, including in phenomena such as Ginzburg-Landau
vortices, Bose-Einstein condensates, Wigner crystals, and superconductivity, see
e.g. \cite{ He80,LuoWei2022,Mu02,Sa12,Ve79}.

In this physical context, many current research results on the critical points
 of the heat kernel are experimental or numerical results. In view of
mathematics, there are only a few rigorous theoretical results, and they 
mainly focus on the extremal points.
Baernstein \cite{Ba95} and Faulhuber \cite{Fa21} proved that
on rectangular and hexagonal tori, the maximum
of $p_{\tau}(\cdot,t)$ occurs at $0$, while the minimum occurs at 
$\frac{1+\tau}{2}$ and $\frac{1+\tau}{3}$, respectively. 
In addition, B\'etermin and Petrache \cite{Be17} proved that, for any deep hole $c$ and any fixed point $z\in E_\tau$ that is not a deep hole,
$p_\tau(c,t)<p_\tau(z,t)$ for sufficiently small $t>0$.
\cite{Re17} established translation inequalities for $p_\tau(z,t)$.
Moreover, \cite{Be21} and \cite{Mo88} proved that, 
after scaling $p_{\tau}$ in \eqref{eq:heat} to volume-one tori,
 the  maxmin and minimax
\begin{equation*}
     \max_{\tau}  \min_{z} p_{\tau}(z,t),   \quad \min_{\tau}  \max_{z} p_{\tau}(z,t) 
\end{equation*}
are attained by the hexagonal torus ($\tau=e^{\frac{\pi i}{3}}$).

For the study of general critical points of \(p_{\tau}(z,t)\),
Burchard and Martínez \cite{Bu25}
 studied the monotonicity of \(p_{\tau}(\cdot,t)\) along geodesics on flat tori and used this property to classify all critical points of 
 \(p_{\tau}(z,t)\) on the equilateral torus. This monotonicity property has also been studied on other manifolds in \cite{Al23}.
Recently, \cite{Ba25} proved that \(z\) is a uniform
critical point for all \(t>0\) if and only if 
every nonempty inhomogeneous shell forms a spherical
1-design.
In addition, the second author, Wei and Zou \cite{LWZ} 
gave a complete classification of the critical points of \(p_{\tau}(z,t)\) on 
rectangular,
hexagonal and rhombic tori, and proved that the critical 
points on rectangular and hexagonal tori are all uniform.
However, for general tori and $t>0$, the number and types of the critical points
 remain completely open at present.

Our main result is the following theorem, which gives a complete 
positive answer to Open Problem \ref{Open Problem}
and characterizes the types of all critical points. The key idea is to use theta functions to estimate $p_{\tau}$
and analyze the types of critical points. This approach is essentially different from that of proving Theorem \ref{thm:Green function} for the Green function.

\begin{theorem}\label{thm:global-classification}
For any $\tau\in\mathbb{H}$ and $t>0$, the 
 heat kernel $p_\tau(\cdot,t)$ on the flat torus $E_\tau$
 has exactly 4 or 6 critical points, i.e.,
 $p_\tau(\cdot,t)$ has at most one pair of nontrivial critical points.
Moreover,  critical points have the following precise properties:
\begin{enumerate}

\item $0$ is always the unique and nondegenerate 
global maximum point of $p_\tau(\cdot,t)$.

\item If $p_\tau(\cdot,t)$ has exactly 4 critical
points, then the half-period points $\frac{1}{2},\frac{\tau}{2},\frac{1+\tau}{2}$
consist of two nondegenerate saddle points and one minimum point.
Furthermore, the minimum point is either nondegenerate
or degenerate with the Hessian $D^2_z p_{\tau}(z,t)$ having
a positive eigenvalue.

\item If $p_\tau(\cdot,t)$ has exactly 6 critical
points, then the half-period points $\frac{1}{2},\frac{\tau}{2},\frac{1+\tau}{2}$
are all nondegenerate saddle points, 
and the unique pair of nontrivial critical points $z_\pm=\pm(\nu+\mu\tau)\notin\{0,\frac{1}{2},\frac{\tau}{2},\frac{1+\tau}{2}\}$ are both
nondegenerate global minimum points.
Furthermore, when $\tau$ lies in the 
fundamental domain \eqref{eq:reduced-region}, the
coordinates $\pm(\nu,\mu)$ lie in $\mathcal T_+$ and
$-\mathcal T_+$, respectively, where
\begin{equation*}
       \mathcal T_+
:=\left\{
(\nu,\mu) :
0<\mu,\nu<\frac12,\quad
\mu+\nu>\frac12
\right\}.
\end{equation*}

\end{enumerate}
\end{theorem}

\begin{remark}\label{remark}
       Here are some comments on Theorem \ref{thm:global-classification}.
\begin{enumerate}
       \item Theorem \ref{thm:global-classification} shows that 
       $p_\tau(\cdot,t)$ has at most one pair of nontrivial critical points,
       which are always nondegenerate minimum points if they exist.
       This result is highly consistent with Theorem \ref{thm:Green function} for the Green function. Furthermore, it was proved in \cite[Theorem 1.3]{Chen2018} that the coordinates $(\nu,\mu)$ of a nontrivial critical point 
$z=\nu+\mu\tau$ of the Green function $G(z;\tau)$ also lie in 
$\pm\mathcal T_+$, which is the same region
 as in Theorem \ref{thm:global-classification} (3).
       However, the heat kernel and the Green function have quite different structures,
       so that their proofs are quite different. An interesting problem is whether
      there is a deeper theory to explain this phenomenon. 
\item Theorem \ref{thm:global-classification}  implies that a degenerate
critical point can only be a trivial critical point,
and there exists at most one such point.
Therefore, when the number of critical points of $p_\tau$
 changes between 4 and 6, a
 splitting phenomenon occurs at the 
 degenerate nonzero trivial critical point.
 Specifically, the degenerate minimum point splits into one
 nondegenerate saddle point and two nondegenerate minimum points.
 We will precisely characterize this degeneracy and
  splitting in Theorem \ref{thm:smooth-degeneracy}.

\end{enumerate}

\end{remark}

There is a deep connection between the heat kernel $p_{\tau}$ and the 
Riemann theta function.
Indeed,
let $\operatorname{Sym}_2^+(\R)$ be the set of real
$2\times2$ positive definite symmetric matrices.
For any $S\in \operatorname{Sym}_2^+(\R)$ and $ q\in\R^2/\Z^2$,
define a two-dimensional associated Riemann theta function
\begin{equation}\label{eq:Riemann theta function}
 \Theta_S(q):=\sum_{\mathbf{n}\in\Z^2}e^{-4\pi^2\mathbf{n}^TS\mathbf{n}}e^{2\pi i \mathbf{n}\cdot q}.
\end{equation}
Then we have the formula
\begin{equation}\label{eq:heat kernel theta function equality}
       p_{\tau}(z,t)=\frac{1}{|\det A|} \Theta_{tA^{-1}A^{-T}}(A^{-1}z),\quad\forall z\in E_{\tau},\;t>0,
\end{equation}
where $\tau=a+bi\in\mathbb{H}$,\begin{equation}\label{eq:AS}A=\begin{pmatrix}
       1&a\\
       0&b
 \end{pmatrix},\qquad S=tA^{-1}A^{-T}=\frac{t}{b^2}\begin{pmatrix}
       a^2+b^2&-a\\
       -a&1
 \end{pmatrix}.\end{equation} The proof of this formula will be given in \eqref{eq:matrix-identification} in Section \ref{subsec heat kernel representation}.
In this paper, we also regard \(z=x+yi\in\mathbb{C}\) as the vector \((x,y)^T\in\R^2\).

This formula \eqref{eq:heat kernel theta function equality} implies that \(z\) is a critical point of \(p_{\tau}(z,t)\)
if and only if \(q=A^{-1}z\) is a critical point of \(\Theta_S(q)\),
and they have the same type.
Moreover, \(z\) is nontrivial if and only if
\(q=A^{-1}z\notin \tfrac12\Z^2/\Z^2\), i.e., \(q\) is not a half-integer point.
Thus, we obtain the following consequence concerning the critical points of the Riemann
 theta function $\Theta_S(q)$.

\begin{corollary}
       For any $S\in \operatorname{Sym}_2^+(\R)$, the Riemann theta 
       function $\Theta_S(q)$ has exactly four or six critical points on
        $q\in\R^2/\Z^2$, with the types and degeneracies described in 
        Theorem~\ref{thm:global-classification}.
\end{corollary}

 \begin{proof}
       For any $S=\begin{pmatrix}
       \alpha&\beta\\
       \beta&\gamma
       \end{pmatrix}\in \operatorname{Sym}_2^+(\R)$, let $\tau=a+bi$ and $t>0$
        be defined by \eqref{eq:coordinate change} via $\alpha,\beta,\gamma$,
         then \eqref{eq:heat kernel theta function equality} holds; see Section 2.1.
        Hence this corollary follows directly from Theorem~\ref{thm:global-classification}.
 \end{proof}

Theorem \ref{thm:global-classification} also resolves a conjecture of Luo and Wei
\cite[Conjecture~2.1 (a)]{LuoWei2022}.
In their study of Mueller--Ho conjecture \cite{Mu02} related to competing systems of Bose-Einstein condensates, they considered
\[
 J(\tau;x,y)
 :=\sum_{m,n\in\mathbb Z}
 e^{-\frac{\pi}{\operatorname{Im}\tau}|m\tau-n|^2}
 \cos\bigl(2\pi(mx+ny)\bigr),
 \quad \tau\in\mathbb H,\quad
 (x,y)\in\mathbb R^2/\mathbb Z^2,
\]
and proposed \cite[Conjecture~2.1]{LuoWei2022} concerning the possible numbers of critical points for $J(\tau;\cdot,\cdot)$. While Mueller--Ho conjecture was solved recently by Luo and Wei \cite{LuoWei2026},
and the (b)-(c) part of  \cite[Conjecture~2.1]{LuoWei2022} was solved by the second author, Wei and Zou \cite{LWZ}, the (a) part of  \cite[Conjecture~2.1]{LuoWei2022}, which is recalled below, remains open so far.

\begin{conjecture}[Luo and Wei \cite{LuoWei2022}, Conjecture~2.1 (a)]\label{conjecture B}
      For each fixed $\tau$, the function
$J(\tau;\cdot,\cdot)$ has exactly four or six critical points.
\end{conjecture}

Here we can 
give a positive answer to this conjecture.

\begin{corollary}\label{cor:luo-wei}
For every $\tau\in\mathbb H$, the function
$J(\tau;\cdot,\cdot)$ has exactly four or six critical points
on $\mathbb R^2/\mathbb Z^2$.
Moreover, any critical point outside
$\frac12\mathbb Z^2/\mathbb Z^2$ is a nondegenerate global
minimum point.
\end{corollary}

\begin{proof}
Write $\tau=a+ib$ and let $A, S$ be defined by \eqref{eq:AS}.
Since $A^{-1}(x+y\tau)=(x,y)^T$ and
$(m,n)S(m,n)^T=\frac{t}{b^2}|m\tau-n|^2,$
formula \eqref{eq:heat kernel theta function equality}
and the definition of $\Theta_S$ give
\[
 b\,p_\tau(x+y\tau,t)
 =\sum_{m,n\in\mathbb Z}
 e^{-\frac{4\pi^2t}{b^2}|m\tau-n|^2}
 e^{2\pi i(mx+ny)}.
\]
Taking $t=b/(4\pi)$ and pairing the terms of
$(m,n)$ and $(-m,-n)$, we obtain
\[
 J(\tau;x,y)
 =b\,p_\tau\left(x+y\tau,\frac{b}{4\pi}\right).
\]
Clearly the change of coordinates $(x,y)\mapsto x+y\tau$ preserves critical points, their types,
and nondegeneracy. The conclusion follows from
Theorem~\ref{thm:global-classification}.
\end{proof}
       
 Note from Remark \ref{remark} (2) that the degeneracy and splitting 
 phenomenon of $p_{\tau}$ occurs only at the degenerate nonzero trivial
  critical point, or equivalently, at the degenerate nonzero half-integer 
  point of $\Theta_S(q)$. For any nonzero half-integer point
 $q\in \left(\tfrac12\Z^2/\Z^2\right)\setminus\{0\}$,
we define the degeneracy set 
\begin{equation}
  \mathcal{D}_q:=\{S\in \operatorname{Sym}_2^+(\R):\det D^2\Theta_S(q)=0\}\subset \operatorname{Sym}_2^+(\R).  
\end{equation}
Our second theorem describes the degeneracy set $\mathcal{D}_q$ and 
the splitting phenomenon near $\mathcal{D}_q$.

\begin{theorem}\label{thm:smooth-degeneracy}
For any nonzero half-integer point $q_0\in\tfrac12\Z^2/\Z^2\setminus\{0\}$, the
degeneracy set $\mathcal{D}_{q_0}\neq \emptyset$ and it is a real analytic
 embedded hypersurface in $\operatorname{Sym}_2^+(\R)$. Moreover,
\begin{equation*}
       \mathcal{D}_{q_0}\cap \mathcal{D}_{\tilde{q}_0}=\emptyset,\quad\forall q_0,\tilde{q}_0\in \tfrac12\Z^2/\Z^2\setminus\{0\},\quad q_0\neq \tilde{q}_0.
\end{equation*}
For any $S_0\in\mathcal D_{q_0}$, there exists a neighborhood $U$ of $S_0$
such that the degeneracy set $\mathcal{D}_{q_0}$ locally divides
$U$ into two regions:
$$U_{+}=\{S\in U \;:\; \det D^2\Theta_S(q_0)>0\},\quad U_{-}=\{S\in U \;:\; \det D^2\Theta_S(q_0)<0\},$$
namely $U\setminus \mathcal D_{q_0}=U_{+}\cup U_{-}$. Furthermore, 
\begin{equation}\label{eq:p tau critical number}
       \begin{aligned}
          &|\{z\in E_{\tau}: \nabla_z p_{\tau}(z,t)=0\}|\\
          = &   
 |\{q\in \R^2/\Z^2 : \nabla_q  \Theta_S(q)=0\}| 
 =
\begin{cases}
6,&S\in U_-,\\
4,&S\in U_+\cup (U\cap \mathcal D_{q_0}).
\end{cases}    
       \end{aligned} 
\end{equation}
Here  $\tau,t,S$ satisfy the expression \eqref{eq:heat kernel theta function equality}.
\end{theorem}

\begin{remark}
       Based on formula \eqref{eq:heat kernel theta function equality},
       the change of $S$ in $\Theta_S$ includes changes in both the torus parameter $\tau$ and the time parameter $t$.
       Therefore, Theorem \ref{thm:smooth-degeneracy}
        describes the degeneracy set and splitting
       phenomenon of $p_{\tau}(\cdot,t)$ with respect to 
       $\tau$ and $t$ simultaneously.
       By fixing one parameter, we can restrict the critical-point count to the other parameter. The corresponding  splitting conclusion holds when this parameter slice is transverse to the degeneracy set.
\end{remark}

\medskip
\noindent\textbf{Outline of the proof.} Now we briefly introduce the main ideas of the proof. The key idea is to prove that all nontrivial critical points are 
nondegenerate minimum points,  calculate the index of each critical point, 
and  apply the Poincar\'e--Hopf theorem. To this goal, 
by \eqref{eq:heat kernel theta function equality}, it suffices to study
the critical points of $\Theta_S$ on $\R^2/\Z^2$, with
$S\in\operatorname{Sym}_2^+(\R)$.
The main idea is to combine local perturbations of $S$ with a global
argument based on $\operatorname{tr}S$.
For every real symmetric matrix $M$, the family $\Theta_S$ satisfies
the heat equation identity
\[
 \left.\frac{d}{ds}\Theta_{S+sM}(q)\right|_{s=0}
 =\sum_{i,j=1}^2 M_{ij}\partial_{q_iq_j}\Theta_S(q).
\]
This identity relates variations in the matrix parameter to spatial
derivatives and is the basis of the perturbation arguments in
Section~\ref{sec:local-degeneracy}.

In Section~\ref{sec:local-degeneracy}, we show that any nontrivial degenerate critical point
can be perturbed into a nearby nontrivial nondegenerate saddle point
while strictly decreasing the matrix trace.
At a degenerate nonzero half-integer point, the local analysis gives
an alternative: either such a perturbation exists, or all nearby
nontrivial critical points are nondegenerate minimum points.

In Section \ref{sec:global-exclusion},
to exclude nontrivial nondegenerate saddle points globally, let
$\mathcal B$ be the set of matrices
$S\in\operatorname{Sym}_2^+(\R)$ for which $\Theta_S$ has such a point.
Suppose that $\mathcal B\ne\varnothing$, and consider
\[
 c_*:=\inf_{S\in\mathcal B}\operatorname{tr}S.
\]
We choose $S_j\in\mathcal B$ with $\operatorname{tr}S_j\to c_*$ and
corresponding nontrivial nondegenerate saddle points $q_j$.
After passing to a subsequence, we have
$S_j\to S_*\ge0$ and $q_j\to q_*$.
The small-time and large-$b$ estimates in
Section~\ref{sec:uniform-estimates} show that $S_*>0$.
Thus, these estimates prevent the sequence from approaching the
boundary of the positive definite cone.

The openness of $\mathcal B$ and the implicit function theorem then
force the limiting critical point $q_*$ to be degenerate.
The local analysis in Section~\ref{sec:local-degeneracy} either
produces a matrix in $\mathcal B$ with trace strictly smaller than
$c_*$, or excludes nearby nontrivial nondegenerate saddle points.
The first possibility contradicts the definition of $c_*$, while
the second contradicts the convergence of $(S_j,q_j)$ to $(S_*,q_*)$.
This proves that $\mathcal B$ is empty, as established in
Lemma~\ref{lem:no-nontrivial-saddle}.

Once nontrivial nondegenerate saddle points have been excluded,
the same perturbation results rule out nontrivial degenerate
critical points. Hence every nontrivial critical point has index
$+1$. Together with the analysis at the half-integer points,
the Poincar\'e--Hopf theorem gives the four-or-six classification
in Theorem~\ref{thm:global-classification}.
Compactness and evenness then show that the nontrivial pair,
when present, consists of nondegenerate global minimum points.
Finally, the positive fourth-order term at a degenerate half-integer
point yields the analytic
structure of the degeneracy sets and the splitting described in
Theorem~\ref{thm:smooth-degeneracy}.

The organization of this paper is as follows. In Section \ref{sec:preliminaries},
we derive some expressions and basic properties of the heat kernel.
In Section \ref{sec:uniform-estimates}, we prove that when $b$ is sufficiently large or $t$ is sufficiently 
small, every nontrivial critical point of $p_\tau(z,t)$ 
is a nondegenerate minimum point.
In Section \ref{sec:local-degeneracy}, we construct nondegenerate saddle points near a degenerate
 critical point by decreasing the matrix trace, which
  will be used to exclude the existence
of nontrivial nondegenerate saddle points.
Finally, in Section \ref{sec:global-exclusion}, we prove the main results Theorem \ref{thm:global-classification}
and Theorem \ref{thm:smooth-degeneracy}.
 Throughout the paper, we denote by $c, c_k, C, C_k>0$ various constants (maybe different in different places) independent of the corresponding varying parameters.

\section{Preliminaries}
\label{sec:preliminaries}

\subsection{Representation of the heat kernel on one and two dimensional tori.}
\label{subsec heat kernel representation}
According to 
\cite[Lemma~2.1]{LWZ}, the heat kernel $p_\tau(z,t)$
 on the torus $E_\tau$ that solves \eqref{eq1-1} is expressed as
\begin{equation}\label{eq2-1}
 p_\tau(z,t)= \frac1{4\pi t}
 \sum_{\omega\in \Lambda_\tau}
 \exp\left(-\frac{|z+\omega|^2}{4t}\right)    =\frac1{4\pi t}
 \sum_{m,n\in\Z}
 \exp\left(-\frac{|z+m+n\tau|^2}{4t}\right).
\end{equation}

The lattice $\Lambda_\tau$ has some important
invariances
\begin{equation}\label{eq:modular invariant}
\Lambda_{\tau+1}=\Lambda_\tau,\qquad \Lambda_{-\frac{1}{\tau}}=\tau^{-1}\Lambda_\tau,\qquad \Lambda_{-\bar\tau}=\overline{\Lambda_\tau},
\end{equation}
which identify the lattices up to rotation, scaling and reflection, as referenced in 
\cite[Section 2]{LongLuoZou26}. Since the critical points of the heat kernel \eqref{eq2-1} 
are in one-to-one correspondence under the 
rotation, scaling ($z$ and $t$ are scaled by $\lambda$ and $\lambda^2$, respectively, for a scaling factor $\lambda>0$), and reflection of the lattice $\Lambda_\tau$,
as referenced in 
\cite[Section 2]{LWZ}, 
without loss of generality, 
we only need to  consider
\begin{equation}\label{eq:reduced-region}
 \tau\in\mathcal F
 :=\{a+ib:0\leq a\leq1/2,\ b>0,\ a^2+b^2\geq1\}.
\end{equation}

By using the Poisson Summation Formula, \cite[Section 2]{LWZ} 
gives an alternative representation of $p_{\tau}$, namely
\begin{equation}\label{eq2-7}
		p_{\tau}(z,t)=\dfrac{1}{4\pi t}\sum_{\omega\in\Lambda_{\tau}} \exp\left(-\frac{|z+\omega|^2}{4t}\right)=\frac{1}{|\Lambda_{\tau}|} \sum_{\omega\in\Lambda_{\tau}^{\bot}}e^{-4\pi^2 t|\omega|^2}e^{2\pi iz\cdot \omega},
\end{equation}
where $|\Lambda_{\tau}|=\mathrm{vol}(\mathbb{C}/\Lambda_{\tau})=\mathrm{vol}(E_\tau)=b$ 
is the covolume of $\Lambda_{\tau}$, equivalently the area of $E_{\tau}$,
and the dual lattice $\Lambda^{\bot}$ of the lattice $\Lambda$ is 
\begin{equation}\label{eq:dual lattice}
\Lambda^{\bot} \coloneqq \{z\in \R^2\;:\; z\cdot \omega\in\Z,\;\forall \omega\in\Lambda  \}.
\end{equation}
Here $\cdot$ is the standard inner product in $\R^2$.

Next we use the Riemann theta function to express the heat kernel
 $p_{\tau}$.
Let $\operatorname{Sym}_2^+(\R)$  be the set of real
$2\times2$ positive definite symmetric matrices.
For any $S\in \operatorname{Sym}_2^+(\R)$ and $ q\in\R^2/\Z^2$,
define an associated Riemann theta function
\begin{equation}\label{eq:def US}
 \Theta_S(q)=\Theta(q;4\pi iS):=\sum_{\mathbf{n}\in\Z^2}e^{-4\pi^2\mathbf{n}^TS\mathbf{n}}e^{2\pi i \mathbf{n}\cdot q}.
\end{equation}
Here,  $$\Theta(q; \Omega):=\sum_{\mathbf{n}\in\Z^2}\exp\left(\pi i \mathbf{n}^T \Omega \mathbf{n}+2\pi i \mathbf{n} \cdot  q\right)$$
denotes the standard Riemann theta function \cite[\S21]{DLMF}.
The function $\Theta_S$ is jointly real analytic with respect to 
 $S$ and $q$.

Suppose $A=\begin{pmatrix}
      1&a\\
      0&b
\end{pmatrix}$ is a basis matrix of $\Lambda_{\tau}$, i.e., 
$\Lambda_{\tau}=A\Z^2=\{A \mathbf{n}: \mathbf{n}\in \Z^2\}$. Then, from the expression \eqref{eq:dual lattice}, the transpose inverse matrix 
 $A^{-T}$ is a basis matrix of the dual lattice $\Lambda_{\tau}^{\bot}$.
Thus, $$\{|\omega|^2:\omega\in \Lambda_{\tau}^{\bot}\}=\{\mathbf{n}^T A^{-1}A^{-T} \mathbf{n}\;:\;\mathbf{n}\in \Z^2 \}$$
and $|\Lambda_{\tau}|=b=|\det A|$. Therefore, using the expressions
\eqref{eq2-7} and \eqref{eq:def US}, we obtain 
the Riemann theta function expression for the heat kernel
\begin{equation}\label{eq:matrix-identification}
       p_{\tau}(z,t)=\frac{1}{|\det A|} \Theta_{tA^{-1}A^{-T}}(A^{-1}z).
\end{equation}

Conversely, for any $S=\begin{pmatrix}
       \alpha&\beta\\
       \beta&\gamma
       \end{pmatrix}\in\operatorname{Sym}_2^+(\R)$ and $q\in\R^2/\Z^2$, 
       we choose 
\begin{equation}\label{eq:coordinate change}
a=-\frac{\beta}{\gamma},\quad
b=\frac{\sqrt{\det S}}{\gamma},\quad
t=\frac{\det S}{\gamma},\quad
A=
\begin{pmatrix}
1&a\\
0&b
\end{pmatrix},
\end{equation}
then $\tau=a+bi$ and $z=Aq$ satisfy $p_{\tau}(z,t)=\frac{1}{|\det A|} \Theta_{S}(q)$.

On the one-dimensional torus 
$\R/\Z$, using the Poisson Summation Formula as in \eqref{eq2-7},
the heat kernel on 
$\R/\Z$ is 
\begin{equation}\label{eq2-9}
       \theta(x;t)=\frac1{\sqrt{4\pi t}}
       \sum_{m\in\Z}e^{-(x-m)^2/(4t)}
       =\sum_{n\in\Z}e^{-4\pi^2tn^2}e^{2\pi inx},
\end{equation}
which is the unique solution satisfying 
\begin{equation}
	\left\{\begin{aligned}
		\partial_t u(x,t)&=\partial_{x}^2 u(x,t),\quad (x,t)\in \R/\Z\times (0,+\infty),\\
		u(x,0)&=\delta_0(x),\quad x\in \R/\Z.
	\end{aligned}\right.
 \end{equation}
This $\theta(x;t)=\theta_3(\pi x;4\pi it)$ is the standard
Jacobi theta function \cite[\S20]{DLMF}. 
 
As in \cite[\S20]{DLMF}, this $\theta(x;t)$ can be rewritten as
\begin{equation}\label{eq:product}
 \theta(x;t)=\prod_{m=1}^{\infty}(1-e^{-4\pi^2 t(2m)})
       (1+2e^{-4\pi^2 t(2m-1)}\cos(2\pi x)+e^{-4\pi^2 t(4m-2)}).
\end{equation}

From expression \eqref{eq:product}, by direct differentiations 
and using the signs of the trigonometric functions,
we obtain sign properties of $\theta(x;t)$ and 
$\partial_x^2\log\theta(x;t)$ for all $t>0$:
\begin{equation}\label{eq:circle-signs}
 \begin{gathered}
       \partial_x^2\log\theta(x;t)<0,~x\in [0,1/4];\\
 \theta_x(x;t)<0,~x\in (0,1/2);\quad  \theta_x(x;t)>0,~x\in (1/2,1);\\
 \theta_{xx}(0;t)<0<\theta_{xx}(1/2;t).
 \end{gathered}
\end{equation}
A similar proof can be found in \cite[Section 3]{LWZ}, and the second inequality was already stated in
\cite[Lemma~3.2]{LWZ}.

\subsection{Schur complement.}
In this subsection, we present the important Schur complement method, 
which will be frequently used in this paper to determine whether
 the Hessian is positive definite.
 The introduction to the Schur complement  can be found in the textbooks 
 \cite[Appendix A.5.5]{Bo2004} and \cite[Chapter 1]{Zhang2005}.

 Suppose that $u(x,y)\in C^{\infty}(\R^2)$ and $(x_0,y_0)$ is a critical point of $u$
 satisfying $u_{yy}(x_0,y_0)\ne 0$. By the implicit function theorem,
 there exists a unique local function $h(x)$ such that 
 $u_y(x,h(x))=0$. Let $f(x)=u(x,h(x))$,
then a direct differentiation yields
 \begin{equation}\label{eq:schur}
       \begin{aligned}
           & f'(x)=u_x(x,h(x)),\quad
  f''(x)=u_{xx}-\frac{u_{xy}^2}{u_{yy}},\\
&  \det D^2u(x,h(x)):=\det D^2 u(x,y)\big|_{y=h(x)}=u_{yy}(x,h(x))f''(x).  
       \end{aligned}
 \end{equation}
 On the curve $(x,h(x))$, $f'(x)=0$
 is equivalent to $\nabla u(x,h(x))=0$.
By \eqref{eq:schur}, the Hessian $D^2u(x,h(x))$ is positive definite if and only if $u_{yy}(x,h(x))>0$ and $f''(x)>0$. We will use this idea in the following proofs.

\section{Nondegenerate minimality of nontrivial critical points in 
two limiting cases}\label{sec:uniform-estimates}

The key step in the proof of Theorem \ref{thm:global-classification} is the 
use of the Poincar\'e--Hopf theorem, so we need to know the properties of the
 Hessian at critical points, especially to prove that all nontrivial critical
  points are nondegenerate minimum points.

In this section, we prove two special cases in Proposition 
\ref{prop:cusp-estimate} and Proposition \ref{prop:uniform-small},
namely that when $b$ is sufficiently large or $t$ is sufficiently small, 
all nontrivial critical points are nondegenerate minimum points.
These two special cases will be used in Section \ref{sec:global-exclusion}
 for the proof of the general case.

\subsection{Nondegenerate minimality when \(b\) is sufficiently large}\label{subsec:reduction}

We first prove the case of sufficiently large $b$
in the following Proposition \ref{prop:cusp-estimate}.

\begin{proposition}\label{prop:cusp-estimate}
Let $J=[\alpha,\beta]\subset(0,\infty)$ be a closed interval.
 There exists $B_J<\infty$ such that for any $t\in J$ and $\tau=a+bi$ with
 $0\leq a\leq1/2$ and $b\geq B_J$,
every nontrivial critical point of
 $p_{\tau}(\,\cdot,t)$ is a nondegenerate local minimum point.
 \end{proposition}

 To prove Proposition \ref{prop:cusp-estimate}, we make a single 
 summation over \(p_{\tau}\) to give a precise characterization of 
 the critical points of \(p_{\tau}\) as \(b\to\infty\).
 For any $z=x+yi\in E_{\tau}$ and $t>0$, using the expressions \eqref{eq2-7}
and \eqref{eq2-9}, we have
\begin{equation}\label{eq:rows}
       \begin{aligned}
              \sqrt{4\pi t} p_\tau(x+yi,t)&=\frac1{\sqrt{4\pi t}}
       \sum_{m,n\in\Z}
       \exp\left(-\frac{|x+yi-m-n(a+bi)|^2}{4t}\right)\\
       &=\sum_{n\in\Z}e^{-\frac{(y-bn)^2}{4t}}\theta(x-an;t).
       \end{aligned}
\end{equation}
In this subsection, we regard \(a,t\) as 
fixed parameters, \(b\) as a varying parameter,
and \(x,y\) as variables.
In the following, we will omit the fixed parameters \(a,t\).

Next, we make a change of variables, replacing \((b,y)\) by
$r=\frac{t}{b^2}$ and $ \mu=\frac yb$,
so that, based on \eqref{eq:rows},
 we can write
\begin{equation}\label{eq:scaled-vertical}
       Z_r(x,\mu)=Z_r(x,\mu;a,t):=\sum_{n\in\Z}e^{-\frac{(\mu-n)^2}{4r}}\theta(x-an;t)=\sqrt{4\pi t} p_\tau(x+yi,t).
\end{equation}
Therefore, $x+yi$ is a critical point of $p_{\tau}(\cdot,t)$ if and only if 
$\nabla_{(x,\mu)} Z_r(x,\mu)=0$.
Thus, in this subsection, we will study the solutions of $\nabla_{(x,\mu)} Z_r(x,\mu)=0$.

According to periodicity, any $z=x+yi\in E_{\tau}$  
can be chosen with $x\in [0,1]$ and $y\in[-\frac{1}{4}b, \frac{3}{4}b]$,
 so we only need to consider  $\mu\in \left[-\frac{1}{4},\frac{3}{4}\right]$.
Note that $t\in J$ is bounded, so as $b\to \infty$, we have
 $r=\frac{t}{b^2}\to 0^+$.

 Now we study the partial derivatives $\partial_{\mu}  Z_r(x,\mu)$, 
 $\partial_{\mu}^2  Z_r(x,\mu)$ and the zero curve satisfying 
 $\partial_{\mu}  Z_r(x,\mu)=0$. 
 
 \begin{lemma}\label{le:partial Zr} 
       There is $r_0>0$ sufficiently small such that for any $r\in (0, r_0)$, there are exactly two 
        smooth functions $\mu_j(\cdot;r): \mathbb{R}/\mathbb Z\to [-\frac{1}{4},\frac{3}{4}]$, $j=0,1$,
       such that  $\partial_{\mu}  Z_r(x,\mu)=0$
       for $\mu\in [-\frac{1}{4},\frac{3}{4}]$ if and only if $\mu=\mu_j(x;r)$ for some $j=0,1$.
       
     Furthermore,
        \begin{equation}\label{eq:mu0xr}
           \partial_{\mu}^2  Z_r(x,\mu_0(x;r)) <0, \quad  \mu_0(x;r)=O(e^{-c/r}),
        \end{equation}
        and
        \begin{equation}\label{eq:3-4}
                \partial_{\mu}^2  Z_r(x,\mu_1(x;r))>0, \quad \mu_1(x;r)=\tfrac12+O(r).     
        \end{equation}
 Here, $c>0$ is a constant independent of $r$, $x\in\R/\Z$, $a\in[0,1/2]$ and $t\in J$, and we omit the fixed parameters $a,t$.

 \end{lemma}

 \begin{proof} Note that $Z_r(x+1,\mu)=Z_r(x,\mu)$, we only need to consider $x\in [0,1]$ below.
 By direct differentiations, we have
 \begin{equation} \label{eq:vertical-first}
       (Z_r)_{\mu}(x,\mu) := \partial_{\mu}  Z_r(x,\mu)=-\frac1{2r}\sum_{n\in\Z} (\mu-n) e^{- \frac{(\mu-n)^2}{4r}}\theta(x-an;t),
 \end{equation}
 \begin{equation}  \label{eq:vertical-second}
       (Z_r)_{\mu\mu}(x,\mu) := \partial_{\mu}^2  Z_r(x,\mu)=\sum_{n\in\Z} \left(\frac{(\mu-n)^2}{4r^2}-\frac1{2r}\right) e^{-\frac{(\mu-n)^2}{4r}}\theta(x-an;t).
 \end{equation}

 We divide $\mu$ into two cases.

 {\bf Case 1: $|\mu|\leq1/4$.} In this case, as $r\to 0^+$, 
 the main part of $(Z_r)_{\mu} (x,\mu)$ is the $n=0$ term.
 Specifically, 
for any $n\neq 0$, we have 
$(\mu-n)^2-\mu^2\geq n^2/2$, thus we obtain the estimate
\begin{equation}\label{eq:decay estimate}
      0< \frac{e^{- \frac{(\mu-n)^2}{4r}}\theta(x-an;t)}{e^{-\frac{\mu^2}{4r}}\theta(x;t)}\leq
       e^{-\frac{n^2}{8r}}
       \frac{\theta(x-an;t)}{\theta(x;t)} 
\le C_J e^{-\frac{n^2}{8r}}.
\end{equation}
Here, the last inequality uses the fact that the smooth positive function
 $\theta(x;t)$
on the compact set $(x,t)\in \R/\Z \times J$ satisfies 
$$0<\min_{x,t}\theta(x;t)<\max_{x,t}\theta(x;t)<\infty$$ and that $C_J$ is a constant 
depending only on the set $J$.

Thus, by using formulas \eqref{eq:vertical-first} and \eqref{eq:decay estimate},
as $r\to 0^+$,
we obtain that 
\begin{equation}\label{eq:first decay estimate}
       \begin{aligned}
            \frac{(Z_r)_{\mu} (x,\mu) }{e^{-\frac{\mu^2}{4r}}\theta(x;t)}
            =&-\frac1{2r}
            \left[\mu+\sum_{n\ne0}(\mu-n)  \dfrac{e^{- \frac{(\mu-n)^2}{4r}}\theta(x-an;t)}{ e^{-\frac{\mu^2}{4r}}\theta(x;t)}\right]\\=&-\frac{\mu}{2r}+O(r^{-1}e^{-c/r}).  
       \end{aligned}
\end{equation}
Similarly, by using formulas \eqref{eq:vertical-second} and \eqref{eq:decay estimate},
as $r\to 0^+$, we have
\begin{equation}\label{eq:second decay estimate}
       \frac{(Z_r)_{\mu\mu}(x,\mu)}{e^{-\frac{\mu^2}{4r}}\theta(x;t)}
 =\frac{\mu^2}{4r^2}-\frac1{2r}
 +O(r^{-2}e^{-c/r}).
\end{equation}
Therefore, for $r>0$ sufficiently small, it follows from \eqref{eq:first decay estimate} that
$$ (Z_r)_{\mu} (x,1/4) <0< (Z_r)_{\mu} (x,-1/4),\quad \forall x\in [0,1],$$
so there is $\mu_0=\mu_0(x;r)\in (-\frac14,\frac14)$ such that $(Z_r)_{\mu} (x,\mu_0)=0$. Using \eqref{eq:first decay estimate} again, 
we have $\mu_0=O(e^{-c/r})$ and then \eqref{eq:second decay estimate} becomes 
\begin{equation}\label{eq:second zero estimate}
      \frac{(Z_r)_{\mu\mu}(x,\mu_0)}{e^{-\frac{\mu_0^2}{4r}}\theta(x;t)}
 =-\frac1{2r}
 +O(r^{-2}e^{-c/r})<0
\end{equation}
by taking $r>0$ smaller if necessary. This shows that every zero of $(Z_r)_{\mu} (x,\cdot)=0$ lies in $|\mu|\leq Ce^{-c/r}$, and $(Z_r)_{\mu\mu}(x,\cdot)<0$ on this interval, so the zero $\mu_0=\mu_0(x;r)$ is unique for every $x\in[0,1]$. This proves \eqref{eq:mu0xr}. The implicit function theorem implies that $\mu_0(\cdot;r)\in C^{\infty}(\mathbb R/\mathbb Z)$.

{\bf Case 2: $1/4\leq\mu\leq3/4$.} In this case,
for any $n\in\Z$, we have $|\mu-n|\geq1/4$. Hence, based on \eqref{eq:vertical-second},
when $r<\frac{1}{32}$,  we obtain 
\begin{equation}\label{eq:derivative increasing}
       (Z_r)_{\mu\mu}(x,\mu)\geq \left( \frac1{64r^2}-\frac1{2r} \right) Z_r(x,\mu)>0,\quad\forall x\in [0,1],\ \frac{1}{4}\leq\mu\leq \frac{3}{4}.
\end{equation}
Thus $(Z_r)_{\mu} (x,\mu)$ is strictly increasing with respect to $\mu$.

In this case, as $r\to 0^+$, 
the main part of $(Z_r)_{\mu} (x,\mu)$ is the $n=0,1$ terms.
In fact, at $\mu=1/4$, 
similar to the estimates \eqref{eq:decay estimate} and \eqref{eq:first decay estimate}
we have 
\[
\begin{aligned}
       \frac{(Z_r)_{\mu} (x,1/4)}{e^{-\frac{1}{64r}}\theta(x;t)}
=-\frac1{8r}+O(r^{-1}e^{-1/(8r)})<0,\quad r\to 0^+.
\end{aligned}
\]
Also, at $\mu=3/4$,
\[
\begin{aligned}
       \frac{(Z_r)_{\mu} (x,3/4)}{e^{-\frac{1}{64r}}\theta(x-a;t)}
=\frac1{8r}+O(r^{-1}e^{-1/(8r)})>0,\quad r\to 0^+.
\end{aligned}
\]
Therefore, for $r>0$ sufficiently small,  $$ (Z_r)_{\mu} (x,1/4) <0< (Z_r)_{\mu} (x,3/4),\quad \forall x\in [0,1].$$
Together with \eqref{eq:derivative increasing}, 
there exists a unique zero point $\mu_1=\mu_1(x;r)\in\left(1/4, 3/4 \right)$ such 
that 
\begin{equation*}
       (Z_r)_{\mu}(x,\mu_1(x;r))=0,\qquad (Z_r)_{\mu\mu}(x,\mu_1(x;r))>0.
\end{equation*}

At the zero point $(x,\mu_1(x;r))$, similarly to estimates 
\eqref{eq:decay estimate} and \eqref{eq:first decay estimate}, 
the part with $n\neq 0,1$ in \eqref{eq:vertical-first} is 
$O(e^{-c/r})$ times the sum of the absolute values of the terms with $n=0,1$.
Specifically, let
$$A_{x,r}=e^{-\frac{\mu_1^2}{4r}}\theta(x;t)>0,\quad B_{x,r}=e^{-\frac{(\mu_1-1)^2}{4r}}\theta(x-a;t)>0,$$
then it follows from \eqref{eq:vertical-first} that
 \begin{equation}\label{eq:n=0,1 derivative}
       \begin{aligned}
            0&=(Z_r)_{\mu} (x,\mu_1(x;r))\\
            &=
        -\frac1{2r}\left( \mu_1 A_{x,r}
        +(\mu_1-1) B_{x,r} + (A_{x,r}+B_{x,r}) O(e^{-c/r})\right).
       \end{aligned}
 \end{equation}
 Since $\mu_1(x;r)\in (\frac{1}{4},\frac{3}{4})$, as $r\to 0^+$, we see from \eqref{eq:n=0,1 derivative}  that
 \begin{equation}
       \frac{\mu_1A_{x,r}}{(1-\mu_1)B_{x,r}}
       =1+O(e^{-c/r}).
\end{equation}
Taking the logarithm of both sides gives
 \begin{equation}\label{eq:mu1 value estimate}
       \frac{2\mu_1(x;r)-1}{4r}
  =\log\frac{\theta(x;t)}{\theta(x-a;t)}
   +\log\frac{\mu_1(x;r)}{1-\mu_1(x;r)}+O(e^{-c/r})=O(1),  
 \end{equation}
 which means that the unique zero curve satisfies $\mu_1(x;r)=\frac{1}{2}+O(r)$ in 
 $\left(1/4, 3/4 \right)$.
The proof is complete.
 \end{proof}
 
 We now prove Proposition \ref{prop:cusp-estimate} by studying the
  critical points on the two zero curves $(x,\mu_0(x;r))$ and $(x,\mu_1(x;r))$.

  \begin{proof}[\bf Proof of Proposition  \ref{prop:cusp-estimate}]
     According to Lemma \ref{le:partial Zr},    
  to study the critical points of $p_\tau(\cdot,t)$, we  only need to consider
  the critical points $(x,\mu)$ of $Z_r(x,\mu)$ on the zero curves 
  $(x,\mu_0(x;r))$ and $(x,\mu_1(x;r))$.

  {\bf Case 1: }
  On the zero curve $(x,\mu_0(x;r))$, we will prove that when $r>0$ is sufficiently small,
  $Z_r(x,\mu)$ has no nontrivial critical points.
  
  In fact,
by using \eqref{eq:scaled-vertical}, \eqref{eq:decay estimate} and the facts that $\theta(x;t)$ 
is uniformly bounded in $x\in [0,1],t\in J$ and 
$\mu_0(x;r)=O(e^{-c/r})$,
  as $r\to 0^+$, there exists a constant $c_0>0$ such that for any $x\in [0,1]$,
  \begin{equation}\label{eq:Zr ck estimate}
         \begin{aligned}
  &|Z_r(x,\mu_0(x;r))-\theta(x;t)|\\
  \leq& \big(|e^{-\mu_0(x;r)^2/(4r)}-1|
  + e^{-\mu_0(x;r)^2/(4r)}O(e^{-c/r})\big)\theta(x;t)\\
  \leq& C  \frac{\mu_0(x;r)^2}{4r}  +  O(e^{-c/r})=O(e^{-c_0/r}).
  \end{aligned} 
  \end{equation}
  Therefore, we have $\|Z_r(\cdot,\mu_0(\cdot;r))-\theta(\cdot;t)\|_{C^0}=O(e^{-c_0/r})$.

Next, we claim that for any $k\geq 1$, 
 \begin{equation}\label{eq:lower-profile}
         \|Z_r(\cdot,\mu_0(\cdot;r))-\theta(\cdot;t)\|_{C^k}=O(e^{-c_k/r}),\quad r\to 0^+.
  \end{equation}

Recalling \eqref{eq:vertical-first}, we can write $(Z_r)_\mu=0$ as
\begin{equation}\label{eq:Rexmu}
 \mu+R_r(x,\mu)=0,\qquad
 R_r(x,\mu):=\sum_{n\ne0}(\mu-n)e^{-\frac{n^2-2n\mu}{4r}}
       \frac{\theta(x-an;t)}{\theta(x;t)}.
\end{equation}
Using \eqref{eq:decay estimate}, $n^2-2n\mu\geq n^2/2$ for $|\mu|\leq1/4$, and the fact that
 all derivatives of $\theta(\cdot;t)$
are uniformly bounded for $t\in J$, we have that, for every fixed $k$,
\[
 \|R_r\|_{C^k_{x,\mu}(\T\times[-1/4,1/4])}
 \leq C_k e^{-c_k/r},\qquad
 1+\partial_\mu R_r\geq\frac12
\]
for small $r>0$ and uniformly in $a\in[0,1/2]$ and $t\in J$.
Here each differentiation only introduces polynomial factors in $n$ and $r^{-1}$,
which are absorbed by the Gaussian decay.
Successively differentiating $\mu_0(x;r)+R_r(x,\mu_0(x;r))=0$
therefore gives $\|\mu_0\|_{C_x^k}\leq C_k e^{-c_k/r}$.
Write
\[
Z_r(x,\mu)
=e^{-\frac{\mu^2}{4r}}\bigl(\theta(x;t)+E_r(x,\mu)\bigr),
\quad
E_r(x,\mu):=\sum_{n\ne0}
e^{-\frac{n^2-2n\mu}{4r}}\theta(x-an;t).
\]
Similarly, we have
\[
 \|E_r\|_{C^k_{x,\mu}(\T\times[-1/4,1/4])}
\leq C_k e^{-c_k/r}.
\]
Together with $\|\mu_0\|_{C_x^k}\leq C_k e^{-c_k/r}$, the chain rule yields
\[
\|E_r(\cdot,\mu_0(\cdot;r))\|_{C^k}
+\|e^{-\mu_0^2/(4r)}-1\|_{C^k}\le C_k e^{-c_k/r},
\]
after decreasing $c_k>0$ to absorb the factor $r^{-1}$.
Applying these $C^k$ estimates to
\[
\begin{aligned}
Z_r(x,\mu_0(x;r))-\theta(x;t)
=\bigl(e^{-\frac{\mu_0(x;r)^2}{4r}}-1\bigr)\theta(x;t)
+e^{-\frac{\mu_0(x;r)^2}{4r}}E_r(x,\mu_0(x;r))
\end{aligned}
\] proves \eqref{eq:lower-profile}.

  Note from \eqref{eq:circle-signs} that
  $ \theta_{xx}(0;t)<0<\theta_{xx}(\frac{1}{2};t)$. Then by using \eqref{eq:lower-profile} with $k=2$,
  there exists $\varepsilon>0$ such that for sufficiently small $r$,
  \begin{equation}
         \left\{\begin{aligned}
             \frac{d^2}{dx^2} Z_r(x,\mu_0(x;r))<0,&\quad \forall |x|<\varepsilon,\\
             \frac{d^2}{dx^2} Z_r(x,\mu_0(x;r))>0, &\quad \forall |x-\frac{1}{2}|<\varepsilon.
         \end{aligned}\right. 
  \end{equation}
The identities \begin{equation}\label{eq:Zrxm}Z_r(-x,-\mu)=Z_r(x,\mu),\quad Z_r(x+1,\mu)=Z_r(x,\mu),\end{equation}
together with the uniqueness of $\mu_0(x;r)$, give
$\mu_0(-x;r)=-\mu_0(x;r)$ and $\mu_0(x+1;r)=\mu_0(x;r)$.
Thus $\mu_0(0;r)=\mu_0(1/2;r)=0$, and $x\mapsto Z_r(x,\mu_0(x;r))$
is even and 1-periodic. Its derivative vanishes at $0$ and $1/2$, so only the trivial points $(x,\mu_0(x;r))=(0,0),(\frac{1}{2},0)$ satisfy
  $\frac{d}{dx} Z_r(x,\mu_0(x;r))=0$ in $x\in B_{\varepsilon}(0)\cup B_{\varepsilon}(\frac{1}{2})$.

  On the compact set $(x,t)\in (\R/\Z)\setminus \left(B_{\varepsilon}(0)\cup B_{\varepsilon}(\frac{1}{2})\right) \times J$, it follows from \eqref{eq:circle-signs} that
  the minimum $\min_{x,t}|\theta_x(x;t)|$ $=c_J>0$. Therefore, 
  based on the uniform convergence \eqref{eq:lower-profile} with $k=1$,
  when $r$ is sufficiently small, we obtain
  \begin{equation*}
         \left|\frac{d}{dx} Z_r(x,\mu_0(x;r))\right|>\frac{c_J}{2}>0,\quad \forall x\in (\R/\Z)\setminus \left(B_{\varepsilon}(0)\cup B_{\varepsilon}(\frac{1}{2})\right).
  \end{equation*}
  In summary, when $r>0$ is sufficiently small, there are no nontrivial critical points $(x,\mu_0(x;r))$
  satisfying $$0=\frac{d}{dx} Z_r(x,\mu_0(x;r))=\partial_x Z_r(x,\mu_0(x;r)),$$
  which means there are no nontrivial critical points on the zero curve
  $(x,\mu_0(x;r))$.

  {\bf Case 2: }
  On the zero curve $(x,\mu_1(x;r))$,
  for any fixed $a\in [0,\frac{1}{2}],t\in J$, let 
  \[
   h_r(x):=\frac{\mu_1(x;r)-1/2}{2r}.
  \]
  Since $\mu_1(x;r)=\frac{1}{2}+O(r)$, we have that $|h_r(x)|$ is bounded. 
  We make the change of variables \(\mu \mapsto h\) and define
  \begin{equation}
  \begin{aligned}
         \widehat Z_r(x,h)
         :=e^{\frac{1}{16r}}Z_r\left(x,\frac{1}{2}+2rh;a,t\right).
         \end{aligned}
  \end{equation}
  This \(\widehat{Z}_r(x,h)\) eliminates the singular term \(e^{-1/(16r)}\).

  Note that as $r\to 0^+$,
   the $n=0,1$ terms form the main part of $\widehat Z_r(x,h)$ 
  and that $n(n-1)\geq2$ for $n\notin\{0,1\}$.
  Since the Gaussian decay term controls the polynomial factors arising from differentiation,
  when $h$ is bounded,
   expression \eqref{eq:scaled-vertical} gives
     \begin{equation}\label{eq:upper-two-rows}
     \begin{aligned}
     \widehat Z_r(x,h)
     &=\sum_{n\in\Z}e^{-n(n-1)/(4r)}
          e^{(n-1/2)h-rh^2}\theta(x-an;t)\\
     &=e^{-rh^2}\bigl(e^{-h/2}\theta(x;t)+e^{h/2}\theta(x-a;t)\bigr)
       +O_{C_{x,h}^5}(e^{-c/r}),
     \end{aligned}
    \end{equation}
  where $O_{C^k}(\cdot)$ denotes a function whose
   $C^k$ norm is bounded by $C\cdot$.
 
  Since $(Z_r)_{\mu} (x,\mu_1(x;r))=0$ when $r>0$ is sufficiently small,
  we obtain 
  $\partial_h \widehat Z_r(x,h_r(x))$ $=0$. Define
$$  \widehat Z_0(x,h):=\lim_{r\to 0^+}  \widehat Z_r(x,h)=e^{-h/2}\theta(x;t)+e^{h/2}\theta(x-a;t).$$
Then 
  \begin{equation*}
         0=\partial_h\widehat Z_0(x,h_0(x))\quad \Longleftrightarrow\quad
   h_0(x)=\log(\theta(x;t)/\theta(x-a;t)).
  \end{equation*}
  At this point,
  \begin{equation*}
         \partial_h^2\widehat Z_0(x,h_0(x))
   =\frac{\sqrt{\theta(x;t)\theta(x-a;t)}}{2}\ge\frac12\min_{\substack{x\in[0,1]\\t\in J}}\theta(x;t)
 =:A>0.
  \end{equation*}
Then
by using the estimate \eqref{eq:upper-two-rows}, $\partial_h^2\widehat Z_0(x,h_0(x))>0$ and the 
implicit function theorem, we can obtain
\begin{equation}\label{eq:h estimate}
     h_r(x)=h_0(x)+O_{C_x^4}(r).  
\end{equation}

Indeed, to justify this $C_x^4$ estimate,
we set $F_r=\partial_h\widehat Z_r$ and choose $M>0$ such that
$|h_r(x)|,|h_0(x)|\leq M$ uniformly in $x$, $t\in J$, and small $r$.
By \eqref{eq:upper-two-rows}, on $|h|\le M$,
\begin{equation}\label{eq:FrF0}
 \|F_r-F_0\|_{C_{x,h}^4}
 \le C\bigl(r+e^{-c/r}\bigr)\le Cr.
\end{equation}
Moreover,
\[
 \begin{aligned}
 \partial_hF_0(x,h)
 &=\frac14\left(e^{-h/2}\theta(x;t)
                +e^{h/2}\theta(x-a;t)\right)\\
 &\ge\frac{\sqrt{\theta(x;t)\theta(x-a;t)}}{2}
 \ge A>0.
 \end{aligned}
\]
Thus $\partial_hF_r\ge A/2$ on $|h|\le M$ for small $r$.
Then the implicit function theorem and successive differentiation of
$F_r(x,h_r(x))=0$ give
\begin{equation}\label{eq:hrc4}
 \|h_r\|_{C_x^4}\le C.
\end{equation}
Indeed, at each order $1\le j\le4$, the highest derivative
$h_r^{(j)}$ has coefficient $\partial_hF_r(x,h_r)$, while all
remaining terms involve only derivatives of $F_r$ of order at most
$j$ and derivatives of $h_r$ of order less than $j$.

Since
$F_r(x,h_r)=F_0(x,h_0)=0$, we have 
\begin{align*}-(F_r(x,h_r)-F_0(x, h_r))&=F_0(x,h_r)-F_0(x,h_0)\\
&=(h_r-h_0)\int_0^1
   \partial_hF_0\bigl(x,(1-s)h_0(x)+sh_r(x)\bigr)\,ds\end{align*}
and so
\[
h_r(x)-h_0(x)
 =
 -\frac{(F_r-F_0)(x,h_r(x))}
 {\int_0^1
   \partial_hF_0\bigl(x,(1-s)h_0(x)+sh_r(x)\bigr)\,ds}.
\]
The numerator is $O_{C_x^4}(r)$ by the chain rule and \eqref{eq:FrF0}-\eqref{eq:hrc4}. The denominator is at least $A$, and it
and its reciprocal have uniformly bounded $C_x^4$ norms, since
$F_0$ is smooth and $h_0,h_r$ are uniformly bounded in $C_x^4$.
Consequently,
\[
\|h_r-h_0\|_{C_x^4}\le Cr,
\]
which proves \eqref{eq:h estimate}.

Substituting \eqref{eq:h estimate} into \eqref{eq:upper-two-rows}
and taking the logarithm, we obtain that 
\begin{equation}\label{eq:widehat Z estimate}
      \log\widehat Z_r(x,h_r(x))
 =\log\left(2\sqrt{\theta(x;t)\theta(x-a;t)}\right)
   +O_{C_x^4}(r).  
\end{equation}

Let $w=x-a/2$ and define
\begin{equation}\label{eq:Lr def}
       L_r(w;a,t):=
       \log\widehat Z_r
       \left(w+\frac a2,h_r\Big(w+\frac a2\Big)\right)=\frac1{16r}+\log Z_r(x,\mu_1(x;r)).
\end{equation}
       
From estimate \eqref{eq:widehat Z estimate}, we have 
\begin{equation}\label{eq:Lr}
       L_r(w;a,t)
 = \log\left(2\sqrt{\theta(w+a/2;t)\theta(w-a/2;t)}\right)+O_{C_w^4}(r).
\end{equation}
Note from $\tau=a+bi$, the periodicity and evenness of $p_{\tau}(\cdot,t)$ and \eqref{eq:scaled-vertical} that
$$Z_r(a-x,1-\mu)=Z_r(x,\mu).$$ Since $\mu_1(x;r)$ is the 
unique zero of $(Z_r)_\mu(x,\mu)=0$ in $[\frac{1}{4},\frac{3}{4}]$, we have 
$$\mu_1(a-x;r)=1-\mu_1(x;r).$$
Thus, based on the definition of $L_r(w;a,t)$,
we get
$$L_r(-w;a,t)=L_r(w;a,t),\qquad L_r(w+1;a,t)=L_r(w;a,t).$$
Then we obtain two fixed critical points
\begin{equation}\label{eq:Lr first derivative}
      L_r'(0;a,t)=L_r'\left(\frac{1}{2};a,t\right)=0,
 \qquad \forall a\in [0,\frac{1}{2}],\;t\in J,\;\text{small}~r>0.  
\end{equation}
Here, $L_r'$ is the derivative of $L_r$ with respect to $w$.
The same symmetry and periodicity give
$\mu_1(a/2;r)=\mu_1((1+a)/2;r)=1/2$.
Hence the points $w=0,1/2\pmod1$ on the upper zero curve correspond exactly to
the two half-period points $\tau/2$ and $(1+\tau)/2$.

We prove this proposition by contradiction.
 Suppose
  there exist
  \[
   b_j\longrightarrow\infty,\quad
   a_j\in[0,1/2],\quad t_j\in J,
  \]
  and non-half-period critical points $z_j$ of
  $p_{a_j+b_j i}(z_j,t_j)$ satisfying
  \begin{equation}\label{eq:fan zheng}
       \nabla_z p_{a_j+ib_j}(z_j,t_j)=0,\qquad
       D^2p_{a_j+ib_j}(z_j,t_j)\not>0.
  \end{equation}
Since
  $r_j=\frac{t_j}{b_j^2}\to 0$, Case 1 implies that $z_j$ lies on the upper zero curve, then by \eqref{eq:3-4},
  $z_j$ can not be a nondegenerate local maximum point.

Since $z_j=x_j+y_j i$ is a nontrivial critical point, 
we have $w_j=x_j-a_j/2\notin\{0,1/2\}\pmod1$.
The formula \eqref{eq:3-4} gives 
$ (Z_{r_j})_{\mu\mu}(x,\mu_1(x;r_j))>0$.
Based on the Schur complement formula
\eqref{eq:schur} and the definition \eqref{eq:Lr def} of $L_r$,
 at the critical point $z_j=x_j+y_j i$,
we have  
\begin{align*}
\det D^2Z_{r_j} (x_j,\mu_1(x_j;r_j))&=(Z_{r_j})_{\mu\mu} (x_j,\mu_1(x_j;r_j))\frac{d^2}{dx^2}Z_{r_j} (x_j,\mu_1(x_j;r_j))\\
&=
Z_{r_j} (x_j,\mu_1) (Z_{r_j})_{\mu\mu} (x_j,\mu_1)L_{r_j}''(w_j),
\end{align*}
Thus, at the critical point  $z_j=x_j+y_j i$,
\[
       D^2Z_{r_j} (x_j,\mu_1(x_j;r_j))>0\iff L_{r_j}''(w_j)>0.
\]
By the contradiction assumption \eqref{eq:fan zheng} and expression \eqref{eq:scaled-vertical},
we have 
\begin{equation}\label{eq:Lr critical}
      L_{r_j}'(w_j;a_j,t_j)=0,\qquad
   L_{r_j}''(w_j;a_j,t_j)\leq0. 
\end{equation}

By evenness and periodicity, we may assume $w_j\in(0,1/2)$. Since $w_j,a_j,t_j$ are bounded,
  up to a subsequence, 
  \[
   (w_j,a_j,t_j)\longrightarrow(w_0,a_0,t_0).
  \]

  By \eqref{eq:Lr} and the continuity of $L_r$ with respect to its parameters, we have
\begin{equation}\label{eq:Lr convergence}
     L_{r_j}(\cdot;a_j,t_j)
 \rightarrow L_0(\cdot;a_0,t_0)=\log\left(2\sqrt{\theta(\cdot+a_0/2;t_0)\theta(\cdot-a_0/2;t_0)}\right)
\end{equation}
 uniformly in $C^4(\T)$, and
\begin{equation*}
     L_0'(w_0;a_0,t_0)=0.
\end{equation*}
 
We now begin to derive a contradiction.
  There are three possibilities for $w_0$.

  If $w_0\notin\{0,1/2\}$, we will prove in Lemma \ref{lem:Lr inequalities} that
  \begin{equation}\label{eq:claim}
         L_0''(w_0;a_0,t_0)>0.
  \end{equation}
  Then by 
  \eqref{eq:Lr convergence}, this contradicts 
  \eqref{eq:Lr critical}.

If $w_0\in\{0,1/2\}$ and
$L_0''(w_0;a_0,t_0)\ne0$, then by
\eqref{eq:Lr convergence}, there exists $\delta>0$ such that 
$L_{r_j}''(\cdot;a_j,t_j)$ has a fixed nonzero sign
in a neighborhood $B_{\delta}(w_0)$ when $j$ is sufficiently large.
Hence $L_{r_j}'(\cdot;a_j,t_j)$ is strictly monotone
in $B_{\delta}(w_0)$ and has at most one zero.
However, $w_0\in\{0,1/2\}$ and \eqref{eq:Lr first derivative} implies that $w_0$ is a zero point of $L_{r_j}'(\cdot;a_j,t_j)$.
 This contradicts \eqref{eq:Lr critical} and $w_j\neq w_0$.
 
 If $w_0\in\{0,1/2\}$ with $L_0''(w_0;a_0,t_0)=0$, 
we will also prove in Lemma \ref{lem:Lr inequalities} that 
 $L_0''(0;a_0,t_0)<0$ and 
 \begin{equation}\label{eq:inequality claim}
        L_0''(1/2;a_0,t_0)=0
        \quad\Longrightarrow\quad
        L_0^{(4)}(1/2;a_0,t_0)>0.
       \end{equation}
 Since $L_0''(0;a_0,t_0)<0$, it follows that 
$w_0=1/2$ and $L_0''(1/2;a_0,t_0)=0$.
 Using \eqref{eq:inequality claim} and \eqref{eq:Lr convergence} again,
 there exist $\delta,c_0>0$ such that, for sufficiently large $j$,
 \begin{equation}\label{eq:Lr fourth}
        L_{r_j}^{(4)}(w;a_j,t_j)\geq c_0,
  \quad  \forall w\in B_{\delta} (w_0)=B_{\delta} (1/2).
 \end{equation}

 Since $L_{r_j}$ is even and 1-periodic, we have $L_{r_j}'(w_0;a_j,t_j)=L_{r_j}'''(w_0;a_j,t_j)=0$. Then by \eqref{eq:Lr fourth}, 
a direct calculation using calculus gives, for any $0<v<\delta$,
\begin{align}\label{eq:even-fourth}
&vL_{r_j}''(w_0+v;a_j,t_j)-L_{r_j}'(w_0+v;a_j,t_j)\notag\\
=&\int_0^v s\int_0^s L_{r_j}^{(4)}(w_0+r;a_j,t_j)\,dr\,ds>0.
\end{align}
Note that as $j\to\infty$, $w_j-w_0\to 0$. 
When $w_j-w_0\in (0,\delta)$,
using 
\eqref{eq:Lr critical} and \eqref{eq:even-fourth} and choosing $v=w_j-w_0$,
we have 
\begin{equation*}
       L_{r_j}''(w_j;a_j,t_j)>0.
\end{equation*}
This contradicts \eqref{eq:Lr critical}.
When $w_j-w_0\in (-\delta,0)$, based on the evenness of $L_r(\cdot)$,
we also obtain $ L_{r_j}''(w_j;a_j,t_j)>0$, again a contradiction.

Therefore, all possible limits have been excluded.
This completes the proof.
\end{proof}

  At the end of this subsection, we need to prove three 
  inequalities that are used in the proof of Proposition \ref{prop:cusp-estimate}.

\begin{lemma}\label{lem:Lr inequalities}
       Let $L_0$ be defined in \eqref{eq:Lr convergence}. For any $a\in [0,\frac{1}{2}]$
       and $t\in (0,\infty)$,  we have   $L_0''(0;a,t)<0$ and 
       \begin{equation}
              L_0'(w;a,t)=0,\quad w\notin\{0,\frac12\}\pmod1
                            \Longrightarrow L_0''(w;a,t)>0,
       \end{equation}
       \begin{equation}
              L_0''(1/2;a,t)=0
              \quad\Longrightarrow\quad
              L_0^{(4)}(1/2;a,t)>0.
             \end{equation}
\end{lemma}

\begin{proof}
First, we recall the Watson identity
 \begin{equation}
      \begin{aligned}\label{eq:addition}
  &\theta(w+a/2;t)\theta(w-a/2;t)\\
= &\frac12\left[
\theta(a/2;t/2)\theta(w;t/2)+
 \theta((1-a)/2;t/2)\theta(w-1/2;t/2)
 \right].
 \end{aligned} 
\end{equation}
This identity can be verified directly using 
       \cite[\S 20.7, (20.7.14)]{DLMF}.
 Let $q=e^{-2\pi^2t}\in(0,1)$. Define the auxiliary function, for any $-1\le c\le1$,
 \begin{equation}\label{eq:G}
 \begin{aligned}
  G(c):=\prod_{n=1}^{\infty}(1-q^{2n})
  \Bigg[
  &\theta(a/2;t/2)
  \prod_{n=1}^{\infty}(1+2q^{2n-1}c+q^{4n-2})\\
  &+
  \theta((1-a)/2;t/2)
  \prod_{n=1}^{\infty}(1-2q^{2n-1}c+q^{4n-2})
  \Bigg].
 \end{aligned}
 \end{equation}
 Then by \eqref{eq:addition} and \eqref{eq:product},  $L_0$ can be represented by 
 \begin{equation}\label{eq:representation}
 L_0(w;a,t)=\log\left(2\sqrt{\theta(w+a/2;t)\theta(w-a/2;t)}\right)
       =\frac12\log\bigl(2G(\cos(2\pi w))\bigr).
      \end{equation}
Therefore, studying the derivative of $L_0$ is 
reduced to studying the derivative of $G$.

For any $c\in[-1,1]$, we have 
$ 1\pm2q^{2n-1}c+q^{4n-2}
 \ge(1-q^{2n-1})^2>0$. Consequently,
 \begin{equation*}
     G(c)>0,\qquad\forall c\in [-1,1].   
 \end{equation*}
 Differentiating \eqref{eq:G} twice and then using 
 $ 1\pm2q^{2n-1}c+q^{4n-2}>0$ again, we also have 
 \begin{equation}\label{eq:G two derivative}
        G''(c)>0,\qquad\forall c\in [-1,1].   
    \end{equation}

Since $0\leq a\leq \frac{1}{2}$, using \eqref{eq:product}, we have 
$ \theta(a/2;t/2)\ge\theta((1-a)/2;t/2)>0$.
Differentiating \eqref{eq:G} at $c=0$ gives
\[
\begin{aligned}
G'(0)={}&\bigl[\theta(a/2;t/2)-\theta((1-a)/2;t/2)\bigr]\\
&\quad\times
\prod_{n=1}^{\infty}(1-q^{2n})(1+q^{4n-2})
\sum_{n=1}^{\infty}\frac{2q^{2n-1}}{1+q^{4n-2}}\geq0.
\end{aligned}
\]
Since $G''>0$ on $[0,1]$, we obtain
\begin{equation}\label{eq:G derivative}
       G'(1)>0. 
\end{equation}
 
Differentiating \eqref{eq:representation} gives
\begin{equation}\label{eq:L-prime}
 L_0'(w;a,t)
 =-\pi\sin(2\pi w)
  \frac{G'(\cos(2\pi w))}{G(\cos(2\pi w))}.
\end{equation}
Therefore, differentiating \eqref{eq:L-prime} and using 
\eqref{eq:G derivative},
 we get the first desired inequality
\[
 \quad L_0''(0;a,t)
 =-2\pi^2\frac{G'(1)}{G(1)}<0.
\]

Suppose $L_0'(w;a,t)=0$ and
$w\notin\{0,\frac12\}\pmod1$. Since $\sin(2\pi w)\ne0$,
\eqref{eq:L-prime} gives $G'(\cos(2\pi w))=0$.
Therefore,
differentiating \eqref{eq:L-prime} and using \eqref{eq:G two derivative},
we obtain the second desired inequality 
\[
\quad L_0''(w;a,t)
 =2\pi^2\sin^2(2\pi w)
  \frac{G''(\cos(2\pi w))}{G(\cos(2\pi w))}>0.
\]

Suppose $L_0''(1/2;a,t)=0$. Since $L_0'(1/2;a,t)=0$,
differentiating \eqref{eq:L-prime} yields
\[
 0=L_0''(1/2;a,t)
  =2\pi^2\frac{G'(-1)}{G(-1)}
 \Longrightarrow G'(-1)=0.
\]

As $v\to0$,
 $\cos\bigl(2\pi(1/2+v)\bigr)
 =-1+2\pi^2v^2+O(v^4).$
Thus, using $G'(-1)=0$,
Taylor's formula  gives
\begin{align*}
 G\!\left(\cos\bigl(2\pi(1/2+v)\bigr)\right)
 &=G(-1)+\frac12G''(-1)(2\pi^2v^2)^2+O(v^6).
\end{align*}
Substituting this back into formula \eqref{eq:representation}, we have 
\begin{align*}
 L_0(1/2+v;a,t)
 &=\frac12\log\bigl(2G(-1)\bigr)
    +\pi^4\frac{G''(-1)}{G(-1)}v^4+O(v^6).
\end{align*}
Comparing the coefficients of $v^4$ and using \eqref{eq:G two derivative}, we obtain
$$L_0^{(4)}(1/2;a,t)
 =24\pi^4\frac{G''(-1)}{G(-1)}>0.$$
The proof is complete.
\end{proof}

\subsection{Nondegenerate minimality when \(t\) is sufficiently small}\label{sec:uniform-small}

In this subsection, we study the case of sufficiently small $t$ in the
 following Proposition \ref{prop:uniform-small}.

\begin{proposition}\label{prop:uniform-small}
       There exists a small $t_*>0$ such that, for any 
       $0<t\leq t_*$ and
       $\tau=a+bi\in\mathcal F$ defined in \eqref{eq:reduced-region}, every nontrivial critical
       point of $p_\tau(\cdot,t)$ is a nondegenerate local minimum point.
 \end{proposition}

 The proof of Proposition \ref{prop:uniform-small} is similar to the proof 
 of Proposition \ref{prop:cusp-estimate}, but uses a more accurate estimate
 since the time $t$ can approach the singular value $0$.

 Using expression \eqref{eq:scaled-vertical}, we can study the critical
  points of $Z_r(x,\mu;a,t)$ in order to study the critical points of
   $p_{\tau}(x+yi,t)$, where $r=\frac{t}{b^2}$ and $\mu=\frac yb$.
   It is worth noting that in this subsection, 
   we regard \(a,b\) as fixed parameters,
   \(t\) as a varying parameter,
   and \(x,y\) as variables, which differs slightly from Section \ref{subsec:reduction}.

Note that for $a+bi\in \mathcal F$, as $t\to 0^+$, we have $r=\frac t{b^2}\le\frac43t \to 0^+$.
Thus \(t\to0^+\) implies \(r\to0^+\) uniformly on \(\mathcal F\), 
and we use the limit \(t\to0^+\) in this subsection.
We first prove that $\partial_{\mu}  Z_r(x,\mu)=0$ also 
has two zero curves.

\begin{lemma}\label{lem:zero curve Zr}
       There is $t_0>0$ sufficiently small such that for any $t\in (0, t_0)$, there are exactly two 
        smooth functions $\mu_j(\cdot;r,t): \mathbb{R}/\mathbb Z\to [-\frac{1}{4},\frac{3}{4}]$, $j=0,1$,
       such that  $\partial_{\mu}  Z_r(x,\mu)=0$
       for $\mu\in [-\frac{1}{4},\frac{3}{4}]$ if and only if $\mu=\mu_j(x;r,t)$ for some $j=0,1$.
       
     Furthermore,
        \begin{equation}
            \partial_{\mu}^2  Z_r(x,\mu_0(x;r,t)) <0, \quad  \mu_0(x;r,t)=O(e^{-c/r}),
        \end{equation}
        and
        \begin{equation}\label{eq:3-4-2}
                \partial_{\mu}^2  Z_r(x,\mu_1(x;r,t))>0, \quad \mu_1(x;r,t)\in \left(\frac{1}{4},\frac{3}{4} \right).     
        \end{equation}
 Here, we omit the fixed parameter $a$.
 \end{lemma}

 \begin{proof}
       Based on the periodicity, similar to Section \ref{subsec:reduction},
 we only need to consider $x\in [0,1]$ and $\mu\in \left[-\frac{1}{4},\frac{3}{4}\right]$.
 We also divide $\mu$ into two cases as in Lemma \ref{le:partial Zr}.

 {\bf Case 1: $|\mu|\leq \frac{1}{4}$}. The $n=0$ term is the main term, and for any $n\neq 0$,
 \begin{equation}\label{eq:3-38}
        \begin{aligned}
          \frac{e^{- \frac{(\mu-n)^2}{4r}}\theta(x-an;t)}{ e^{-\frac{\mu^2}{4r}}\theta(x;t)}&\leq
        e^{-\frac{(\mu-n)^2-\mu^2}{4r}}
        \frac{\theta(x-an;t)}{\theta(x;t)}\\
        &\le Ce^{-\frac{n^2}{8r}+\frac{1}{16t}}.    
        \end{aligned}
  \end{equation}
  The proof of this inequality is similar to that of \eqref{eq:decay estimate}
and uses \begin{equation}\label{eq:thetatheta}\frac{\theta(x-an;t)}{\theta(x;t)}\leq \frac{\theta(0;t)}
  {\theta(1/2;t)}
  \le Ce^{\frac{1}{16t}},\end{equation} which is the only difference between this Lemma and Lemma 
  \ref{le:partial Zr}. Here the first inequality of \eqref{eq:thetatheta} follows from \eqref{eq:product} and the second inequality follows from \eqref{eq2-9}.
Note from $r\leq4t/3$ that for $n\ne0$,
\[
 -\frac{n^2}{8r}+\frac1{16t}
 \leq-\frac{n^2}{8r}+\frac1{12r}
 \leq-\frac{n^2}{24r},
\]
so the quantity in \eqref{eq:3-38} is exponentially small uniformly on $\mathcal F$.
  Therefore, using \eqref{eq:3-38}, $r\le\frac43t$, and arguing similarly to the process in
  \eqref{eq:first decay estimate}, we also have 
  \begin{equation}
         \begin{aligned}
               &\frac{(Z_r)_\mu(x,\mu)}{e^{-\frac{\mu^2}{4r}}\theta(x;t)} 
  =-\frac{\mu}{2r}+O(r^{-1}e^{-c/r}), \\
  &\frac{(Z_r)_{\mu\mu}(x,\mu)}{e^{-\frac{\mu^2}{4r}}\theta(x;t)}
  =\frac{\mu^2}{4r^2}-\frac1{2r}
  +O(r^{-2}e^{-c/r}). 
         \end{aligned} 
  \end{equation}
  Repeating the proof of Lemma \ref{le:partial Zr}, we obtain the unique
   $\mu_0(x;r,t)=O(e^{-c/r})$ satisfying $\partial_{\mu}  Z_r(x,\mu_0(x;r,t)) =0$ and $\partial_{\mu}^2  Z_r(x,\mu_0(x;r,t)) <0 $.

   {\bf Case 2: $\mu\in (\frac{1}{4},\frac{3}{4})$}.
   This case is also similar to that in Lemma \ref{le:partial Zr}, and we only need to use 
   the new inequalities
   \eqref{eq:thetatheta} and $r\le\frac43t$.
More explicitly, \eqref{eq:3-38} at $\mu=1/4$, and its reflected version at $\mu=3/4$, give
\[
\begin{aligned}
&\frac{(Z_r)_\mu(x,1/4)}{e^{-1/(64r)}\theta(x;t)}
=-\frac1{8r}+O(r^{-1}e^{-c/r}),\\
&\frac{(Z_r)_\mu(x,3/4)}{e^{-1/(64r)}\theta(x-a;t)}
=\frac1{8r}+O(r^{-1}e^{-c/r}).
\end{aligned}
\]
These have opposite signs uniformly on $\mathcal F$ for small $t$.
Together with \eqref{eq:derivative increasing}, we obtain the unique zero curve $\mu_1(x;r,t)$.
   It is worth noting that we cannot obtain 
   \(\mu_1(x;r,t)=\tfrac12+O(r)\) 
   because \(\log\frac{\theta(x;t)}{\theta(x-a;t)}\) is unbounded as \(t\to 0^+\).
 \end{proof}

 Now, we prove that on the zero curve $(x,\mu_0(x;r,t))$, $Z_r(x,\mu)$ 
 has no nontrivial critical points.

 \begin{lemma}\label{le:no nontrivial on mu0}
       There exists $t_1>0$ such that for any $0<t<t_1$,
       \[
              \frac{\dd}{\dd x}Z_r(x,\mu_0(x;r,t))=0
              \quad\Longleftrightarrow\quad
              x\in\left\{0,\frac12\right\}\pmod1.
              \]
Consequently, if $(x_0,\mu_0(x_0;r,t))$ is a critical point of $Z_r(x,\mu)$, 
then $z_0=x_0+b\mu_0(x_0)i$
       is a half-period point, i.e., it is not a nontrivial critical point.
 \end{lemma}

Here and below, we often omit the parameters $r,t$ in $\mu_0(x_0;r,t)$.

\begin{proof} First, we claim that for $t>0$ small, 
\begin{equation}\label{eq:normalized}
\frac{1}{\theta(x;t)}\frac{\dd}{\dd x}Z_r(x,\mu_0(x);a,t)
=\partial_x\log\theta(x;t)+O_{C_x^1}(e^{-c/t}).
\end{equation}

Indeed, differentiating gives
\[
 \left|\partial_x^k e^{-\frac{(x-m)^2}{4t}}\right|
 \leq C_k t^{-k}(1+|x-m|^k)e^{-\frac{(x-m)^2}{4t}},
 \quad 0<t\leq t_0.
\]
Since by the periodicity, we only need to consider $|x|\leq 1/2$, we have $(x-m)^2-x^2\geq m^2-|m|\geq0$ and then
\[
\begin{aligned}
 \sum_{m\in\mathbb Z}(1+|x-m|^k)e^{-\frac{(x-m)^2}{4t}}
 &\leq e^{-\frac{x^2}{4t}}
 \sum_{m\in\mathbb Z}
 \bigl(1+(|m|+1)^k\bigr)e^{-\frac{m^2-|m|}{4t_0}}\\
 &\leq C_k e^{-\frac{x^2}{4t}}.
\end{aligned}
\]
Thus, by \eqref{eq2-9}, we get
\begin{equation}\label{eq:patheta}
 |\partial_x^k\theta(x;t)|
 \leq \frac{C_k t^{-k}}{\sqrt{4\pi t}}e^{-\frac{x^2}{4t}}
 \leq C_k t^{-k}\theta(x;t),\quad 0<t\leq t_0,
\end{equation}
where $C_k$ depends only on $k$ and $t_0$. Consequently,
\[\bigg|\partial_x^k\bigg(\frac{\theta(x-an;t)}{\theta(x;t)}\bigg)\bigg|\leq C_k t^{-k}\frac{\theta(x-an;t)}{\theta(x;t)},\quad 0<t\leq t_0.\]
Recalling $R_r$ defined in \eqref{eq:Rexmu}, then
\eqref{eq:3-38} and $t^{-1}\leq4/(3r)$ give
\[
 \|R_r\|_{C^k_{x,\mu}(\T\times[-1/4,1/4])}
 \leq C_k e^{-c_k/r}.
\]
Implicit differentiation of $\mu_0+R_r(x,\mu_0)=0$ as before gives
$\|\mu_0\|_{C_x^2}\leq Ce^{-c/r}$. Then a similar argument as \eqref{eq:lower-profile} yields
\[
 Z_r(x,\mu_0(x);a,t)=\theta(x;t)(1+\eta_r(x)),\qquad
 \|\eta_r\|_{C_x^2}\leq Ce^{-c/r}.
\]
Since $\|\partial_x\log\theta\|_{C_x^1}\leq Ct^{-2}$,
differentiation gives
\[
 \frac{1}{\theta}\frac{d}{dx}Z_r(x,\mu_0(x);a,t)
 =(1+\eta_r)\partial_x\log\theta+\eta_r'.
\]
Using $r\leq4t/3$ and reducing $c>0$ if necessary, we obtain \eqref{eq:normalized}.

The same argument as \eqref{eq:Zrxm} implies that $x\mapsto Z_r(x,\mu_0(x))$
is even and 1-periodic, and $\mu_0(0)=\mu_0(1/2)=0$.
Thus we have two trivial critical points
\begin{equation}\label{eq:exact}
\frac{\dd}{\dd x}Z_r(0,0)=\frac{\dd}{\dd x}Z_r(1/2,0)=0\quad \textrm{with}\quad \mu_0(0)=\mu_0(1/2)=0.
\end{equation}
       
Suppose the lemma is false.
Then by the evenness and period-one periodicity of  $x\mapsto Z_r(x,\mu_0(x))$,
there exist 
$t_j\to 0^+$, $0<r_j\le4t_j/3$, $a_j\in[0,1/2]$, and 
non-half-period zeros $x_j\in(0,1/2)$ satisfying 
$$\frac{\dd}{\dd x}Z_{r_j}(x_j,\mu_0(x_j;r_j,t_j);a_j,t_j)=0.$$

Note that, for any $x\in (0,\beta_0]\subset (0,\frac{1}{2})$,
as $t\to 0^+$, the $n=0$ term $\frac1{\sqrt{4\pi t}} e^{-\frac{x^2}{4t}}$
is the main term of $ \theta(x;t)$ in the definition
\eqref{eq2-9}:
$$\theta(x;t)=\frac{1}{\sqrt{4\pi t}}e^{-\frac{x^2}{4t}}\bigg(1+\sum_{m\in\mathbb Z\setminus\{0\}}e^{\frac{2mx-m^2}{4t}}\bigg).$$
 Thus, by
using the Gaussian decay term to control the polynomial term,
 we obtain the uniform estimate
\begin{equation}\label{eq:log theta estimate}
       \partial_x\log\theta(x;t)=-\frac{x}{2t}+O_{C_x^1}(e^{-c/t}),\quad t\to 0^+.
\end{equation}

Therefore, if up to a subsequence, 
$x_j\to x_*\in (0,\frac{1}{2})$ as $j\to \infty$, then by using 
\eqref{eq:normalized} and \eqref{eq:log theta estimate}, 
we obtain
\[
0=\frac1{\theta(x_j;t_j)}
\frac{d}{dx}Z_{r_j}(x_j,\mu_{0}(x_j))=-\frac{x_j}{2t_j}+O(e^{-c/t_j})<0,
~\textmd{for $j$ large},
\]
a contradiction. 
Hence, up to a subsequence,
$x_j\to 0$ or $x_j\to\frac12.$ We will derive a contradiction
in each of these two cases.

{\bf Case 1: $x_j\to 0$.} 
For small $\delta>0$ and $x\in (0,\delta)$, similarly to the proof of 
\eqref{eq:log theta estimate}, we obtain the uniform estimate
\begin{equation}\label{eq:two derivative theta}
   \partial_x^2\log\theta(x;t)=-\frac1{2t}+O(e^{-c/t}),\quad~t\to 0^+.   
\end{equation}

When $t>0$ is small enough,
differentiating \eqref{eq:normalized} therefore gives
\[
\frac{\dd}{\dd x}\left[
\frac{1}{\theta(x;t)}\frac{\dd}{\dd x}Z_r(x,\mu_0(x);a,t)
\right]
=-\frac1{2t}+O(e^{-c/t})<0.
\]
Thus, based on \eqref{eq:exact},
  $x=0$ is the unique solution of $\frac{\dd}{\dd x}Z_r(x,\mu_0(x))=0$
 for $x\in [0,\delta)$. This contradicts $\frac{\dd}{\dd x}Z_{r_j}(x_j,\mu_0(x_j))=0$
when $j$ is sufficiently large.

{\bf Case 2: $x_j\to \frac12$.}
Note that near the point $\frac{1}{2}$,
as $t\to 0^+$, the $n=0,1$ terms
are the main terms of $ \theta(\cdot;t)$ in the definition
\eqref{eq2-9}:
$$\theta(x;t)=\frac{1}{\sqrt{4\pi t}}\bigg(e^{-\frac{x^2}{4t}}+e^{-\frac{(x-1)^2}{4t}}+\sum_{m\in\mathbb Z\setminus\{0,1\}}e^{-\frac{(x-m)^2}{4t}}\bigg).$$
Thus,
similar to the proof of \eqref{eq:log theta estimate},
when $|x|<\delta$ for a sufficiently small $\delta>0$, we have a
uniform estimate
\begin{equation}\label{eq:theta partial x}
       \partial_x \log\theta(1/2+x;t)
=-\frac x{2t}
+\frac1{4t}\tanh\frac x{4t}
+O_{C_x^1}(e^{-c/t}),\quad~t\to 0^+.
\end{equation}
      Using the fact that $x_j$ is a critical point and $x_j\to \frac12$, together with formulas \eqref{eq:normalized} and \eqref{eq:theta partial x}, we obtain
$$-\frac {x_j-\frac12}{2t_j}
+\frac1{4t_j}\tanh\frac {x_j-\frac12}{4t_j}\to 0$$
and so $(x_j-1/2)/t_j\to0$.
On the other hand, again by \eqref{eq:normalized} and \eqref{eq:theta partial x}, 
we get that for $|x|\leq t_j$ and sufficiently large $j$,
\begin{align*}
&\frac{\dd}{\dd x}\left[
\frac{1}{\theta(1/2+x;t_j)}\frac{\dd}{\dd x}
Z_{r_j}(1/2+x,\mu_{0}(1/2+x))\right]\\
=&\frac{1}{4t_j^2(e^{\frac{x}{4t_j}}+e^{-\frac{x}{4t_j}})^2}-\frac1{2t_j}
+O(e^{-c/t_j})>0.
\end{align*}
Thus, based on \eqref{eq:exact}, $x=\frac12$ is the unique solution of $\frac{\dd}{\dd x}Z_{r_j}(x,\mu_0(x);a_j, t_j)=0$
 for $x\in [\frac12-t_j,\frac12+t_j]$, 
which contradicts $\frac{\dd}{\dd x}
Z_{r_j}(x_j,\mu_{0}(x_j);a_j,t_j)=0$ and $x_j\in (\frac12-t_j, \frac12+t_j)$ for $j$ large. This completes the proof.
\end{proof}

Now, we begin the proof of Proposition \ref{prop:uniform-small}, which is 
the most delicate proof in this paper.
This proof is similar to the proof of Proposition \ref{prop:cusp-estimate} 
but uses some more accurate estimates.

\begin{proof}[\bf Proof of Proposition  \ref{prop:uniform-small}]
According to Lemma \ref{lem:zero curve Zr} and Lemma \ref{le:no nontrivial on mu0},  
when $t>0$ is sufficiently small, 
the nontrivial critical points of $Z_r(x,\mu)$
only exist on the zero curve  $(x,\mu_1(x;r,t))$.

Suppose this proposition is false.
Then there exist $t_j\to 0^+$, $\tau_j=a_j+ib_j\in\mathcal{F}$, and $r_j=t_j/b_j^2\leq \frac43 t_j$
such that 
\begin{equation}\label{eq:counter}
\nabla p_{\tau_j}(z_j,t_j)=0,\quad D^2p_{\tau_j}(z_j,t_j)\not>0,
\end{equation}
where $z_j=x_j+y_j i$ is not a half-period point. 
Here we use \eqref{eq:3-4-2}
to exclude the possibility that $z_j$ is a nondegenerate local maximum point.

By Lemma \ref{lem:zero curve Zr}, we have 
$\mu_1(x_j;r_j,t_j)\in (\frac{1}{4},\frac{3}{4})$.
We recall the function $L_r$ defined in \eqref{eq:Lr def}
\begin{equation}\label{eq:reduced}
 L_r(w;a,t)=\log Z_r\bigl(w+a/2,\mu_1(w+a/2);a,t\bigr)+\frac1{16r}.
\end{equation}
Let $w_j=x_j-a_j/2$. As in the proof of Proposition 
\ref{prop:cusp-estimate},
we again obtain the properties \eqref{eq:Lr first derivative} and
\eqref{eq:Lr critical}:
\begin{equation}\label{eq:Lr first derivative small}
      L_r'(0;a,t)=L_r'\left(\frac{1}{2};a,t\right)=0,
\end{equation}
\begin{equation}\label{eq:bad}
       L_{r_j}'(w_j;a_j,t_j)=0,\qquad
    L_{r_j}''(w_j;a_j,t_j)\leq0. 
 \end{equation}
 We show that \eqref{eq:bad} leads to a contradiction as $j\to \infty$ in four steps.

By the evenness and periodicity of $L_r$, we have 
$$L_{r_j}(-w;a_j,t_j)=L_{r_j}(w;a_j,t_j)=L_{r_j}(w+1;a_j,t_j).$$
Thus, without loss of generality, we can assume that
$1/2<w_j<1$.

{\bf Step 1: We prove $w_j\in (\frac{1}{2},\frac{4}{5})$ for $j$ large.}
We first give an estimate for $L_{r_j}''(w;a_j,t_j)$ as $j\to \infty$.
Let $\mu_{1,j}(x):=\mu_1(x;r_j,t_j)\in (\frac{1}{4},\frac{3}{4})$, with $a=a_j$,
and $w=x-\frac{a_j}{2}\in (\frac{1}{2},1)$.
Then, at the zero point $(x,\mu_{1,j}(x))$,
the \(n=0,1\) terms are the main part of $Z_{r_j} (x,\mu_{1,j}(x))$. Thus,
as $j\to \infty$, we  have an estimate similar to  
\eqref{eq:n=0,1 derivative}:
\begin{equation*}
       \begin{aligned}
            0&=(Z_{r_j})_{\mu} (x,\mu_{1,j}(x))\\
&= -\frac1{2r_j}\left( \mu_{1,j} A_{x,j}
            +(\mu_{1,j}-1) B_{x,j} + (A_{x,j}+B_{x,j}) O(e^{-c/t_{j}})\right),
       \end{aligned}
 \end{equation*}
where $$A_{x,j}=e^{- \frac{\mu_{1,j}(x)^2}{4r_j}}\theta(x;t_j)>0,\quad B_{x,j}=e^{- \frac{(\mu_{1,j}(x)-1)^2}{4r_j}}\theta(x-a_j;t_j) >0.$$
Letting $j\to \infty$, we also obtain the asymptotic equivalence
\begin{equation}\label{eq:main term compare}
       \frac{\mu_{1,j}(x)A_{x,j}}
       {(1-\mu_{1,j}(x))B_{x,j}}
       =1+O(e^{-c/t_j}),\quad~j\to \infty.
\end{equation}
From here, we claim that
\begin{equation}\label{eq:small-upper-derivative}
\begin{aligned}
L_{r_j}'(x-a_j/2;a_j,t_j)
={}&(1-\mu_{1,j}(x))\partial_x\log\theta(x;t_j)\\
 &+\mu_{1,j}(x)\partial_x\log\theta(x-a_j;t_j)
 +O(e^{-c/t_j}).
\end{aligned}
\end{equation}

Indeed, since $(Z_{r_j})_\mu(x,\mu_{1,j}(x))=0$, we have
\[
 L_{r_j}'(x-a_j/2;a_j,t_j)
 =\frac{(Z_{r_j})_x(x,\mu_{1,j}(x))}
        {Z_{r_j}(x,\mu_{1,j}(x))}.
\]
Note from \eqref{eq:patheta} that
$|\theta_x(x;t_j)|\le Ct_j^{-1}\theta(x;t_j)$.
Since the remainder terms with \(n\neq 0,1\) are 
\(O(e^{-c/t_j})\) times the \(n=0,1\) terms in the expression of $Z_{r_j} (x,\mu_{1,j}(x))$, i.e.,
\[ Z_{r_j}(x,\mu_{1,j}(x))
 =(A_{x,j}+B_{x,j})\bigl(1+O(e^{-c/t_j})\bigr),\]
we get
\[
 \begin{aligned}
 (Z_{r_j})_x(x,\mu_{1,j}(x))
 &=A_{x,j}\partial_x\log\theta(x;t_j)
   +B_{x,j}\partial_x\log\theta(x-a_j;t_j)\\
 &\quad +(A_{x,j}+B_{x,j})O(t_j^{-1}e^{-c/t_j}),
 \end{aligned}
\]
where the derivative is taken with $\mu$ fixed before setting
$\mu=\mu_{1,j}(x)$.
On the other hand, \eqref{eq:main term compare} and
$1/4<\mu_{1,j}(x)<3/4$ imply
\[
 \begin{aligned}
 \frac{B_{x,j}}{A_{x,j}+B_{x,j}}
 &=\mu_{1,j}(x)+O(e^{-c/t_j}),\\
 \frac{A_{x,j}}{A_{x,j}+B_{x,j}}
 &=1-\mu_{1,j}(x)+O(e^{-c/t_j}).
 \end{aligned}
\]
Substituting these identities into the quotient above proves
\eqref{eq:small-upper-derivative}, after absorbing the factor $t_j^{-1}$
into the exponential by decreasing $c>0$.

We will also use differentiated versions of these estimates.
For fixed $k,\ell\ge0$ with $k+\ell\le4$, uniformly for
$1/4\le\mu\le3/4$,
\begin{equation}\label{eq:scaled-tail-estimate}
 \left|(t_j\partial_x)^k(r_j\partial_\mu)^\ell
 \left[
 \frac{Z_{r_j}(x,\mu)}
 {e^{-\frac{\mu^2}{4r_j}}\theta(x;t_j)
  +e^{-\frac{(\mu-1)^2}{4r_j}}\theta(x-a_j;t_j)}-1
 \right]\right|\le Ce^{-c/r_j}.
\end{equation}
Indeed, by \eqref{eq:thetatheta},  the terms $n=-k$ and $n=1+k$, $k\ge1$, relative to
$n=0$ and $n=1$, respectively, are bounded by
$$\frac{e^{-\frac{(\mu+k)^2}{4r_j}}\theta(x+a_jk;t_j)+e^{-\frac{(\mu-1-k)^2}{4r_j}}\theta(x-a_j(k+1);t_j)}{e^{-\frac{\mu^2}{4r_j}}\theta(x;t_j)
  +e^{-\frac{(\mu-1)^2}{4r_j}}\theta(x-a_j;t_j)}\leq Ce^{-\frac{k^2+k/2}{4r_j}+\frac1{16t_j}}.$$
The inequality $r_j\le4t_j/3$ leaves a fixed positive Gaussian
decay rate in $1/r_j$. Also,
$|\partial_x^k\theta(x;t_j)|\le C_k t_j^{-k}\theta(x;t_j)$
by \eqref{eq:patheta}; hence the scaled derivatives above only introduce
polynomial factors, which are absorbed by decreasing $c$.
Thus, \eqref{eq:scaled-tail-estimate} holds.

By \eqref{eq:scaled-tail-estimate}, we claim that
\begin{equation}\label{eq:partialmu1}
 |\partial_x^k\mu_{1,j}(x)|\le C_k r_jt_j^{-k},\quad 1\le k\le 3.
\end{equation}
To prove \eqref{eq:partialmu1}, fix $x_0\in[0,1]$ and set
$\mu_*:=\mu_{1,j}(x_0)\in(1/4,3/4)$. Introduce the scaled variables
\[
 x=x_0+t_jX,\qquad \mu=\mu_*+r_jY,
 \qquad
 F_j(X,Y):=2r_j\frac{(Z_{r_j})_\mu(x,\mu)}{Z_{r_j}(x,\mu)}.
\]
We first estimate the partial derivatives of $F_j$.
Set
\[
 A=e^{-\frac{\mu^2}{4r_j}}\theta(x;t_j),\quad
 B=e^{-\frac{(\mu-1)^2}{4r_j}}\theta(x-a_j;t_j),\quad
 \varepsilon_j=\frac{Z_{r_j}(x,\mu)}{A+B}-1.
\]
By \eqref{eq:scaled-tail-estimate}, all $(X,Y)$ derivatives of
$\varepsilon_j$ of order at most four are $O(e^{-c/r_j})$.
Since $\partial_Y=r_j\partial_\mu$, direct differentiation gives
\[
 \begin{aligned}
 F_j(X,Y)
 &=2\partial_Y\log(A+B)+2\partial_Y\log(1+\varepsilon_j)\\
 &=-\mu+p_j(X,Y)+O_{C^3_{X,Y}}(e^{-c/r_j}),
 \quad\text{where}\quad p_j:=\frac{B}{A+B}.
 \end{aligned}
\]
To estimate the derivatives of $p_j$, write
\[
 p_j=\frac{e^{H_j}}{1+e^{H_j}},\quad\text{where}\quad
 H_j:=\log\frac BA
 =\frac{2\mu-1}{4r_j}
   +\log\frac{\theta(x-a_j;t_j)}{\theta(x;t_j)}.
\]
Clearly $\partial_YH_j=1/2$, all derivatives
involving $Y$ of total order at least two vanish, and \eqref{eq:patheta} gives
$|\partial_X^mH_j|\le C_m$ for $1\le m\le3$.
The function $u\mapsto e^u/(1+e^u)$ and its first three derivatives
are bounded on $\mathbb R$. The chain rule therefore gives
\[
 |\partial_X^k\partial_Y^\ell F_j(0,0)|\le C,
 \qquad 1\le k+\ell\le3,
\]
uniformly in $j$ and $x_0$.
Moreover, $F_j(0,0)=0$ implies
$p_j(0,0)=\mu_*+O(e^{-c/r_j})$.
Since $\partial_Yp_j=p_j(1-p_j)/2$, we obtain
\[
 \begin{aligned}
 (F_j)_Y(0,0)
 &=-r_j+\frac12p_j(0,0)(1-p_j(0,0))+O(e^{-c/r_j})\\
 &=-r_j+\frac12\mu_*(1-\mu_*)+O(e^{-c/r_j})
 \ge c_0>0
 \end{aligned}
\]
for large $j$, because $\mu_*(1-\mu_*)\ge3/16$.

Now define, for $X$ near $0$,
\[
 Y_j(X):=\frac{\mu_{1,j}(x_0+t_jX)-\mu_*}{r_j}.
\]
Then $Y_j(0)=0$ and $F_j(X,Y_j(X))=0$.
Differentiating this identity yields
\[
 \begin{aligned}
 (F_j)_Y Y_j'={}&-(F_j)_X,\\
 (F_j)_Y Y_j''={}&-\bigl((F_j)_{XX}
       +2(F_j)_{XY}Y_j'+(F_j)_{YY}(Y_j')^2\bigr),\\
 (F_j)_Y Y_j'''={}&-\bigl((F_j)_{XXX}
       +3(F_j)_{XXY}Y_j'+3(F_j)_{XYY}(Y_j')^2\\
 &\qquad +(F_j)_{YYY}(Y_j')^3
       +3\bigl((F_j)_{XY}+(F_j)_{YY}Y_j'\bigr)Y_j''\bigr).
 \end{aligned}
\]
All derivatives of $F_j$ are evaluated at $(X,Y_j(X))$.
Evaluating at $X=0$ and using the preceding bounds successively gives
$|Y_j'(0)|+|Y_j''(0)|+|Y_j'''(0)|\le C$.
Finally,
\[
 \partial_x^k\mu_{1,j}(x_0)
 =r_jt_j^{-k}Y_j^{(k)}(0),\qquad 1\le k\le3.
\]
Since $x_0$ was arbitrary, this proves \eqref{eq:partialmu1}.

For the second derivative of $L_r$,
by using the Schur complement formula \eqref{eq:schur}, we have 
\begin{equation}\label{eq:L_r schur}
       L_{r_j}''(x-a_j/2;a_j,t_j)
=(\log Z_{r_j})_{xx}-\left.\frac{(\log Z_{r_j})_{x\mu}^2}{(\log Z_{r_j})_{\mu\mu}}\right|_{(x,\mu_{1,j}(x))}.
\end{equation}
Similar to the proof of \eqref{eq:small-upper-derivative} and using \eqref{eq:scaled-tail-estimate} for $k+\ell\le2$, we have
\[
 \begin{aligned}
 (\log Z_{r_j})_{xx}
 ={}&(1-\mu_{1,j})\partial_x^2\log\theta(x;t_j)
    +\mu_{1,j}\partial_x^2\log\theta(x-a_j;t_j)\\
 &+\mu_{1,j}(1-\mu_{1,j})
 \bigl[\partial_x\log\theta(x;t_j)
       -\partial_x\log\theta(x-a_j;t_j)\bigr]^2
 +O(e^{-c/t_j}),\\
 (\log Z_{r_j})_{x\mu}
 ={}&-\frac{\mu_{1,j}(1-\mu_{1,j})}{2r_j}
 \bigl[\partial_x\log\theta(x;t_j)
       -\partial_x\log\theta(x-a_j;t_j)\bigr]\\
 &+O(r_j^{-1}e^{-c/t_j}),\\
 (\log Z_{r_j})_{\mu\mu}
 ={}&\frac{\mu_{1,j}(1-\mu_{1,j})-2r_j}{4r_j^2}
 +O(r_j^{-2}e^{-c/t_j}).
 \end{aligned}
\]
Since $\mu_{1,j}(1-\mu_{1,j})\ge3/16$, the last denominator in
\eqref{eq:L_r schur} is bounded below by $c/r_j^2$.
Then inserting these estimates into \eqref{eq:L_r schur} gives
\begin{equation}\label{eq:Lr'' estimate}
 \begin{aligned}
&L_{r_j}''(w;a_j,t_j)\\
=&(1-\mu_{1,j})\partial_x^2\log\theta(w+\frac{a_j}{2};t_j)+\mu_{1,j}\partial_x^2\log\theta(w-\frac{a_j}{2};t_j)\\
&-\frac{2r_j\mu_{1,j}(1-\mu_{1,j})}
{\mu_{1,j}(1-\mu_{1,j})-2r_j}\cdot
\left[
\partial_x\log\theta(w+\frac{a_j}{2};t_j)
-\partial_x\log\theta(w-\frac{a_j}{2};t_j)
\right]^2\\
&+O(e^{-c/t_j}),~j\to\infty.
\end{aligned}      
\end{equation}
This estimate for $w=x-\frac{a_j}{2}\in (\frac{1}{2},1)$ holds uniformly.

For any $w\in [\frac{4}{5},1)$,
since $a_j\in [0,\frac{1}{2}]$, we obtain 
$w\pm \frac{a_j}{2}-1\in [-\frac{9}{20},\frac{1}{4})$.
Therefore, the $n=1$ term is the main term of $\theta\left(w\pm \frac{a_j}{2};t_j\right)$.
Then, similar to the proof of 
\eqref{eq:two derivative theta}, 
when $w\in [\frac{4}{5},1)$, we obtain the uniform estimate
\begin{equation}\label{eq:variable derivative theta}
\partial_x^2\log\theta\left(w\pm \frac{a_j}{2};t_j\right)
=-\frac1{2t_j}+O(e^{-c/t_j}),~j\to\infty.
\end{equation}
However, substituting \eqref{eq:variable derivative theta} into 
 \eqref{eq:Lr'' estimate} and using $\mu_{1,j}\in(\frac{1}{4},\frac{3}{4})$,
for sufficiently large $j$, we have 
\[
L_{r_j}''(w;a_j,t_j)\le-\frac1{2t_j}+O(e^{-c/t_j})<0,\quad \forall w\in [4/5,1).
\]
Thus $w=1$ is the unique zero point of
$L_{r_j}'(w;a_j,t_j)$ 
on $w\in [\frac{4}{5},1]$.
 Hence $w_j\in (\frac{1}{2},\frac{4}{5})$
when $j$ is sufficiently large.

{\bf Step 2:}
For any $j\in\Z_+$, we define
$v_j:=\frac{w_j-1/2}{4t_j}>0$.
We need to give a delicate estimate of $v_j$,
and in this step we prove that 
$\frac{a_j}{8t_j}-v_j\to+\infty$
as $j\to \infty$.

Note from $w_j\in (\frac{1}{2},\frac{4}{5})$ that $w_j\pm \frac{a_j}{2}=\frac12+4t_j(v_j\pm \frac{a_j}{8t_j}) \in (\frac{1}{4},\frac{21}{20})$, so the $n=0,1$ terms are 
the main part of 
$\theta(\cdot;t_j)$ in the definition
\eqref{eq2-9}.  
Thus, by a similar proof as
\eqref{eq:theta partial x}, we obtain the uniform first and second derivative estimates
\begin{align}\label{eq:theta first partial x}
       &4t_j\partial_x \log\theta\left(w_j\pm \frac{a_j}{2};t_j\right)\notag\\
  =&\tanh \left(v_j\pm \frac{a_j}{8t_j} \right)-8t_j \left(v_j\pm \frac{a_j}{8t_j} \right)+O(e^{-c/t_j}),~j\to\infty ,
\end{align}
\begin{equation}\label{eq:theta second partial x}
 16t_j^2\partial_x^2\log\theta\left(w_j\pm \frac{a_j}{2};t_j\right)
=\operatorname{sech}^2 \left(v_j\pm \frac{a_j}{8t_j} \right)-8t_j+O(e^{-c/t_j}),~j\to \infty.  
\end{equation}

Note from \eqref{eq:patheta} that $
|\partial_x\log\theta(x;t)|\le\frac Ct.
$
Substituting \eqref{eq:theta second partial x} into
\eqref{eq:Lr'' estimate}
and using $\mu_{1}(x_j)\in (\frac{1}{4},\frac{3}{4})$
and \eqref{eq:bad}, 
there exists a constant $C>0$ such that  
\begin{equation}\label{eq:3-59}
\begin{aligned}
0\geq 16t_j^2L_{r_j}''(w_j)\ge{}&
(1-\mu_{1}(x_j))\operatorname{sech}^2
\left(\frac{a_j}{8t_j}+v_j\right)\\
&+\mu_{1}(x_j)\operatorname{sech}^2
\left(\frac{a_j}{8t_j}-v_j\right)-Ct_j,~j\to \infty.
\end{aligned}
\end{equation}
Letting $j\to \infty$, we obtain that 
\begin{equation}\label{eq:vj infty}
  \left|\frac{a_j}{8t_j}\pm v_j\right|\longrightarrow\infty.     
\end{equation}

Similarly, substituting \eqref{eq:theta first partial x} into  
\eqref{eq:small-upper-derivative} with $x=w_j+\frac{a_j}{2}$
and using \eqref{eq:bad}, we obtain
\begin{equation}\label{eq:vj infty estimate}
       \begin{aligned}
          0=&4t_j L_{r_j}'(w_j;a_j,t_j)\\
=&(1-\mu_1(x_j))\tanh\left(\frac{a_j}{8t_j}+v_j\right)
          -\mu_1(x_j)\tanh\left(\frac{a_j}{8t_j}-v_j\right)\\
        &-8t_jv_j+a_j(2\mu_1(x_j)-1)+O(v_je^{-c/t_j}),~j\to\infty.      
       \end{aligned}    
\end{equation}
Here, we need the remainder term $O(v_je^{-c/t_j})$ instead of $O(e^{-c/t_j})$ for later usage. The factor $v_j$ in $O(v_je^{-c/t_j})$ follows from an exact reflection.
Indeed, note from \eqref{eq:scaled-vertical} and uniqueness of $\mu_{1,j}(x)$ that
\[
 Z_{r_j}(1+a_j-x,1-\mu)=Z_{r_j}(x,\mu),\qquad
 \mu_{1,j}(1+a_j-x)=1-\mu_{1,j}(x).
\]
Thus, under $x=(1+a_j)/2+4t_jv$, the function
$\mu_{1,j}(x)-1/2$ is odd in $v$.
Both the function $v\mapsto 4t_jL_{r_j}'(1/2+4t_jv;a_j,t_j)$ and the main terms in
\eqref{eq:vj infty estimate}, viewed as functions of $v$, are odd, so their difference vanishes
at $v=0$. Estimates \eqref{eq:scaled-tail-estimate} and
\eqref{eq:partialmu1} with $k=1$ bound its $v$ derivative by
$Ce^{-c/t_j}$ uniformly on $[0,v_j]$. Thus,
integration from $0$ to $v_j$ gives $O(v_je^{-c/t_j})$.

   Since $a_j\geq 0,t_j,v_j>0$, \eqref{eq:vj infty} gives 
   $\frac{a_j}{8t_j}+v_j\to+\infty$.
   If, up to a subsequence,
   $\frac{a_j}{8t_j}-v_j\to -\infty$, then 
   \eqref{eq:vj infty estimate} can be rewritten as 
   \[
   8t_jv_j=1+a_j(2\mu_1(x_j)-1)+o(1),\quad j\to\infty.
   \]
   However, from Step 1, $w_j\in (\frac{1}{2},\frac{4}{5})$,
   and using $a_j\in [0,\frac{1}{2}]$
   and $\mu_1\in (\frac{1}{4},\frac{3}{4})$, we obtain a contradiction
   \begin{equation*}
          \frac{3}{5}>2w_j-1=8t_jv_j=1+a_j(2\mu_1(x_j)-1)+o(1)\geq \frac{3}{4}+o(1),~j\to\infty.
   \end{equation*}
   Thus, based on \eqref{eq:vj infty}, we have $\frac{a_j}{8t_j}-v_j\to+\infty$.

{\bf Step 3: We prove that 
$v_j=O(1)$
as $j\to \infty$.}

Similar to the proof of \eqref{eq:n=0,1 derivative} and \eqref{eq:mu1 value estimate}, 
taking  logarithms in \eqref{eq:main term compare}, we obtain
\begin{equation}\label{eq:prop33-vertical}
       \frac{\mu_1(x_j)-1/2}{2r_j}
       -\log\frac{\mu_1(x_j)}{1-\mu_1(x_j)}
       =\log\frac{\theta(x_j;t_j)}{\theta(x_j-a_j;t_j)}
       +O(v_je^{-c/t_j}).
       \end{equation}
Here, the improved remainder $O(v_je^{-c/t_j})$ follows in the same way: the difference
between the two sides without the remainder is odd in $v$,
vanishes at $v=0$, and has $v$ derivative $O(e^{-c/t_j})$
by \eqref{eq:scaled-tail-estimate} and \eqref{eq:partialmu1} with $k=1$.

Similar to the proof of 
\eqref{eq:theta first partial x} and using $w_j\pm \frac{a_j}{2}=\frac12+4t_j(v_j\pm \frac{a_j}{8t_j})$, we have the uniform estimate
\begin{equation}\label{eq:log theta w estimate}
       \begin{aligned}
&\log\theta\left(w_j\pm\frac{a_j}{2};t_j\right)
=-\frac1{16t_j}-\frac12\log(\pi t_j)
-4t_j\left(v_j\pm\frac{a_j}{8t_j}\right)^2\\
&\qquad\qquad\qquad
+\log\cosh\left(v_j\pm\frac{a_j}{8t_j}\right)
+O(e^{-c/t_j}),~j\to\infty.
\end{aligned}
\end{equation}
Inserting it into the right-hand side of 
\eqref{eq:prop33-vertical} with $x_j=w_j+\frac{a_j}{2}$, we have 
\begin{equation}\label{eq:mu1 lower bound}
    \begin{aligned}
       &\frac{\mu_1(x_j)-1/2}{2r_j}
       -\log\frac{\mu_1(x_j)}{1-\mu_1(x_j)}\\=&\log\frac{\cosh(\frac{a_j}{8t_j}+v_j)}
{\cosh(\frac{a_j}{8t_j}-v_j)}
-2a_jv_j+O(v_je^{-c/t_j})\\
\ge& 2v_j\left[
\tanh\left(\frac{a_j}{8t_j}-v_j\right)
-a_j-Ce^{-c/t_j}\right]>0.
\end{aligned}   
\end{equation}
Here the inequalities use $\frac{d}{dx}\log\cosh x=\tanh x$, 
the fact that $\tanh(\cdot)$ is increasing, 
$\frac{a_j}{8t_j}-v_j\to+\infty$ as proved in Step 2, and $a_j\le 1/2$.

Since $\frac{\mu-1/2}{2r_j}
       -\log\frac{\mu}{1-\mu}$
is strictly increasing in $\mu\in(1/4,3/4)$ for small 
$r_j>0$
and vanishes at $\mu=1/2$,  \eqref{eq:mu1 lower bound}
gives $\mu_1(x_j)>1/2$ when $j$ is sufficiently large.
Together with 
$\frac{a_j}{8t_j}\pm v_j\to +\infty$, 
the identity
 \eqref{eq:vj infty estimate} becomes
\[
\begin{aligned}
&\frac12\left[
\tanh\left(\frac{a_j}{8t_j}+v_j\right)
-\tanh\left(\frac{a_j}{8t_j}-v_j\right)
\right]-8t_jv_j\\
=&(2\mu_1(x_j)-1)
\left[
\frac{\tanh(\frac{a_j}{8t_j}+v_j)+
\tanh(\frac{a_j}{8t_j}-v_j)}2-a_j
\right]+O(v_je^{-c/t_j})\\
\geq& -Cv_je^{-c/t_j}.
\end{aligned}
\]
Therefore, when $j$ is sufficiently large, using
 $\frac{a_j}{8t_j}\pm v_j>0$, we obtain
\begin{equation}\label{eq:vj estimate}
\begin{aligned}
4t_jv_j
&\le\frac12\left[
\tanh\left(\frac{a_j}{8t_j}+v_j\right)
-\tanh\left(\frac{a_j}{8t_j}-v_j\right)
\right]\\
&\le\frac12\operatorname{sech}^2
\left(\frac{a_j}{8t_j}-v_j\right)
\le Ct_j,
\end{aligned}
\end{equation}
where the last inequality uses \eqref{eq:3-59} and $\mu_1(x_j)\in(\frac{1}{4},\frac{3}{4})$.
Thus $v_j=O(1)$ as $j\to\infty$.

{\bf Step 4: We prove that $L_{r_j}''(w_j;a_j,t_j)>0$, which contradicts
\eqref{eq:bad}.}

By estimating \eqref{eq:vj estimate} again
and using $v_j=O(1)$,
when $j$ is sufficiently large, we obtain
\[
\begin{aligned}
4t_jv_j
\le\frac12\int_{\frac{a_j}{8t_j}-v_j}^{\frac{a_j}{8t_j}+v_j}
\operatorname{sech}^2 u\dd u\le4v_je^{-2(\frac{a_j}{8t_j}-v_j)}
=4v_je^{-a_j/(4t_j)} O(1).
\end{aligned}
\]
Thus, $e^{-a_j/(4t_j)}>c t_j$ for some constant $c>0$.
On the other hand, using \eqref{eq:3-59}, 
there exists a constant $C>0$ such that 
\[
e^{-a_j/(4t_j)}\le
e^{-2(\frac{a_j}{8t_j}-v_j)}
\le\operatorname{sech}^2\left(\frac{a_j}{8t_j}-v_j\right)
\le Ct_j.
\]
Therefore, when $j$ is sufficiently large,  
\begin{equation}\label{eq:prop33-scales}
 ct_j\leq e^{-a_j/(4t_j)}\le Ct_j\quad\Longrightarrow\quad a_j\sim -4t_j \log(t_j).
\end{equation}

Now we make a more delicate change of variables for $x$
with 
\[
x=\frac{1+a_j}{2}+4t_js,
\]
where  $|s|\le M:=\sup_{j\geq 1} v_j $.
For any fixed $\delta_0>0$, \eqref{eq:prop33-scales} gives 
$x\in(1/2-\delta_0,1/2+\delta_0)$ uniformly for $|s|\le M$ and sufficiently large $j$.
Using
$\mu_{1,j}((1+a_j)/2)=1/2$ and \eqref{eq:partialmu1}, we get
\[
 \left|\mu_{1,j}\left(\frac{1+a_j}{2}+4t_js\right)-\frac12\right|
 \le \left|\int_0^s
 \partial_\rho\mu_{1,j}\left(\frac{1+a_j}{2}+4t_j\rho\right)\,d\rho\right|
 \le Cr_j|s|\le CMr_j,
\]
\[
 \left|\partial_s^k\mu_{1,j}\left(\frac{1+a_j}{2}+4t_js\right)\right|
 =(4t_j)^k\left|\partial_x^k\mu_{1,j}(x)\right|\le C_k r_j,\quad 1\leq k\leq 3.
\]
This proves \begin{equation}\label{eq:mu1 C3 estimate}
       \left\|\mu_{1,j}\left(\frac{1+a_j}{2}+4t_j s\right)-\frac12\right\|_{C_s^3([-M,M])}=O(r_j).
\end{equation}

Moreover, applying the chain rule with \eqref{eq:scaled-tail-estimate}
and \eqref{eq:partialmu1}, we can improve the remainder term in
\eqref{eq:prop33-vertical} to $O_{C_s^3}(e^{-c/t_j})$ and obtain
\begin{equation}\label{eq:prop33-vertical-s}
\frac{\mu_{1,j}(x)-1/2}{2r_j}
-\log\frac{\mu_{1,j}(x)}{1-\mu_{1,j}(x)}
=\log\frac{\theta(x;t_j)}{\theta(x-a_j;t_j)}
+O_{C_s^3}(e^{-c/t_j}).
\end{equation}

In the same way, we can also improve the remainder 
term of \eqref{eq:log theta w estimate} to $O_{C_v^3}(e^{-c/t_j})$
and obtain
\[
\begin{aligned}
\log\theta\left(\frac12+4t_j v;t_j\right)=-\frac1{16t_j}-\frac12\log(\pi t_j)
-4t_jv^2+\log\cosh v
+O_{C_v^3}(e^{-c/t_j}).
\end{aligned}
\]
Note that $x=\frac{1+a_j}{2}+4t_js$.
By choosing $v=s\pm \frac{a_j}{8t_j}$, we obtain
the uniform estimate for $|s|\le M$,
\begin{equation}\label{eq:prop33-ratio}
\begin{aligned}
\log\frac{\theta(x;t_j)}{\theta(x-a_j;t_j)}
&=-2a_js+
\log\frac{\cosh(\frac{a_j}{8t_j}+s)}
{\cosh(\frac{a_j}{8t_j}-s)}
+O_{C_s^3}(e^{-c/t_j})\\
&=2(1-a_j)s+O_{C_s^3}(t_j).
\end{aligned}
\end{equation}
Here, the last equality 
follows from \eqref{eq:prop33-scales} and $|s|\leq M$.

Since \eqref{eq:mu1 C3 estimate} implies 
$\log\frac{\mu_{1,j}(x)}{1-\mu_{1,j}(x)}
=O_{C_s^3}(r_j)$, together with \eqref{eq:prop33-ratio} and $r_j=O(t_j)$ because $r_j\le4t_j/3$,
formula \eqref{eq:prop33-vertical-s} implies 
\begin{equation}\label{eq:prop33-mu}
\mu_{1,j}(x)=\frac12+4r_j(1-a_j)s
+O_{C_s^3}(t_j^2).
\end{equation}

Substituting \eqref{eq:prop33-mu} and \eqref{eq:theta first partial x}
(with $w_j$ replaced by $\frac12+4t_js$) into 
\eqref{eq:small-upper-derivative}, and using
\eqref{eq:scaled-tail-estimate} and \eqref{eq:partialmu1} to control
the remainder in $C_s^3$, we obtain
\begin{equation}\label{Lr taylor}
   \begin{aligned}
&4t_jL_{r_j}'\left(\frac12+4t_js;a_j,t_j\right)\\
=&-8t_js+\frac12\left[
\tanh\left(\frac{a_j}{8t_j}+s\right)
-\tanh\left(\frac{a_j}{8t_j}-s\right)
\right]\\
-&4r_j(1-a_j)s
\left[
\tanh\left(\frac{a_j}{8t_j}+s\right)
+\tanh\left(\frac{a_j}{8t_j}-s\right)
-2a_j
\right]+O_{C_s^3}(t_j^2).
\end{aligned}    
\end{equation}
Using \eqref{eq:prop33-scales} and $|s|\leq M$,
we have 
$$\tanh\left(\frac{a_j}{8t_j}\pm s\right)
=1-2e^{-a_j/(4t_j)}e^{\mp2s}
+O_{C_s^3}(t_j^2),$$ then \eqref{Lr taylor} becomes
\begin{equation}\label{eq:prop33-Lprime}
\begin{aligned}
&4t_jL_{r_j}'\left(\frac12+4t_js;a_j,t_j\right)\\
=&-8\bigl[t_j+r_j(1-a_j)^2\bigr]s
+2e^{-a_j/(4t_j)}\sinh(2s)
+O_{C_s^3}(t_j^2).
\end{aligned}
\end{equation}

Differentiating \eqref{eq:prop33-Lprime} three times with respect to $s$
and using \eqref{eq:prop33-scales} yields, for all $|s|\leq M$ and sufficiently large $j$,
\begin{equation}\label{eq:prop33-fourth}
\begin{aligned}
(4t_j)^4 L_{r_j}^{(4)}
\left(\frac12+4t_js;a_j,t_j\right)=16e^{-a_j/(4t_j)}\cosh(2s)+O(t_j^2)
\ge c_1t_j>0.
\end{aligned}
\end{equation}
In particular, $L_{r_j}^{(4)}(w;a_j,t_j)>0$ for any $\frac{1}{2}\le w\le w_j$, because $w_j=\frac12+4t_jv_j$ with $0<v_j\leq M$.

Since $L_{r_j}$ is even and 1-periodic, by a direct calculation similar to
\eqref{eq:even-fourth},  for any $v\in (0,w_j-\frac{1}{2}]$,
\[
\begin{aligned}
&vL_{r_j}''(1/2+v;a_j,t_j)
-L_{r_j}'(1/2+v;a_j,t_j)\\
=&\int_0^{v}u\int_0^u
L_{r_j}^{(4)}\left(\frac12+r;a_j,t_j\right)
\,d r\,du>0.
\end{aligned}
\]
Choosing $v=w_j-\frac{1}{2}$ and using
$L_{r_j}'(w_j;a_j,t_j)=0$, we obtain
$L_{r_j}''(w_j;a_j,t_j)>0$, which contradicts
\eqref{eq:bad}. This completes the proof.
 \end{proof}

 \section{The splitting phenomenon at degenerate critical points}
\label{sec:local-degeneracy}
The purpose of this section is to construct nondegenerate saddle points near
 a degenerate critical point by perturbation. These results will exclude 
 the existence of nontrivial nondegenerate saddle points of $\Theta_S$ and prove 
 that all nontrivial critical points are minimum points in Section 
 \ref{sec:global-exclusion}.

In this section, based on the identity \eqref{eq:matrix-identification},
we use $\Theta_S$ defined in \eqref{eq:def US}
to express the heat kernel $p_{\tau}$.
Specifically, for any $\tau=a+bi\in\mathcal{F}$ and $z\in E_{\tau}$,
let 
\begin{equation*}
A=\begin{pmatrix}1&a\\0&b\end{pmatrix},
\quad
S=tA^{-1}A^{-T},
\quad
q=A^{-1}z.
\end{equation*}
Using \eqref{eq:matrix-identification} and the chain rule,
by direct calculation
we have 
\begin{equation}\label{eq:theta S first derivative}
       \nabla_zp_\tau(z,t)
=\frac1{|\det A|}
A^{-T}\nabla_q\Theta_S(q)
\end{equation}
and 
\begin{equation}\label{eq:theta S second derivative}
   D_z^2p_\tau(z,t)
=\frac1{|\det A|}
A^{-T}D_q^2\Theta_S(q)A^{-1}.    
\end{equation}
Therefore, $z$ is a critical point of $p_{\tau}(\cdot,t)$ 
if and only if $q=A^{-1}z$ is a critical point of $\Theta_S(q)$,
and their Hessians $D_z^2p_{\tau}$ and $D_q^2\Theta_S$
have eigenvalues with the same sign.
Moreover, $z\in 
     \{  0,\frac12,
       \frac{\tau}{2},
       \frac{1+\tau}{2} \}
       \pmod{\Lambda_\tau}$ is a trivial critical point if and only if 
$q=A^{-1}z\in\{ (0,0),
       \left(\frac12,0\right),
       \left(0,\frac12\right),
       \left(\frac12,\frac12\right) 
       \}\pmod{\Z^2}$ is a half-integer point.

Note that for any real symmetric matrix $M=(M_{ij})_{2\times 2}\in \operatorname{Sym}_2(\R)$,
 for sufficiently small $s$, $S+sM$ is positive definite.
 By direct calculation,
for any $m\geq 0$, the \(m\)-th
order directional derivative of $\Theta_S$ satisfies
the important identity
\[
 \frac{d^m}{ds^m}\Theta_{S+sM}(q)
 =\left(\sum_{i,j=1}^2M_{ij}\partial_{q_iq_j}\right)^m
 \Theta_{S+sM}(q).
 \]
Specifically,
the first variation of $\Theta_S$ satisfies
\begin{equation}\label{eq:matrix-heat}
 \left.\frac{d}{ds}\Theta_{S+sM}(q)\right|_{s=0}
 =\sum_{i,j=1}^2M_{ij}\partial_{q_iq_j}\Theta_S(q).
\end{equation}

\subsection{The degenerate critical point with Hessian of rank one}\label{subsec:localization}
In this subsection, we study the case of rank-one Hessian
in Proposition \ref{prop:reverse-off-half}
and Proposition \ref{prop:reverse-half}. To prove Proposition \ref{prop:reverse-off-half},
we need the following lemma.

\begin{lemma}\label{lem:no-straight}
Let $S\in\operatorname{Sym}_2^+(\R)$, $q_0\in\R^2/\Z^2$, and
$v\in\R^2\setminus\{0\}$. For every
nonempty open interval $I=(\alpha,\beta)\subset\R$,
\[
 \nabla\Theta_S(q_0+sv)\not\equiv0,\quad s\in I.
\]
\end{lemma}

\begin{proof}
Suppose $\nabla\Theta_S(q_0+sv)\equiv0$ on $I$, then
 $\Theta_S(q_0+sv)$ is constant on $I$. 
 Since $\Theta_S$ is real analytic, it follows that
\begin{equation}\label{eq:q0 segment}
      \Theta_S(q_0+sv)=\Theta_S(q_0),\quad\forall s\in\R.
  \end{equation}
       
 Suppose $v=(v_1,v_2)^T$, if $v_2\neq 0$ and $v_1/v_2\in \R\setminus\mathbb{Q}$,
 then  by using 
 Kronecker's density theorem, $q_0+\R v$ is dense in $\R^2/\Z^2$, i.e.,
 \begin{equation}\label{eq:segment density}
      \overline{\{q_0+sv\pmod{\Z^2}:s\in\R\}}
 =\R^2/\Z^2.
\end{equation}
  Therefore,  \eqref{eq:q0 segment} and \eqref{eq:segment density} imply
  $\Theta_S(q)\equiv \Theta_S(q_0),q\in\R^2/\Z^2$.
  However, this leads to a contradiction
 \begin{align*}
      0&=\int_{[0,1]^2} \Theta_S(q_0) e^{-2\pi i (1,0)\cdot q}\,dq =\int_{[0,1]^2}\Theta_S(q)e^{-2\pi i(1,0)\cdot q}\,dq
 \\
& =e^{-4\pi^2 (1,0)S(1,0)^T}>0.
  \end{align*}
       
 In the remaining case, $v$ is a nonzero real multiple of a
 primitive integer vector $p=(p_1,p_2)^T$. Rescale the line
  parameter so that its direction is $p$. Take $a,b\in\Z$ such that $p_1b-p_2a=1$. Then
  \[
   M:=\begin{pmatrix}p_1&a\\p_2&b\end{pmatrix}\in SL(2,\Z),
   \quad M^{-1}p=(1,0)^T,
  \]
and the identity
 $\Theta_S(q)=\Theta_{M^{-1}SM^{-T}}(M^{-1}q)$ reduces the
 problem to the direction $(1,0)^T$.
       
 When $v=(1,0)^T$, 
 \eqref{eq:q0 segment} becomes
 $\Theta_{S}(q_0+s(1,0))\equiv \Theta_{S}(q_0),\forall s\in\R$. Let $q_0=(q_{0,1},q_{0,2})$,
 then for any $j\in \Z\setminus\{0\}$, based on the orthogonality of trigonometric functions and \eqref{eq2-9},
 \begin{equation}\label{eq:theta line zero}
      \begin{aligned}
       &0=\int_0^1 \Theta_{S}(q_{0,1},q_{0,2}) e^{-2\pi i j x}\,dx=\int_0^1 \Theta_{S}(x,q_{0,2})e^{-2\pi ij x}\,dx\\
       &=\sum_{n\in\mathbb Z}
e^{-4\pi^2(S_{11}j^2+2S_{12}jn+S_{22}n^2)}
e^{2\pi inq_{0,2}}\\
&=e^{-4\pi^2j^2S_{11}}
\sum_{n\in\mathbb Z}e^{-4\pi^2S_{22}n^2}
e^{2\pi in(q_{0,2}+4\pi iS_{12}j)}\\
&=e^{-4\pi^2j^2S_{11}}
\theta(q_{0,2}+4\pi iS_{12}j;S_{22}).
 \end{aligned}   
 \end{equation} 
According to \cite[\S20.2(iv)]{DLMF}, 
 the zeros of $\theta(z;t)$ with $t>0$ are
 $z=\left(\frac12+k\right)+ \left( \frac12+\ell\right) 4\pi it,\forall k,\ell\in\Z$.
 Consequently, \eqref{eq:theta line zero} implies
 \[
        q_{0,2},~
 \frac{jS_{12}}{S_{22}}\in\frac12+\Z,\quad \forall j\in\mathbb Z\setminus\{0\}.
 \]
 However, $j=1,2$
 cannot both satisfy $\frac{S_{12}}{S_{22}} j\in \frac{1}{2}+\Z$
 simultaneously, which also leads to a contradiction.
       \end{proof}

Now we prove that if $\Theta_S$ has a non-half-integer degenerate critical 
point with rank-one Hessian, then after a perturbation, it 
can split into a nondegenerate saddle point.

\begin{proposition}\label{prop:reverse-off-half}
       Suppose $S_0\in\operatorname{Sym}_2^+(\R)$ and
       $q_0\notin\tfrac12\Z^2/\Z^2$ is a degenerate critical point 
       of $\Theta_{S_0}$ such that 
        \[
         \nabla\Theta_{S_0}(q_0)=0,\quad
         \operatorname{rank}D^2\Theta_{S_0}(q_0)=1.
        \]
       Then there exist $S_j\in\operatorname{Sym}_2^+(\R)$ and
       $q_j\notin\tfrac12\Z^2/\Z^2$ such that
       $S_j\to S_0$, $q_j\to q_0$ and 
       \[
        \nabla\Theta_{S_j}(q_j)=0,\quad
        \det D^2\Theta_{S_j}(q_j)<0,\quad \operatorname{tr}S_j<\operatorname{tr}S_0.
       \]
       It means that $q_j$ is a nondegenerate saddle point of $\Theta_{S_j}$.
\end{proposition}

\begin{proof}
Since $\operatorname{rank}D^2\Theta_{S_0}(q_0)=1$, there exist
orthonormal eigenvectors $v,w\in\R^2$ and $\sigma\in\{-1,1\}$ such that
\begin{equation}\label{eq:degenerate derivative}
     D^2\Theta_{S_0}(q_0)v=0,\quad
 \sigma D_w^2\Theta_{S_0}(q_0)=\lambda>0.  
\end{equation}
Here $D_w^2\Theta_{S_0}(q_0)$ denotes the second directional
derivative along $w$, namely,
\[
 D_w^2\Theta_{S_0}(q_0)
 :=\left.\frac{d^2}{ds^2}\Theta_{S_0}(q_0+sw)\right|_{s=0}
 =w^T D^2\Theta_{S_0}(q_0)w.
\]
For small parameters $s,z$, we define the perturbation of 
$\Theta_{S_0}$ in the directions
$vv^T$ and $vw^T+wv^T$ by
 \begin{equation}\label{eq:def u sz}
     u_{s,z}(x,y):=\sigma\Theta_{S_0+svv^T+z(vw^T+wv^T)}
         (q_0+xv+yw),~x,y\in\R.
 \end{equation}

 Since $$(u_{0,0})_y(0,0)=0,\qquad (u_{0,0})_{yy}(0,0)=\lambda>0,$$
the analytic implicit function theorem gives
 $h_{s,z}(x)$ near the point $(s,z,x,y)=(0,0,0,0)$
satisfying
\begin{equation}\label{eq:rank-one-normal}
 (u_{s,z})_y(x,h_{s,z}(x))=0,\quad
 (u_{s,z})_{yy}(x,h_{s,z}(x))\ge\lambda/2.
\end{equation}
Let $$f_{s,z}(x):=u_{s,z}(x,h_{s,z}(x)).$$ By using \eqref{eq:degenerate derivative} and \eqref{eq:rank-one-normal}, 
we have 
\begin{equation}\label{eq:h f expansion}
        h_{0,0}(0)=h_{0,0}'(0)=f_{0,0}'(0)=f_{0,0}''(0)=0.
\end{equation}

Now, for sufficiently small $\epsilon>0$, we choose 
\begin{equation}\label{eq:local-trace-control}
       S_\epsilon:=S_0-\epsilon^2vv^T+z_\epsilon(vw^T+wv^T)>0,
\end{equation}
where $z_{\epsilon}\to 0$ as $\epsilon\to 0$ is a parameter to be determined.
Since $|v|=1$ and $v\cdot w=0$, we have
$\operatorname{tr}(vv^T)=1$ and
$\operatorname{tr}(vw^T+wv^T)=2v\cdot w=0$.
Hence $\operatorname{tr}S_\epsilon=\operatorname{tr}S_0-\epsilon^2
<\operatorname{tr}S_0$.
The Schur complement formula \eqref{eq:schur} gives 
\begin{equation}
\begin{aligned}
          f'_{s,z}(x)=0
\iff
\nabla u_{s,z}(x,h_{s,z}(x))=0,\quad
\det D^2u_{s,z}
=(u_{s,z})_{yy}f''_{s,z}.
 \end{aligned}
\end{equation}
Therefore, to prove this proposition, we only need to find 
$z_\epsilon\to 0$ and a point $x_{\epsilon}\to 0$ such that 
\begin{equation}\label{eq:reduced-saddle-test}
 f_{-\epsilon^2,z_\epsilon}'(x_\epsilon)=0,\quad
 f_{-\epsilon^2,z_\epsilon}''(x_\epsilon)<0
\end{equation}
and let $S_j=S_\epsilon,\ q_j=q_{\epsilon}:=q_0+x_\epsilon v+h_{-\epsilon^2,z_{\epsilon}}(x_\epsilon)w$.

To prove \eqref{eq:reduced-saddle-test}, we first 
give some equalities.
The matrix heat identity \eqref{eq:matrix-heat} yields
the identities for $u_{s,z}$
\begin{equation}\label{eq:two-parameter-heat}
 \partial_su_{s,z}(x,y)=\partial_x^2u_{s,z}(x,y),\quad
 \partial_zu_{s,z}(x,y)=2\partial_x\partial_yu_{s,z}(x,y).
\end{equation}
For brevity, write 
\begin{equation}\label{eq:uhf}u=u_{0,0}=\sigma\Theta_{S_0}
         (q_0+xv+yw),\quad h=h_{0,0},\quad f=f_{0,0}=u(x,h(x)).\end{equation}
Using \eqref{eq:h f expansion}, we have the Taylor expansion
\begin{align}
 h(x)&=\frac{h^{(r)}(0)}{r!}x^r+O(x^{r+1}),\quad h^{(r)}(0)\ne0,\quad r\ge2,
 \label{eq:valley-order}\\
 f(x)&=f(0)+\frac{f^{(m)}(0)}{m!}x^m+O(x^{m+1}),\quad f^{(m)}(0)\ne0,\quad m\ge3.
 \label{eq:reduced-order}
\end{align}
Here we set $r=\infty$ if $h\equiv0$, and $m=\infty$ if $f\equiv f(0)$. The corresponding Taylor formula above is used only when its order is finite.
Note that $r,m$ cannot both be
$\infty$ simultaneously. Otherwise,
it follows from \eqref{eq:rank-one-normal} and \eqref{eq:uhf} that $u_x(x,0)=u_y(x,0)=0$ on a neighborhood $x\in(-\delta,\delta)$,
which means
$\nabla u(x,0)=0$. Thus, $ \nabla\Theta_{S_0}(q_0+xv)=0,
\forall x\in(-\delta,\delta)$, contradicting
Lemma~\ref{lem:no-straight}.

 In addition, using the Taylor expansion along the curve \(y=h(x)\),
 we have 
\begin{equation}\label{eq:analytic-valley}
       \begin{aligned}
           u(x,y)&=u(x,h(x))+\sum_{n=2}^{\infty} \frac{1}{n!} (\partial_y)^n u(x,h(x))
 (y-h(x))^n\\
 &=f(x)+\tfrac12(y-h(x))^2a(x,y).    
       \end{aligned}
\end{equation}
Moreover, the analytic remainder part satisfies $a(x,y)=\lambda+O(|x|+|y|)$.

By \eqref{eq:two-parameter-heat}, we have
\begin{equation}\label{eq:parameter-coefficient-identities}
 \left.\partial_s^k\partial_z^\ell u_{s,z}(x,y)\right|_{s=z=0}
 =2^\ell\partial_x^{2k+\ell}\partial_y^\ell u(x,y),
 \qquad k,\ell\ge0.
\end{equation}
Thus, the Taylor expansion
in $(s,z)$ is
\begin{equation}\label{eq:parameter-taylor-expansion}
 u_{s,z}(x,y)
 =\sum_{k,\ell\ge0}\frac{2^\ell s^kz^\ell}{k!\ell!}
 \partial_x^{2k+\ell}\partial_y^\ell u(x,y),
\end{equation}
for $(s,z,x,y)$ sufficiently close to $(0,0,0,0)$.
Now we begin to find $z_{\epsilon},x_{\epsilon}$ satisfying 
\eqref{eq:reduced-saddle-test} in two cases.

{\bf Case 1: $m<2r$.} In this case, $f(x)$ is the dominant term in the 
expression for $u(x,y)$ in \eqref{eq:analytic-valley}.

Choose 
\begin{equation*}
       z_\epsilon=0
\end{equation*}
 and let $G(x,s):=(u_{s,0})_x(x,0)$, then using
 \eqref{eq:two-parameter-heat}--\eqref{eq:analytic-valley}
and $m<2r$,
we obtain
\begin{equation}\label{eq:G equation0}
    G_s(x,s) =G_{xx}(x,s),\end{equation}
\begin{align}\label{eq:G equation}
 G(x,0)&=u_x(x,0)=f'(x)+h(x)h'(x)a(x,0)
+\frac12h(x)^2a_x(x,0)
\\
\quad&=f'(x)+O(x^m)= \frac{f^{(m)}(0)}{(m-1)!}x^{m-1}+O(x^m),~x\to 0. \notag 
\end{align}
Here the correction to $f'(x)$ is $O(x^{2r-1})$ when $r<\infty$, and is identically zero when $h\equiv0$. Since $m<2r$, it is $O(x^m)$ in either case.

Take  $x=\epsilon X, s=-\epsilon^2$, then 
\begin{equation}\label{eq:4-15}
      \begin{aligned}
f'_{-\epsilon^2,0}(\epsilon X)-G(\epsilon X,-\epsilon^2)
&=(u_{-\epsilon^2,0})_x(\epsilon X,h_{-\epsilon^2,0}(\epsilon X))
-(u_{-\epsilon^2,0})_x(\epsilon X,0)\\
&=\int_0^{h_{-\epsilon^2,0}(\epsilon X)}
(u_{-\epsilon^2,0})_{xy}(\epsilon X,y)\,dy.
\end{aligned} 
\end{equation}

Using the Taylor expansion of $G(x,s)$ with respect to $s\to 0$
and scaling $x=\epsilon X, s=-\epsilon^2$, we claim that as $\epsilon\to 0$,
\begin{equation}\label{eq: Q convergence}
     \epsilon^{-(m-1)}G(\epsilon X,-\epsilon^2)
 \to Q_0(X)\quad\hbox{in }C^1_{\mathrm{loc}}(\R),   
\end{equation}
where 
\begin{equation}\label{eq:Q expression}
 Q_0(X):=
 f^{(m)}(0)
 \sum_{j=0}^{\lfloor(m-1)/2\rfloor}
 \frac{(-1)^j X^{m-1-2j}}
 {j!(m-1-2j)!},
 \qquad X\in\mathbb R.
\end{equation}

Indeed, set $N=\lfloor(m-1)/2\rfloor$. Since $G_s=G_{xx}$,
Taylor's formula in $s$ gives, uniformly for $X$ in any fixed
compact interval $I\subset\mathbb R$,
\[
 G(\epsilon X,-\epsilon^2)
 =\sum_{j=0}^{N}\frac{(-\epsilon^2)^j}{j!}
 \partial_x^{2j}G(\epsilon X,0)
 +O_{C^1(I)}(\epsilon^{2N+2}).
\]
By analyticity and \eqref{eq:G equation}, for $0\le j\le N$,
\[
 \partial_x^{2j}G(\epsilon X,0)
 =\frac{f^{(m)}(0)}{(m-1-2j)!}
 \epsilon^{m-1-2j}X^{m-1-2j}
 +O_{C^1(I)}(\epsilon^{m-2j}).
\]
Substituting this into the preceding formula and using
$2N+2\ge m$, we obtain
\[
 \epsilon^{-(m-1)}G(\epsilon X,-\epsilon^2)
 =f^{(m)}(0)\sum_{j=0}^{N}
 \frac{(-1)^jX^{m-1-2j}}{j!(m-1-2j)!}
 +O_{C^1(I)}(\epsilon).
\]
This proves \eqref{eq: Q convergence}.
By \cite[(18.5.13)]{DLMF},
\[
 Q_0(X)=\frac{f^{(m)}(0)}{(m-1)!}H_{m-1}(X/2),
\]
where $H_{m-1}$ is the Hermite polynomial.
It has $m-1$ distinct real roots; see \cite[\S18.2(vi)]{DLMF}.
Since $m-1\ge2$, there exists $X_0\in\mathbb R$ such that
\begin{equation}\label{eq:Q zero point}
    Q_0(X_0)=0,\qquad Q_0'(X_0)<0.   
\end{equation}

When $r<\infty$, differentiating formula \eqref{eq:analytic-valley} 
with respect to $y$ at $y=0$ and using 
\eqref{eq:valley-order}, we obtain the estimate
\begin{equation}\label{eq:u partial y at 0}
       u_y(x,0)=-\lambda \frac{h^{(r)}(0)}{r!}x^r+O(x^{r+1}),\quad x\to 0.
\end{equation}
Moreover, using the heat identity $\partial_su_{s,0}=\partial_x^2u_{s,0}$ in 
\eqref{eq:two-parameter-heat}
and \eqref{eq:u partial y at 0}, together with
\eqref{eq:parameter-taylor-expansion}, we take 
the Taylor expansion of $(u_{\cdot,0})_y(\cdot,0)$ to obtain
\begin{equation}\label{eq:u partial y at perturbation}
       (u_{-\epsilon^2,0})_y(\epsilon X,0)
=O_{C_X^2}(\epsilon^r), \quad\epsilon\to 0.
\end{equation}
Here we use $O_{C_X^j}$ to denote $O_{C^j_{\mathrm{loc}}}$ with respect to the variable $X$.

 Based on the implicit function theorem for equation \eqref{eq:rank-one-normal}, 
 and the above estimate \eqref{eq:u partial y at perturbation}, we obtain the estimate
 \begin{equation}\label{eq:4-18}
        h_{-\epsilon^2,0}(\epsilon X)=O_{C_X^2}(\epsilon^r),~\epsilon\to 0.
 \end{equation}

 Similarly, taking the Taylor expansion of $(u_{-\epsilon^2,0})_{xy}(\epsilon X,y)$
 with respect to $y$ and using \eqref{eq:u partial y at perturbation},
 for any $|y|\leq C\epsilon^r$ with a constant $C>0$, we have 
 \begin{equation}\label{eq:4-19}
        (u_{-\epsilon^2,0})_{xy}(\epsilon X,y)
        =O_{C_X^1}(\epsilon^{r-1}),~\epsilon\to 0.
 \end{equation}

Consequently, substituting estimates \eqref{eq:4-18} and \eqref{eq:4-19} into 
\eqref{eq:4-15}, we have 
\begin{equation}\label{eq:reduction-error}
f_{-\epsilon^2,0}'(\epsilon X)-G(\epsilon X,-\epsilon^2)=O_{C_X^1}(\epsilon^{2r-1}),~\epsilon\to 0.
\end{equation}

If $r=\infty$, then $h\equiv 0$ in \eqref{eq:valley-order},
which means $(u_{0,0})_{y}(x,0)\equiv 0$.
Then, using \eqref{eq:parameter-coefficient-identities}, we have 
\[
\begin{aligned}
\left.\partial_s^j(u_{s,0})_y(x,0)\right|_{s=0}
=\partial_x^{2j}(u_{0,0})_y(x,0)=0,\quad \forall j\geq 1.
\end{aligned}
\]
Thus, the Taylor expansion \eqref{eq:parameter-taylor-expansion}
gives $(u_{s,0})_y(x,0)\equiv0$ and then
$h_{-\epsilon^2,0}\equiv0$, so \eqref{eq:4-15} 
directly implies that 
the left-hand side of \eqref{eq:reduction-error} is identically zero.

Therefore, in this case $m<2r$, 
together with \eqref{eq: Q convergence} and 
\eqref{eq:reduction-error}, we obtain
\begin{equation}\label{eq:low-order-comparison}
 \epsilon^{-(m-1)}f_{-\epsilon^2,0}'(\epsilon X)
       =Q_0(X)+o_{C_X^1}(1).
\end{equation}
By \eqref{eq:Q zero point}, there exists $\delta>0$ such that 
$Q_0(X_0-\delta)>0,
Q_0(X_0+\delta)<0$ and $Q_0'(x)<0,~ \forall |x-X_0|<\delta$.
Therefore, \eqref{eq:low-order-comparison} and the intermediate value theorem imply that for 
sufficiently small $\epsilon$,
there exists 
 $x_\epsilon=\epsilon (X_0+o(1))\in \epsilon(X_0-\delta,X_0+\delta)$ such that 
\[
 \begin{gathered}
 f_{-\epsilon^2,0}'(x_\epsilon)=0,\quad
 f_{-\epsilon^2,0}''(x_\epsilon)
      =\epsilon^{m-2}\bigl(Q_0'(x_\epsilon /\epsilon)+o(1)\bigr)<0.
 \end{gathered}
\]
Thus \eqref{eq:reduced-saddle-test} holds.

{\bf Case 2: $m\ge2r$.}
In this case, the contribution of \(h(x)^2\)
to  $u(x,y)$ in \eqref{eq:analytic-valley} cannot be neglected.
The proof of this case is similar to that of Case 1
and we omit some of the similar estimation details.

 Choosing  
$$z_\epsilon=\frac{h^{(r)}(0)}{r!}\epsilon^{\max\{r,3\}}$$
 and setting 
$G(x,s)=(u_{s,0})_{xxy}(x,0)$, then similarly to \eqref{eq:G equation}
and using \eqref{eq:u partial y at 0}, we have 
 \begin{equation}\label{eq:mixed-zero-data}
  G_s(x,s)=G_{xx}(x,s),\quad
  G(x,0)=-\lambda \frac{h^{(r)}(0)}{r!}r(r-1)x^{r-2}+O(x^{r-1}).
 \end{equation}

 If $r\geq 3$, similarly to \eqref{eq:reduction-error} and \eqref{eq: Q convergence}, we have 
 \begin{equation}\label{eq:4-32}
      \epsilon^{-(2r-2)}
 f'_{-\epsilon^2,\frac{h^{(r)}(0)}{r!}\epsilon^r}(\epsilon X)-2 \frac{h^{(r)}(0)}{r!}\epsilon^{-(r-2)}
 G(\epsilon X,-\epsilon^2)
 =O_{C_X^1}(\epsilon),~\epsilon\to 0,
 \end{equation}
Indeed, since $m\ge2r$, \eqref{eq:analytic-valley} gives
$u_x(x,0)=O(x^{2r-1})$. Applying
\eqref{eq:parameter-taylor-expansion} with $z=0$ yields
\[
 (u_{-\epsilon^2,0})_x(\epsilon X,0)
 =O_{C_X^1}(\epsilon^{2r-1}).
\]
The estimates \eqref{eq:4-18}--\eqref{eq:4-19} still hold,
so the same integration in $y$ as in \eqref{eq:4-15} gives
\[
 f'_{-\epsilon^2,0}(\epsilon X)
 =(u_{-\epsilon^2,0})_x(\epsilon X,0)
   +O_{C_X^1}(\epsilon^{2r-1})
 =O_{C_X^1}(\epsilon^{2r-1}).
\]
Since $(u_{s,z})_y(x,h_{s,z}(x))=0$, the chain rule and
\eqref{eq:two-parameter-heat} give
\[
 \begin{aligned}
 \partial_z f_{s,z}(x)
 &=(\partial_z u_{s,z})(x,h_{s,z}(x))
  +(u_{s,z})_y(x,h_{s,z}(x))\partial_z h_{s,z}(x)\\
 &=2(u_{s,z})_{xy}(x,h_{s,z}(x)),
 \end{aligned}
\]
and hence
\[
 \partial_z f'_{s,z}(x)
 =2(u_{s,z})_{xxy}(x,h_{s,z}(x))
  +2(u_{s,z})_{xyy}(x,h_{s,z}(x))h_{s,z}'(x).
\]
At $s=-\epsilon^2$, $z=0$, $x=\epsilon X$, this equals
$2G(\epsilon X,-\epsilon^2)+O_{C_X^1}(\epsilon^{r-1})$.
On the other hand, Taylor's formula for $f'_{s,z}(x)$ in $z$ at $z=0$,
with $s,x$ fixed, gives
\[
 f'_{s,z}(x)
 =f'_{s,0}(x)
 +z\left.\partial_z f'_{s,z}(x)\right|_{z=0}
 +z^2\int_0^1(1-\tau)
 (\partial_z^2 f')_{s,\tau z}(x)\,d\tau.
\]
By joint analyticity, $\partial_z^2 f'$ and
$\partial_x\partial_z^2 f'$ are bounded in a fixed neighborhood
of $(s,z,x)=(0,0,0)$. Since $\partial_X=\epsilon\partial_x$
at $x=\epsilon X$, the remainder is $O_{C_X^1}(z^2)$,
uniformly for small $\epsilon$ and $z$.
Thus, the preceding estimates yield
\[
 \begin{aligned}
 f'_{-\epsilon^2,z}(\epsilon X)
 =O_{C_X^1}(\epsilon^{2r-1})
 +z\bigl[2G(\epsilon X,-\epsilon^2)
       +O_{C_X^1}(\epsilon^{r-1})\bigr]
 +O_{C_X^1}(z^2).
 \end{aligned}
\]
Substituting $z=z_\epsilon=\frac{h^{(r)}(0)}{r!}\epsilon^r$
and dividing by $\epsilon^{2r-2}$, we obtain
 \eqref{eq:4-32}.

Applying
\eqref{eq:parameter-coefficient-identities}--\eqref{eq:parameter-taylor-expansion}
to $G$ and using \eqref{eq:mixed-zero-data}, a similar argument as \eqref{eq: Q convergence} implies that
 \begin{equation}\label{eq:f convergence Q}
    \epsilon^{-(2r-2)}f_{-\epsilon^2,\frac{h^{(r)}(0)}{r!}\epsilon^r}'(\epsilon X)
  \longrightarrow Q(X)\quad\hbox{in }C^1_{\mathrm{loc}}(\R),~\epsilon\to 0.   
 \end{equation}
Here,
\[
\begin{aligned}
 Q(X)
 &=-2\lambda\left(\frac{h^{(r)}(0)}{r!}\right)^2r(r-1)
 \sum_{j=0}^{\lfloor(r-2)/2\rfloor}
 \frac{(-1)^j(r-2)!}{j!(r-2-2j)!}X^{r-2-2j}\\
 &=-2\lambda\left(\frac{h^{(r)}(0)}{r!}\right)^2
 r(r-1)H_{r-2}(X/2).
\end{aligned}
\]
Thus $Q$ has degree $r-2$, has $r-2$ distinct real roots,
and has negative leading coefficient.
Then the largest real root $X_0\in\R$ of $Q$ satisfies $Q(X_0)=0$
 and $Q'(X_0)<0$. The rest of the proof is the same as in Case 1,
 with $x_\epsilon=\epsilon(X_0+o(1))$.

 If $r=2$, the main term of weighted order four in \eqref{eq:analytic-valley} is
 \begin{equation}\label{eq:main term u(x,y)}
        u(x,y)-f(0)\sim \frac{f^{(4)}(0)}{24}x^4+\frac{\lambda}{2}\left(y-\frac{h^{(r)}(0)}{r!}x^2\right)^2.
 \end{equation}
 Similarly to the proofs above in
 \eqref{eq: Q convergence} and \eqref{eq:f convergence Q},
 we can obtain
 \begin{equation}\label{eq:r=2 explicit compute}
        \epsilon^{-3}f'_{-\epsilon^2,\frac{h^{(r)}(0)}{r!}\epsilon^3}(\epsilon X)
 \longrightarrow Q_1(X) \quad\hbox{in }C^1_{\mathrm{loc}}(\R),~\epsilon\to 0, 
 \end{equation}
 where 
 \begin{equation}\label{eq: Q express at r=2}
        Q_1(X)= \frac{f^{(4)}(0)}{6}(X^3-6X)-8\lambda \frac{h^{(r)}(0)}{r!}^2X-4\lambda \frac{h^{(r)}(0)}{r!}^2.
 \end{equation}

To verify \eqref{eq:r=2 explicit compute}, substitute
$x=\epsilon X$, $y=\epsilon^2Y$, $s=-\epsilon^2$,
and $z=h''(0)\epsilon^3/2$ into
\eqref{eq:parameter-taylor-expansion}.
A monomial $x^ay^bs^cz^d$ then has order
$\epsilon^{a+2b+2c+3d}$.
We retain all terms of order $\epsilon^4$;
the higher-order terms, after division by $\epsilon^4$,
tend to zero, together with any fixed number of derivatives,
on compact sets of $(X,Y)$, by joint analyticity.
Then, 
\eqref{eq:parameter-coefficient-identities}--\eqref{eq:parameter-taylor-expansion}
applied to \eqref{eq:main term u(x,y)} give the leading terms
\[
 h_{-\epsilon^2,z_\epsilon}(\epsilon X)
 =\frac{h''(0)}2\epsilon^2(X^2-2)+O_{C_X^2}(\epsilon^3),
\]
and
\[
 \begin{aligned}
 \epsilon^{-4}\bigl(f_{-\epsilon^2,z_\epsilon}(\epsilon X)-f(0)\bigr)
 \longrightarrow{}&
 \frac{f^{(4)}(0)}{24}(X^4-12X^2+12)\\
 &-4\lambda\left(\frac{h''(0)}2\right)^2(X^2+X-1)
 \quad\hbox{in }C^2_{\mathrm{loc}}.
 \end{aligned}
\]
Differentiation in $X$ yields \eqref{eq:r=2 explicit compute}.

 A direct calculation gives $$Q_1(-\sqrt6)=4\lambda \frac{h^{(r)}(0)}{r!}^2(2\sqrt6-1)>0,\quad
   Q_1(0)=-4\lambda \frac{h^{(r)}(0)}{r!}^2<0.$$
  Thus, $Q_1$ has a sign-changing root $X_0\in(-\sqrt6,0)$.
  Moreover,
  this root is simple. In fact, if $f^{(4)}(0)=0$, 
  then $Q_1$ is a nonconstant linear polynomial;
   if $f^{(4)}(0)\neq 0$,
  we have $Q_1''(X)=f^{(4)}(0)X\neq 0$ on $(-\sqrt6,0)$.
  Hence, $Q_1$ also has a root $X_0\in\R$ such that $Q_1(X_0)=0$ 
  and $Q_1'(X_0)<0$. The remaining proof is the same as above.

The proof is complete.
\end{proof}

Proposition \ref{prop:reverse-off-half} proves the splitting 
of non-half-integer degenerate critical points, while the case of degenerate
half-period
 critical points has an alternative property.

\begin{proposition}
       \label{prop:reverse-half}
       Suppose $S_0\in\operatorname{Sym}_2^+(\R)$ and
       $q_0\in\tfrac12\Z^2/\Z^2\setminus\{0\}$ is 
       a nonzero half-integer point and is a degenerate critical point 
       of $\Theta_{S_0}$ such that 
        \[
         \nabla\Theta_{S_0}(q_0)=0,~
         \operatorname{rank}D^2\Theta_{S_0}(q_0)=1,~ \operatorname{spec}(D^2\Theta_{S_0}(q_0))=\{0,\lambda\} ~\textmd{with}~\lambda>0.
        \]
Let $u_{s,z},h_{s,z},f_{s,z}$ be defined in the proof of Proposition 
\ref{prop:reverse-off-half} as in \eqref{eq:def u sz}. Then exactly one 
of the following alternatives holds:
\begin{enumerate}
       \item If $f_{0,0}^{(4)}(0)\le 0$, 
        then there exist $S_j\in\operatorname{Sym}_2^+(\R)$ and
        $q_j\notin\tfrac12\Z^2/\Z^2$ such that
        $S_j\to S_0$, $q_j\to q_0$, and
        \[
         \nabla\Theta_{S_j}(q_j)=0,\quad
         \det D^2\Theta_{S_j}(q_j)<0,\quad \operatorname{tr}S_j<\operatorname{tr}S_0.
        \]
        \item If $f_{0,0}^{(4)}(0)>0$, then there exists a neighborhood $\mathcal U$ of
        $(S_0,q_0)$ such that
        \begin{equation}\label{eq:prop 4.3 condition}
            (S,q)\in\mathcal U,\quad
         q\notin\tfrac12\Z^2/\Z^2,\quad\nabla\Theta_S(q)=0
         \quad\Longrightarrow\quad D^2\Theta_S(q)>0.   
        \end{equation}
       \end{enumerate}
 \end{proposition}

\begin{proof}
(1) Consider the case $f_{0,0}^{(4)}(0)\le 0$. The proof
is similar to Proposition
 \ref{prop:reverse-off-half}, and the only difference is that
 $q_0\in\tfrac12\Z^2/\Z^2\setminus\{0\}$ is 
 a half-integer point.

 As in the proof of Proposition
\ref{prop:reverse-off-half}, to prove this case, 
 we only need to find 
 $z_\epsilon\to 0$ and a point $x_{\epsilon}>0,\ x_{\epsilon}\to 0$ such that 
\begin{equation}\label{eq:f second derivative negative}
  f_{-\epsilon^2,z_\epsilon}'(x_\epsilon)=0,\quad
 f_{-\epsilon^2,z_\epsilon}''(x_\epsilon)<0
 \end{equation}
 and set $$S_j=S_{\epsilon}=S_0-\epsilon^2vv^T+z_\epsilon(vw^T+wv^T),$$ $$q_j=q_{\epsilon}=q_0+x_\epsilon v+h_{-\epsilon^2,z_{\epsilon}}(x_\epsilon)w.$$
 Note that we need $x_{\epsilon}>0$ in this proposition
  such that $q_j\notin\tfrac12\Z^2/\Z^2$.

 Since $q_0$ is 
 a half-integer point, by the periodicity and evenness of $\Theta_{S_0}$,
 we have 
 \[
 u_{0,0}(-x,-y)=u_{0,0}(x,y),~h_{0,0}(-x)=-h_{0,0}(x),~
 f_{0,0}(-x)=f_{0,0}(x).
 \]
In fact, these identities also hold for every sufficiently small $(s,z)$,
since $q_0$ remains a half-integer point and $u_{s,z}(-x,-y)=u_{s,z}(x,y)$.
Uniqueness in the implicit function theorem makes $h_{s,z}$ odd
and $f_{s,z}$ even.

 Thus, the Taylor orders in
 \eqref{eq:valley-order} and \eqref{eq:reduced-order} satisfy
 \begin{equation}\label{eq:half-orders}
  r\in\{3,5,7,\ldots\}\cup\{\infty\},\quad
  m\in\{4,6,8,\ldots\}\cup\{\infty\}.
 \end{equation}
 Moreover, $ (r,m)\ne(\infty,\infty)$ by the same argument as in Proposition 
 \ref{prop:reverse-off-half}.

 {\bf Case 1: $m<2r$.} We also choose $z_{\epsilon}=0$ and we also have
 \eqref{eq:low-order-comparison}
by \eqref{eq:parameter-taylor-expansion} and the same estimates, where $Q_0(X)$ satisfies \eqref{eq:Q expression}
 and has \(m-1\) distinct real roots.

 Because $m\geq 4$ is even, $Q_0(X)$ defined in \eqref{eq:Q expression}
 is an odd polynomial. Hence $Q_0(X)$ has at least \(\frac{m-2}{2}\)
 positive roots.
 
 If $m\geq 6$, then $Q_0(X)$ has at least \(2\)
 positive roots. Thus, there exists $X_0>0$ such that
 \begin{equation}\label{eq:half-low-negative-root}
  Q_0(X_0)=0,\qquad Q_0'(X_0)<0.
 \end{equation}
 The choice of $x_{\epsilon}=\epsilon (X_0+o(1))$ and the remaining proof are the same as in Proposition
\ref{prop:reverse-off-half}.

If $m=4$, then by
$f_{0,0}^{(4)}(0)\le 0$, we have $f_{0,0}^{(4)}(0)<0$.
Note from \eqref{eq:Q expression} that
$ Q_0(X)=\frac{f_{0,0}^{(4)}(0)}{6}(X^3-6X)$.
So, by choosing $X_0=\sqrt{6}$, we also obtain 
\eqref{eq:half-low-negative-root}.

 {\bf Case 2: $m\geq 2r$.}
 We choose 
 \begin{equation*}
               z_\epsilon=\begin{cases}\frac{h^{(r)}(0)}{r!}\epsilon^r,~&r\geq 5,\\
            0,~&r=3.
        \end{cases}
 \end{equation*}

 If $r\geq 5$, we also have \eqref{eq:f convergence Q}
by \eqref{eq:parameter-taylor-expansion} and the same estimates, where
 $Q(X)$ has $r-2\ge 3$ distinct real roots and $\deg Q=r-2$,
 with negative leading coefficient.
Since $Q(X)$ is odd,  
 the largest real root $X_0>0$ of $Q$ satisfies $Q(X_0)=0$
 and $Q'(X_0)<0$. The remaining proof is the same.

 If $r=3$, we also explicitly compute the limiting polynomial.
 Similar to \eqref{eq:main term u(x,y)},
 the main term of weighted order six in 
 \eqref{eq:analytic-valley} is
  \begin{equation}
         u(x,y)-f(0)\sim \frac{f^{(6)}(0)}{6!}x^6+\frac{\lambda}{2}\left(y-\frac{h^{(3)}(0)}{3!}x^3\right)^2.
  \end{equation}
 Then, similarly to the process in \eqref{eq:r=2 explicit compute}
 and \eqref{eq: Q express at r=2},
 we can explicitly compute 
  \begin{equation}
        \epsilon^{-5}  f'_{-\epsilon^2,0}(\epsilon X)
 \longrightarrow Q_2(X) \quad\hbox{in }C^1_{\mathrm{loc}}(\R),~\epsilon\to 0, 
 \end{equation}
 where 
 \begin{equation}
        Q_2(X)
        =6X\left[
               \frac{f^{(6)}(0)}{6!}(X^4-20X^2+60)
        -6\lambda \left(\frac{h^{(3)}(0)}{3!}\right)^2(X^2-4)
        \right].
\end{equation}


Let $0<X_{-}=10-2\sqrt{10}<X_{+}=10+2\sqrt{10}$ 
be the roots of 
\begin{equation*}
X^2-20X+60=0,
\end{equation*}
 we have $Q_2(\sqrt{X_{-}})>0>Q_2(\sqrt{X_{+}})$.
Since $Q_2(X)/(6X)$ is a polynomial of degree
at most two in  $X^2$, it has a simple sign-changing root
$X_0\in(\sqrt{X_-},\sqrt{X_+})$ with $X_0>0$ satisfying $Q_2(X_0)=0$
 and $Q_2'(X_0)<0$. The remaining proof is the same.

(2) Consider the case $f_{0,0}^{(4)}(0)>0$.
Similarly to the definition of $u_{s,z},h_{s,z},f_{s,z}$, we define
\[
 \begin{gathered}
 u_S(x,y)=\Theta_S(q_0+xv+yw),~
 (u_S)_y(x,h_S(x))=0,~f_S(x)=u_S(x,h_S(x)).
 \end{gathered}
\]
Here, the orthonormal eigenvectors $v,w\in\R^2$ satisfy
\eqref{eq:degenerate derivative} as before.

The implicit function theorem implies that $h_S(x),f_S(x),(u_S)_{yy}(x,h_S(x))$ are analytic 
with respect to $S,x$.
Thus, there exists a neighborhood $\mathcal U_0$ of $(S_0,0)$
such that for any $(S,x)\in \mathcal U_0$,
\begin{equation}\label{eq:fS fourth derivative}
    f_S^{(4)}(x)\ge \frac{f_{S_0}^{(4)}(0)}{2}>0,\quad
 (u_S)_{yy}(x,h_S(x))\ge\frac{(u_{S_0})_{yy}(0,0)}{2}=\frac{\lambda}{2}>0.    
\end{equation}

By the periodicity and evenness of $\Theta_{S}$, 
$f_S(\cdot)$ is even. Similar to \eqref{eq:even-fourth}, a
direct calculation gives
\begin{equation}\label{eq:half-fourth-identity}
       \begin{aligned}
       xf_S''(x)-f_S'(x)
       =\int_0^x t\int_0^t  f_S^{(4)}(s)\,ds\,dt.
       \end{aligned}
      \end{equation}

In a neighborhood \(\mathcal{U}\) of \((S_0,q_0)\), 
 for any \((S,q)\in \mathcal{U}\), \(q\) can be decomposed orthogonally as
\(q=q_0+xv+yw\).
 Choose \(\mathcal{U}\) small enough so that 
\((S,x)\in \mathcal U_0\) and the implicit function $h_S(x)$ is unique
for all \((S,q)\in \mathcal{U}\).

 Then, if $(S,q)\in\mathcal U$,
$q\notin\tfrac12\Z^2/\Z^2$, and $\nabla\Theta_S(q)=0$,
the implicit function theorem implies
\(y=h_S(x)\), \(x\neq 0\) and $ f_S'(x)=0$.
Hence,
it follows from 
\eqref{eq:fS fourth derivative}-\eqref{eq:half-fourth-identity} that
 \[
  f_S''(x)>0.
\]
 By the Schur complement formula \eqref{eq:schur}, we have 
\(D^2\Theta_S(q)>0\).

The proof is complete.
\end{proof}

\subsection{The degenerate critical point with vanishing Hessian}
In this subsection, we study the case of vanishing Hessian, i.e.,
$D^2\Theta_{S_0}(q_0)=0$, in Proposition \ref{prop:reverse-rank-zero}. 
Since the Hessian has no nonzero eigenvalues, we cannot use the proof
as in Section \ref{subsec:localization}.

Here we use a new expression 
 of the heat kernel in terms of three theta functions.
For any $\tau\in\mathcal{F}$ with $a>0$ and $t>0$, we define the  
time-scaling variables
\begin{equation}
       t_1(a,b,t):=\frac{(1-a)t}{b^2},\quad
 t_2(a,b,t):=\frac{(a^2+b^2-a)t}{b^2},\quad
 t_3(a,b,t):=\frac{at}{b^2}.
\end{equation}
For any $x_1,x_2,x_3\in\R$,
we define the integral 
\begin{equation}\label{eq:global-overlap}
 I(x_1,x_2,x_3)
 :=\int_0^1\prod_{j=1}^3\theta(\xi-x_j;t_j)\,d\xi.
\end{equation}

Using the expression \eqref{eq:matrix-identification}, for any 
 $z=\nu+\mu\tau\in E_{\tau}$ with $\nu,\mu\in [0,1]$,
we obtain 
 \begin{equation}\label{eq:global-overlap-representation}
       \begin{aligned}
           &bp_\tau(\nu+\mu\tau,t)=\sum_{m,n\in\Z}
       e^{-4\pi^2t[m^2-2amn+(a^2+b^2)n^2]/b^2}
       e^{2\pi i(m\mu+n\nu)}\\
       =&\sum_{m,n\in\Z}e^{-4\pi^2(t_1m^2+t_2n^2+t_3(m-n)^2)}e^{2\pi i(m\mu+n\nu)}=I(\mu,-\nu,0),   
       \end{aligned} 
\end{equation}
Indeed, expanding each theta function by \eqref{eq2-9}, we obtain
\[
\begin{aligned}
I(\mu,-\nu,0)
={}&\sum_{k_1,k_2,k_3\in\mathbb Z}
e^{-4\pi^2(t_1k_1^2+t_2k_2^2+t_3k_3^2)}
e^{2\pi i(-k_1\mu+k_2\nu)}
\int_0^1e^{2\pi i(k_1+k_2+k_3)\xi}\,d\xi\\
={}&\sum_{k_1+k_2+k_3=0}
e^{-4\pi^2(t_1k_1^2+t_2k_2^2+t_3k_3^2)}
e^{2\pi i(-k_1\mu+k_2\nu)}.
\end{aligned}
\]
Taking
$(k_1,k_2,k_3)=(-m,n,m-n)$ gives the last equality of \eqref{eq:global-overlap-representation}.

Moreover, since the theta function $\theta$ has period 1,
the change of variables $\eta=\xi-c$ in \eqref{eq:global-overlap}
gives the translation invariance
\begin{equation}\label{eq:I identity}
     I(x_1,x_2,x_3)=I(x_1+c,x_2+c,x_3+c),\quad\forall x_1,x_2,x_3,c\in\R.  
\end{equation}

Therefore, choosing $P(\mu,\nu,t)=bp_\tau(\nu+\mu\tau,t)$, 
the equalities \eqref{eq:global-overlap-representation} and \eqref{eq:I identity}
imply that 
\begin{equation}\label{eq:global-center-derivatives}
 \begin{gathered}
\partial_{x_1}I+\partial_{x_2}I+\partial_{x_3}I=0,\quad\textrm{at}~ (x_1,x_2,x_3)\in\R^3,\\
\begin{cases} \partial_{x_1}I=\partial_\mu P,~
 \partial_{x_2}I=-\partial_\nu P\\
 \partial_{x_3}I=-\partial_\mu P+\partial_\nu P\end{cases}~\textrm{at}~ (x_1,x_2,x_3)=(\mu,-\nu,0).
 \end{gathered}
\end{equation}
In particular, $\nabla_z p_{\tau}(z,t)=0$ if and only if $$\partial_{x_1}I(\mu,-\nu,0)=\partial_{x_2}I(\mu,-\nu,0)=\partial_{x_3}I(\mu,-\nu,0)=0.$$

Note that \cite[Theorem 1.1]{LWZ} proved that the heat kernel has no
 nontrivial critical point on rectangular tori, so when we study 
 nontrivial critical points, we can assume $a>0$.
The next lemma characterizes the locations of the nontrivial critical points.

\begin{lemma}\label{lem:global-triangle}
For any $\tau\in\mathcal{F}$ and $t>0$,
suppose $\pm (\nu+\mu\tau)$ are    
 nontrivial critical points satisfying
 $(\mu,\nu)\notin\tfrac12\mathbb Z^2/\mathbb Z^2$ and
 $\nabla_z p_{\tau}(\pm (\nu+\mu\tau),t)=0$.
Then exactly one of this pair of 
critical points lies in 
\begin{equation}\label{eq:global-triangle}
      \mathcal T_+
      =\{(\nu,\mu) \pmod{\Z^2}:0<\mu,\nu<1/2,\ \mu+\nu>1/2\}.
 \end{equation}
\end{lemma}

\begin{proof}
For any nontrivial critical point $z=\nu+\mu\tau\in E_{\tau}$ with $\nu,\mu\in [0,1]$,
we let $(x_1,x_2,x_3)=(\mu,-\nu,0)$.

{\bf Step 1:} We first prove that there does not exist an interval
$J=[\alpha,\alpha+\frac{1}{2}]$
such that $x_1,x_2,x_3 \pmod{1}\in J$.
Suppose for contradiction that such $J$ exists. Then  up to modulo 1,
 we choose representatives $x_1,x_2,x_3\in J$.
Let $x_i=\max\{x_1,x_2,x_3\}$, then 
\[
 x_i-x_j\in[0,1/2],\quad \forall j\ne i.
\]

From \eqref{eq:product},
we directly obtain that
$\theta(x,t)=\theta(-x,t)=\theta(x+1,t)$ for all $x\in\R,t>0$,
and \eqref{eq:circle-signs} implies that $\theta(x,t)$ is strictly decreasing on $x\in(0,1/2)$.
Therefore, for any $0<\xi<1/2$,
  \begin{equation}\label{eq:theta compare}
          \theta(\xi+x_i-x_j;t_j)\le\theta(\xi-x_i+x_j;t_j).
  \end{equation}    
Moreover, equality holds if and only if $x_i-x_j\in\{0,1/2\}$.    
  
Taking the derivative of \eqref{eq:global-overlap} and using the evenness 
and periodicity of $\theta$ together with \eqref{eq:theta compare}, 
we obtain
\[
       \begin{aligned}
       \partial_{x_i}I=&-\int_0^1
       \theta_x(\xi-x_i;t_i)
       \prod_{j\ne i}\theta(\xi-x_j;t_j)\,d\xi\\
       =&-\int_{-1/2}^{1/2}
\theta_x(s;t_i)
\prod_{j\ne i}\theta(s+x_i-x_j;t_j)\,ds\\
       =&-\int_0^{1/2}\theta_x(s;t_i)
       \bigg[
       \prod_{j\ne i}\theta(s+ x_i-x_j;t_j)
       -\prod_{j\ne i}\theta(s- x_i+x_j;t_j)
       \bigg]\,d s\le0.
       \end{aligned}
       \]
 Since $\theta_x(s;t_i)<0$ for $0<s<1/2$, it follows that
$\partial_{x_i}I=0$ if and 
 only if $ x_i-x_j\in\{0,1/2\}$, $\forall j\ne i$.
However, \eqref{eq:global-center-derivatives} gives that $\partial_{x_i}I(\mu,-\nu,0)=0$,
which means $\mu,\nu\in\{0,1/2\}$ and yields a contradiction.

{\bf Step 2:}
Step 1 implies that $0,\mu,1-\nu$ modulo $1$
cannot lie in any closed semicircle and so they are distinct. Choose representatives $\mu,\nu\in(0,1)$.
Replacing $\nu+\mu\tau$ by its paired critical point
$(1-\nu)+(1-\mu)\tau$ if necessary, we may assume
\[
 0<\mu<1-\nu<1.
\]
The three consecutive gaps around the circle are therefore
\[
 \mu,\qquad 1-\nu-\mu,\qquad \nu.
\]
Each gap is positive and strictly less than $1/2$.
Indeed, if one gap is at least $1/2$, the complementary
closed arc would have length at most $1/2$ and contain all
three points, contradicting Step 1.
Consequently,
\[
 0<\mu,\nu<\frac12,
 \qquad
 \mu+\nu>\frac12.
\]
Thus $(\nu,\mu)\in\mathcal T_+$ and $(-\nu,-\mu)\in -\mathcal T_+$. Clearly $\mathcal T_+\cap-\mathcal T_+=\emptyset$.
\end{proof}

Now we prove that if $\Theta_S$ has a non-half-period
degenerate critical point with vanishing Hessian, then after a perturbation, it can split into a
nondegenerate saddle point like Proposition \ref{prop:reverse-off-half}.

\begin{proposition}\label{prop:reverse-rank-zero}
       Suppose $S_0\in\operatorname{Sym}_2^+(\R)$ and
       $q_0\notin\tfrac12\Z^2/\Z^2$ is 
       not a half-integer point and is a degenerate critical point 
       of $\Theta_{S_0}$ such that 
       \[
         \nabla\Theta_{S_0}(q_0)=0,\quad
         D^2\Theta_{S_0}(q_0)=0.
        \]
        Then there exist $S_j\in\operatorname{Sym}_2^+(\R)$ and
        $q_j\notin\tfrac12\Z^2/\Z^2$ such that
        $S_j\to S_0$, $q_j\to q_0$ and 
        \[
         \nabla\Theta_{S_j}(q_j)=0,\quad
         \det D^2\Theta_{S_j}(q_j)<0,\quad \operatorname{tr}S_j<\operatorname{tr}S_0.
        \]
       \end{proposition}

\begin{proof}
{\bf Step 1:} 
We first find a unit vector $v$ such that the derivative of 
$ D^2\Theta_{S_0}(q_0)$ in the direction $v$ satisfies
\begin{equation}\label{eq:cubic-hess-det}
 \det D_vD^2\Theta_{S_0}(q_0)<0.
\end{equation}

Based on expression \eqref{eq:matrix-identification} 
and using the change of variables $\tau=a+bi,t>0$ and  $z=Aq_0$
defined in \eqref{eq:coordinate change}, 
together with formulas
\eqref{eq:theta S first derivative} and 
\eqref{eq:theta S second derivative}, we obtain that 
\begin{equation}\label{eq:p two derivative vanish}
      \nabla_z p_\tau(z,t)=0,\quad D_z^2 p_\tau(z,t)=0.
\end{equation}

By the modular invariance in \eqref{eq:modular invariant},
we may change the lattice basis and rescale space and time so that
$\tau=a+ib\in\mathcal F$. The corresponding parameters $t>0$
and $z\in E_\tau$ satisfy \eqref{eq:p two derivative vanish} and
 \begin{equation}\label{eq:bpthetaU}b p_{\tau}(z,t)=\Theta_{tA^{-1}A^{-T}}(A^{-1}z)=\Theta_{U^{-1}S_0U^{-T}}(U^{-1}q_0)=\Theta_{S_0}(q_0),\end{equation} where 
  $U\in GL(2,\mathbb Z)$ is the change-of-basis matrix satisfying
  \begin{equation}\label{eq:modular matrix}
              tA^{-1}A^{-T}=U^{-1}S_0U^{-T},
              \quad
              A^{-1}z=U^{-1}q_0.
  \end{equation}

  Write $z=\nu+\mu\tau$. Since $q_0=U(\nu,\mu)^T\notin\tfrac12\Z^2/\Z^2$,
  we obtain that 
  $z$ is a nontrivial critical point.
  \cite[Theorem 1.1]{LWZ} proved that the heat kernel on 
  rectangular tori has no nontrivial critical point, so $a>0$.
  Note that this proposition is equivalent under $q\to -q$,
  without loss of generality,
    Lemma \ref{lem:global-triangle}
  and
  \eqref{eq:global-overlap-representation} 
 imply that
  \begin{equation}\label{eq:I expression}
         P(\mu,\nu,t)=bp_\tau(\nu+\mu\tau,t)=I(\mu,-\nu,0),\quad
  (\mu,\nu)\in\mathcal T_+.
  \end{equation}

  Let $(x_1,x_2,x_3)=(\mu,-\nu,0)$. The formulas \eqref{eq:global-center-derivatives},
  \eqref{eq:p two derivative vanish},
  and \eqref{eq:I expression} imply that 
  \begin{equation}\label{eq:I derivative zero}
       \partial_{x_i}I(\mu,-\nu,0)=0,\quad
   \partial_{x_i}\partial_{x_j}I(\mu,-\nu,0)=0,\quad\forall 1\leq i,j\leq 3.  
  \end{equation}
   
  Now we define
  \[
  G_1(s):=\partial_{x_3}I(s,-\nu,0),\qquad   G_2(s):=\partial_{x_3}I(\mu,s,0).
  \]
  Then \eqref{eq:I derivative zero} implies that 
  \begin{equation}\label{eq:two multiple zero}
      G_1(\mu)=G_1'(\mu)=0,\qquad
  G_2(-\nu)=G_2'(-\nu)=0.   
  \end{equation}

  Note that 
  \begin{equation}\label{eq:expression G1}
       \begin{aligned}
          G_1(s)&=\int_0^1
  \theta(s-\xi;t_1)
  \bigl[-\theta_x(\xi;t_3)
  \theta(\xi+\nu;t_2)\bigr]\,d\xi\\
  &= \theta(\cdot;t_1)\ast(-\theta_x(\cdot;t_3)\theta(\cdot+\nu;t_2)),  
       \end{aligned}  
  \end{equation}
  which is the solution of the heat equation at time $t_1$
  \[
  \begin{cases}
  \partial_t u(x,t)=\partial_x^2u(x,t),
  & x\in\mathbb R/\mathbb Z,\quad t>0,\\
  u(x,0)=-\theta_x(x;t_3)\theta(x+\nu;t_2),
  & x\in\mathbb R/\mathbb Z.
  \end{cases}
  \]
The initial function $-\theta_x(x;t_3)\theta(x+\nu;t_2)$ has exactly two simple zeros on $x\in\R/\Z$, and the heat flow is never identically zero because its Fourier multipliers are nonzero.
By \cite[Section~2]{Angenent91}, each zero of order $k$ at a positive time splits into $k$ simple zeros at slightly earlier times.
Since the number of zeros cannot increase under the periodic heat flow, $G_1$ has at most two zeros counting multiplicity.

Hence, \eqref{eq:two multiple zero} implies that $G_1''(\mu)\ne0$ and 
that $\mu$ is the unique zero point of $G_1$.
Similarly, $G_2''(-\nu)\ne0$ and $-\nu$ is the unique zero point of $G_2$.

Using the evenness and period-one periodicity of $\theta(\cdot,t)$,
formula
\eqref{eq:expression G1} gives
\begin{equation}\label{eq:G1 value at 0}
     \begin{aligned}
G_1(0)=-\int_0^{1/2}
\theta_x(\xi;t_3)\theta(\xi;t_1)
\bigl[
\theta(\xi+\nu;t_2)
-\theta(\xi-\nu;t_2)
\bigr]\,d\xi<0.
\end{aligned}  
\end{equation}
Here, the last inequality uses $\nu\in (0,\frac{1}{2})$
and the strict monotonicity of $\theta(\cdot,t)$.
Therefore, the formulas \eqref{eq:two multiple zero},
\eqref{eq:G1 value at 0}, and $G_1''(\mu)\neq 0$ imply that $G_1''(\mu)<0$.
For $G_2$, the same argument gives
\[
 G_2(0)=-\int_0^{1/2}\theta_x(\xi;t_3)\theta(\xi;t_2)
 \bigl[\theta(\xi-\mu;t_1)-\theta(\xi+\mu;t_1)\bigr]\,d\xi>0.
\]
Since $-\nu$ is its unique zero of order two, $G_2''(-\nu)>0$.
Thus
\begin{equation}\label{eq:rank-zero-two-signs}
G_1''(\mu)
=\partial_{x_3}\partial_{x_1}^2I(\mu,-\nu,0)
<0<
\partial_{x_3}\partial_{x_2}^2I(\mu,-\nu,0)
=G_2''(-\nu).
\end{equation}

Now take vector $e=(-1,1)^T$, then \eqref{eq:global-center-derivatives} and 
\eqref{eq:rank-zero-two-signs}
 imply
\begin{equation}\label{eq:det P positive}
      \det(D_eD^2P)|_{(\mu,\nu,t)}
=\left((\partial_{x_3}\partial_{x_1}^2I)
(\partial_{x_3}\partial_{x_2}^2I)-
(\partial_{x_3}\partial_{x_1}\partial_{x_2}I)^2\right)|_{(\mu,-\nu,0)}
<0. 
\end{equation}

Note that for any $\eta=(\eta_1,\eta_2)^T\in\R^2$, by \eqref{eq:bpthetaU}, we have 
\begin{align}\label{eq:P transfer to Theta}
      P(\eta,t)&=bp_{\tau}(\eta_2+\eta_1\tau,t)=\Theta_{U^{-1}S_0U^{-T}}\bigl((\eta_2,\eta_1)^T\bigr)\\
      &=\Theta_{S_0}(q_0+UJ(\eta-(\mu,\nu)^T)),\notag
\end{align}
where $J= \begin{pmatrix}0&1\\1&0\end{pmatrix}$ and $q_0=U(\nu,\mu)^T=UJ(\mu,\nu)^T$ is used.
Taking $v=\frac{UJe}{|UJe|}$,  formula \eqref{eq:det P positive}
and \eqref{eq:P transfer to Theta} imply that 
\[
\begin{aligned}
0>\det(D_eD^2P(\mu,\nu,t))=|UJe|^2(\det (UJ))^2\det(D_vD^2\Theta_{S_0}(q_0)).
\end{aligned}
\]
This completes Step 1.

{\bf Step 2:} We construct $S_j,q_j$ satisfying this proposition.

Let $v$ be the unit vector obtained in
Step 1, and let $M=D_vD^2\Theta_{S_0}(q_0)$ be such that $\det M<0$.
Similarly to Proposition \ref{prop:reverse-off-half}, 
for any $\eta\in\R^2$,
we define the perturbation
\begin{equation}\label{eq:perturbation on vanishing Hessian}
     S(\eta)=S_0+\frac12(v\eta^T+\eta v^T),
\quad
F(\epsilon,\eta)
=\nabla\Theta_{S(\eta)}(q_0+\epsilon v).  
\end{equation}

Since 
\(\nabla\Theta_{S_0}(q_0)=0\), we have $F(0,0)=0$.
Moreover, for any $h\in\R^2$,
the matrix heat
identity \eqref{eq:matrix-heat} implies
\begin{equation}\label{eq:F first derivative}
      \begin{aligned}
&D_\eta F(0,0)h=\left.\frac{d}{dt}\right|_{t=0}F(0,th)
=\left.\frac{d}{dt}\right|_{t=0}  \nabla\Theta_{S_0+\frac12(tvh^T+th v^T)}(q_0)\\
=& \frac12\nabla \left( \sum_{i,j=1}^2
(v_i h_j+h_i v_j)
\partial_{q_iq_j}\Theta_{S_0}\right)(q_0)=\nabla D_vD_h\Theta_{S_0}(q_0)
=Mh.
\end{aligned} 
\end{equation}
Thus $D_\eta F(0,0)=M$ is invertible.
Based on the implicit function theorem, there exists 
an analytic function $\eta(\epsilon)$ for $\epsilon\in (-\delta,\delta)$ such that 
\begin{equation}\label{eq:implicit F}
     \eta(0)=0,\quad
F(\epsilon,\eta(\epsilon))=0.  
\end{equation}

Differentiating the identity \eqref{eq:implicit F} with respect to $\epsilon$ at
$\epsilon=0$ and using \eqref{eq:F first derivative}, we obtain
\begin{equation*}
       0=F_\epsilon(0,0)
       +D_\eta F(0,0)\eta'(0)=M\eta'(0)+F_\epsilon(0,0).
\end{equation*}
Since $F_\epsilon(0,0)
=D^2\Theta_{S_0}(q_0)v=0$, we have $\eta'(0)=0$.

Differentiating \eqref{eq:implicit F} once more yields
\[
M\eta''(0)+F_{\epsilon\epsilon}(0,0)=0.
\]
Note that $F_{\epsilon\epsilon}(0,0)
=D_v^2\nabla\Theta_{S_0}(q_0)=Mv$, hence $\eta''(0)=-v$. Therefore,
\begin{equation}\label{eq:eta estimate}
    \eta(\epsilon)
=-\frac{\epsilon^2}{2}v+O(\epsilon^3).   
\end{equation}

We 
choose sequences $S_j=S_{\epsilon_j},q_j=q_{\epsilon_j}$ for $\epsilon_j\to 0$, where
\[
S_\epsilon:=S_0+\frac12(v\eta(\epsilon)^T+\eta(\epsilon) v^T),
\quad
q_\epsilon:=q_0+\epsilon v.
\]
For sufficiently small $\epsilon$, we have $S_\epsilon\in \operatorname{Sym}_2^+(\R)$,
 $S_\epsilon\to S_0$, and
$q_\epsilon\to q_0\notin\tfrac12\Z^2/\Z^2$.
Moreover, \eqref{eq:perturbation on vanishing Hessian} and \eqref{eq:implicit F}
imply that $\nabla\Theta_{S_\epsilon}(q_\epsilon)=0$. Formula
\eqref{eq:eta estimate} implies that 
\[
\operatorname{tr}S_\epsilon
=\operatorname{tr}S_0-\frac{\epsilon^2}{2}
+O(\epsilon^3)
<\operatorname{tr}S_0,~\epsilon\to 0.
\]
Furthermore, using the Taylor expansion and estimate \eqref{eq:eta estimate}, we have
\[
D^2\Theta_{S_\epsilon}(q_\epsilon)
=D^2\Theta_{S_0}(q_0)
+\epsilon D_vD^2\Theta_{S_0}(q_0)
+O(\epsilon^2)=\epsilon M+O(\epsilon^2).
\]
Therefore, 
$\det D^2\Theta_{S_\epsilon}(q_\epsilon)
=\epsilon^2\det M+O(\epsilon^3)<0$ as $\epsilon\to 0$.
These sequences $S_{\epsilon_j},q_{\epsilon_j}$ satisfy the proposition. The proof is complete.
\end{proof}

 \section{Proofs of the main results}
 \label{sec:global-exclusion}
 
 In this final section, we prove Theorem \ref{thm:global-classification}
 and Theorem \ref{thm:smooth-degeneracy}.
 We need the following two lemmas concerning the properties of the
 Hessian at critical points.

\begin{lemma}\label{lem:all-half-positive}
       For any $\tau=a+bi\in\mathcal{F}$ and $t>0$,
       let $z\in\{ \frac{1}{2},\frac{\tau}{2},\frac{1+\tau}{2}\}$
      be a nonzero half-period on the torus $E_{\tau}$,
       then $D^2p_\tau(z,t)$ has at least one positive eigenvalue.
       \end{lemma}

\begin{proof}
       If $a=0$, then
       \eqref{eq:matrix-identification} yields
$p_{bi}(\nu+\mu\tau,t)
=\frac1b\theta(\nu;t)\theta\left(\mu;\frac{t}{b^2}\right).$
       Thus, $\theta''(\frac{1}{2};t)>0$ in \eqref{eq:circle-signs} implies that the Hessian has 
       at least one positive eigenvalue.

       If $a\neq 0$, for any $z=\nu+\mu\tau\in\{ \frac{1}{2},\frac{\tau}{2},\frac{1+\tau}{2}\}$,
       there is the following correspondence with $(\mu,-\nu,0)$.
       \[
\begin{array}{c|c}
(\mu,\nu)&(\mu,-\nu,0)\\ \hline
(\frac12,0)&(\frac12,0,0)\\
(0,\frac12)&(0,-\frac12,0)\\
(\frac12,\frac12)&(\frac12,-\frac12,0)
\end{array}
\]
Thus, based on the period-one periodicity in $\nu,\mu$
and the translation invariance
\eqref{eq:I identity},
we can choose $x_i=x_j=0, x_k=\frac12$ with $\{i,j,k\}=\{1,2,3\}$
such that 
\begin{equation}\label{eq:5-1}
       P(\mu,\nu,t)=bp_\tau(\nu+\mu\tau,t)=I(\mu,-\nu,0)=I(x_1,x_2,x_3).
\end{equation}

Integrating by parts and using the evenness and periodicity, we have 
\[
\begin{aligned}
&\partial_{x_k}^2I (x_1,x_2,x_3)
=\int_0^1 \theta(s;t_i)\theta(s;t_j)\theta_{xx}(s-\tfrac12;t_k)\,ds\\
=&-2\int_0^{1/2}
\dfrac{d}{ds}(\theta(s;t_i)\theta(s;t_j))\theta_x(s-\tfrac12;t_k)\,ds>0.
\end{aligned}
\]
Here, the last inequality follows from the sign \eqref{eq:circle-signs} of $\theta_x$.
     
By \eqref{eq:global-center-derivatives} and \eqref{eq:5-1},
 we have
$D_e^2P(\mu,\nu,t)=\partial_{x_k}^2I (x_1,x_2,x_3)>0$ for some vector $e\ne0$.
Therefore, $D^2 p_{\tau}(\nu+\mu\tau,t)$ also has a positive 
eigenvalue.
\end{proof}

By modular reduction and the Hessian congruence in
\eqref{eq:theta S second derivative}, the same conclusion as Lemma \ref{lem:all-half-positive} holds for
$D^2\Theta_S(q)$ at every nonzero half-integer point $q$ and every
$S\in\operatorname{Sym}_2^+(\mathbb R)$.

\begin{lemma}\label{lem:no-nontrivial-saddle}
       For any $S\in \operatorname{Sym}^+_2(\R)$ and $q\notin\tfrac12\mathbb Z^2/\Z^2$,
    we have
        \begin{equation}
               \nabla\Theta_S(q)=0\ \Longrightarrow\
       \det D^2\Theta_S(q)\ge0.
        \end{equation}
        That is, any nontrivial critical point is not a nondegenerate saddle point.
\end{lemma}

\begin{proof}
       Assume for contradiction and define
\[
\mathcal B=
\left\{
S\in\operatorname{Sym}_2^+(\mathbb R):
\exists q\notin\tfrac12\mathbb Z^2/\mathbb Z^2,\
\nabla\Theta_S(q)=0,\
\det D^2\Theta_S(q)<0
\right\}.
\]
Then $\mathcal B\ne\varnothing$.
By the implicit function theorem for $\nabla \Theta_S(q)$,
the set $\mathcal B$ is open.
Let 
\[
c_*=\inf_{S\in\mathcal B}\operatorname{tr}S\geq 0.
\]
Then there exist a sequence $S_j\in\mathcal B$ and corresponding
 $q_j$ such that 
\[
\operatorname{tr}S_j\to c_*,
\quad
\nabla\Theta_{S_j}(q_j)=0,
\quad
\det D^2\Theta_{S_j}(q_j)<0.
\]

By compactness,
up to a subsequence, there exist $S_*\geq 0$ and $q_*\in\R^2/\Z^2$
such that 
\begin{equation*}
       S_j\to S_*,\quad q_j\to q_*,\quad \operatorname{tr}S_*= c_*.
\end{equation*}

{\bf Step 1.}
We will prove that $S_*>0$. 
If this is not true, since $S_*\ge0$,
we have
\begin{equation}\label{eq:S* det}
      \det S_j\to \det S_*=0. 
\end{equation}

For any $j\geq 1$, let 
\begin{equation}\label{eq:tj estimate}
     \begin{aligned}
\ell_j^2=\min_{\mathbf{m}\in\mathbb Z^2\setminus\{0\}}
\mathbf{m}^TS_j^{-1}\mathbf{m}=\min_{\mathbf{m}\in\mathbb Z^2\setminus\{0\}}
|S_j^{-1/2}\mathbf{m}|^2,\quad t_j:=\ell_j^{-2}.
\end{aligned}  
\end{equation}

Define the lattice $\Lambda_j=\ell_j^{-1}S_j^{-1/2}\Z^2$, then we have 
$$t_j  (\ell_j^{-1}S_j^{-1/2})^{-1}(\ell_j^{-1}S_j^{-1/2})^{-T}=S_j.$$
Thus, formulas \eqref{eq:matrix-identification},\eqref{eq:theta S first derivative}
and \eqref{eq:theta S second derivative} imply that the heat kernel on the torus $\R^2/\Lambda_j$ has 
a nontrivial nondegenerate saddle point. 
It is worth mentioning that \eqref{eq:matrix-identification}
remains valid for any lattice basis $\omega_1,\omega_2$ instead of $1,\tau$.

Note that $\Lambda_j$ has shortest vector of length $1$, so there exists
 $\tau_j=a_j+b_j i\in\mathcal{F}$ 
 such that, 
 up to integer basis transformations, rotations, and reflections,
       $ \begin{pmatrix}
1&a_j\\
0&b_j
\end{pmatrix}\Z^2$ equals $\Lambda_j$.
Thus, based on the discussion in Section \ref{subsec heat kernel representation}, 
 $p_{\tau_j}(\cdot;t_j)$ also has a nontrivial nondegenerate saddle point.

 Since  $\operatorname{tr}S_j$ are bounded, there exists $C>0$ such that 
 $0<S_j< CI$. Thus $S_j^{-1}>C^{-1}I$ and
 $\ell_j^2\ge C^{-1}$, so $t_j\le C$.
 Then
 Proposition \ref{prop:uniform-small} implies that 
 \begin{equation}
        t_j\in [t_*,C]\subset (0,\infty).
 \end{equation}
 Together with 
  Proposition \ref{prop:cusp-estimate}, there exists a constant $B>0$ such that  
 $b_j\leq B$.

 However, \eqref{eq:S* det} and $b_j=|\Lambda_j|=\frac{t_j}{\sqrt{\det S_j}}$
 lead to a contradiction
 \begin{equation*}
        0=\det S_*=\lim_{j\to \infty} \det S_j=\lim_{j\to\infty} \frac{t_j^2}{b_j^2}\geq \frac{t_*^2}{B^2}.
 \end{equation*}

 {\bf Step 2:} We use $S_*>0$ to derive a contradiction.
 By continuity,
 \[
 \nabla\Theta_{S_*}(q_*)=0,\quad
 \det D^2\Theta_{S_*}(q_*)\le0.
 \]
 We will prove in \eqref{eq:zero point-determinant} that 
  \(D^2\Theta_{S_*}(0)<0\),
 so \(q_*\ne0\).

 Since $\mathcal B$ is open, if $S_*\in\mathcal B$,
 then for sufficiently small $\delta>0$, we have $S_*-\delta I\in\mathcal B$,
 which leads to a contradiction
 \begin{equation*}
        \operatorname{tr}(S_*-\delta I)<c_*=\inf_{S\in\mathcal B}\operatorname{tr}S.
 \end{equation*}
 Thus, $S_*\notin\mathcal B$. Because of $S_*>0$ and $S_*\notin \mathcal{B}$, 
 there are two cases for $S_*,q_*$.

 {\bf Case 1:}
 \(q_*\notin\tfrac12\Z^2/\Z^2\) and \(\det D^2\Theta_{S_*}(q_*)=0\).
 Based on Proposition \ref{prop:reverse-off-half} and Proposition
 \ref{prop:reverse-rank-zero}, there exists
 \(S'\in\mathcal B\) satisfying
 \[
 \operatorname{tr}S'
 <\operatorname{tr}S_*=c_*,
 \]
 which leads to a contradiction.

{\bf Case 2:} $q_*\in
\tfrac12\Z^2/\Z^2\setminus\{0\}$.
If $\det D^2\Theta_{S_*}(q_*)<0$, then by the implicit function theorem for 
$\nabla\Theta_S(q)=0$ and $\nabla\Theta_S(q_*)\equiv 0$, we deduce that there 
exists a neighborhood $U\times V$ of $(S_*,q_*)$ such that, for any $S\in U$, 
the equation $\nabla\Theta_S(q)=0$ has the unique solution $q_*$ in $V$.
Thus, $q_j=q_*$ for all sufficiently large $j$, which contradicts
  $q_j\notin\tfrac12\mathbb Z^2/\mathbb Z^2$.
Hence
\begin{equation}\label{eq:S_* degenerate}
      \det D^2\Theta_{S_*}(q_*)=0. 
\end{equation}
Furthermore, based on Lemma \ref{lem:all-half-positive} and \eqref{eq:S_* degenerate},
\[
\operatorname{spec}D^2\Theta_{S_*}(q_*)
=\{0,\lambda\} ~\textmd{with}~\lambda>0.
\]
Therefore, from Proposition \ref{prop:reverse-half},
either there exists some $S'\in\mathcal B$ with
\[
\operatorname{tr}S'<c_*,
\]
or there exists a neighborhood of $(S_*,q_*)$ in which $\Theta_S(q)$ has no nontrivial nondegenerate
saddle points.
These contradict the definition of $c_*$ and $S_j\in\mathcal{B}$, respectively. The proof is complete.
\end{proof}

At last, we now prove Theorem \ref{thm:global-classification}
and Theorem \ref{thm:smooth-degeneracy}.

\begin{proof}[\bf Proof of Theorem \ref{thm:global-classification}]

Based on formulas \eqref{eq:matrix-identification} ,\eqref{eq:theta S first derivative} and 
\eqref{eq:theta S second derivative}, we only need to study 
the critical points of $\Theta_S(q)$ and their properties, where 
$S=tA^{-1}A^{-T}\in\operatorname{Sym}_2^+(\R),q=A^{-1}z$ and 
$A=\begin{pmatrix}
       1&a\\
       0&b
\end{pmatrix}$.

{\bf Case 1:} If $q_0\notin\tfrac12\Z^2/\Z^2$
is a nontrivial critical point,
by Lemma \ref{lem:no-nontrivial-saddle}, we have 
\[
\det D^2\Theta_S(q_0)\ge0.
\]
If $\det D^2\Theta_S(q_0)=0$, then $\operatorname{rank}D^2\Theta_{S}(q_0)\in\{0,1\}$.
Therefore, by Proposition \ref{prop:reverse-off-half}
and Proposition \ref{prop:reverse-rank-zero}, there exist 
$\tilde{S}\in \operatorname{Sym}^+_2(\R)$ and 
 $\tilde{q}\notin\tfrac12\mathbb Z^2/\Z^2$ such that 
        \begin{equation}
               \nabla\Theta_{\tilde{S}}(\tilde{q})=0,\quad
       \det D^2\Theta_{\tilde{S}}(\tilde{q})<0,
        \end{equation}
which contradicts Lemma \ref{lem:no-nontrivial-saddle}.
Thus, any nontrivial critical point $q_0$ satisfies
\begin{equation}\label{eq:nontrivial-positive-determinant}
       \det D^2\Theta_S(q_0)>0,
\quad \ind_{q_0}\nabla\Theta_S=1.
\end{equation}

{\bf Case 2:} At the origin $q=0$, formula \eqref{eq:def US} implies that 
\begin{equation*}
v^TD^2\Theta_S(0)v
=-4\pi^2\sum_{\mathbf{m}\in\mathbb Z^2}
(\mathbf{m}\cdot v)^2e^{-4\pi^2\mathbf{m}^TS\mathbf{m}}<0,~\forall v\in\R^2\setminus\{0\}.
\end{equation*}
Therefore,
\begin{equation}\label{eq:zero point-determinant}
       D^2\Theta_S(0)<0,~
\ind_0\nabla\Theta_S=1,~ \text{$0$ is a nondegenerate maximum point}.
\end{equation}
Moreover, the Fourier series gives
\[
 \Theta_S(0)-\Theta_S(q)
 =\sum_{\mathbf{m}\in\mathbb Z^2}e^{-4\pi^2\mathbf{m}^TS\mathbf{m}}
   \bigl(1-\cos(2\pi \mathbf{m}\cdot q)\bigr)>0,
 \; q\not\equiv0\;(\operatorname{mod}\;\mathbb Z^2).
\]
Thus $0$ is the unique global maximum point.
This proves Conclusion (1).

{\bf Case 3:}
For any nonzero half-period $q_1\in \tfrac12\mathbb Z^2/\Z^2\setminus\{0\}$,
Lemma \ref{lem:all-half-positive} implies that $D^2\Theta_S(q_1)$ has a positive 
eigenvalue.
If $\det D^2\Theta_S(q_1)<0$, then $q_1$ is a nondegenerate saddle point.
If $\det D^2\Theta_S(q_1)>0$, then $q_1$ is a nondegenerate minimum point.

If $\det D^2\Theta_S(q_1)=0$, then 
Proposition \ref{prop:reverse-half} and Lemma 
\ref{lem:no-nontrivial-saddle} imply that 
there exist perturbations $$u(x,y)=u_{0,0}(x,y)=\Theta_{S}
(q_1+xv+yw),\; h(x)=h_{0,0}(x),\; f(x)=f_{0,0}(x)$$
in the neighborhood $V$ of $(x,y)=(0,0)$, as defined in
Proposition \ref{prop:reverse-half}, such that 
\begin{equation}\label{eq:u f inequalities}
 u_y(x,h(x))=0,\quad
 u_{yy}(x,h(x))>0,\quad f(x)=u(x,h(x)),
\end{equation}
and $f^{(4)}(0)>0$. Here $v,w$ satisfy \eqref{eq:degenerate derivative} with $\sigma=1$.
After shrinking $V$, continuity also gives $u_{yy}>0$ throughout $V$, so $u(x,y)\ge u(x,h(x))$, with equality only when $y=h(x)$.

Since
$f$ is even and $f''(0)=0$, we have
\begin{equation}\label{eq:f fourth derivative}
     f(x)-f(0)
=\frac{f^{(4)}(0)}{24}x^4+O(x^6).  
\end{equation}

Therefore, for $(x,y)\in V\setminus\{(0,0)\}$ sufficiently close to $(0,0)$, formulas
\eqref{eq:u f inequalities} and \eqref{eq:f fourth derivative}
imply that 
\begin{equation*}
       \begin{aligned}
           &\Theta_S(q_1+xv+yw)-\Theta_S(q_1)=u(x,y)-u(0,0)\\
           =&(u(x,y)-u(x,h(x)))+(f(x)-f(0))> 0. 
       \end{aligned}
\end{equation*}
Thus, $q_1$ is a strict local minimum point.
When $x$ is sufficiently close to $0$,
it follows from \eqref{eq:f fourth derivative} that $f'(x)=0$ if and
 only if $x=0$, so $q_1$ is an isolated minimum point.
The classical index theorem for 
isolated local minimum points gives $\operatorname{ind}_{q_1}\nabla\Theta_S=1.$

Therefore, there are three cases for \(q_1\):
\begin{equation}\label{eq:half point index}
       \left\{ \begin{aligned}
              &\det D^2\Theta_S(q_1)<0,~ \operatorname{ind}_{q_1}\nabla\Theta_S=-1,~  q_1\text{ is a nondegenerate saddle,}\\
              &\det D^2\Theta_S(q_1)=0,~ \operatorname{ind}_{q_1}\nabla\Theta_S=1,~  q_1\text{ is a degenerate minimum,}\\
              &\det D^2\Theta_S(q_1)>0,~ \operatorname{ind}_{q_1}\nabla\Theta_S=1,~  q_1\text{ is a nondegenerate minimum.}
       \end{aligned}
       \right.
\end{equation}

Let 
\begin{equation*}
       N:=|\{q_0\notin\tfrac12\Z^2/\Z^2: q_0 \text{ is a nontrivial critical point}\}|,
\end{equation*}
\begin{equation*}
       s:=|\{q_1\in \tfrac12\mathbb Z^2/\Z^2\setminus\{0\} : q_1 \text{ is a saddle point}\}|.
\end{equation*}

All critical points are isolated, so the critical points of $\Theta_S$ are finite.
By the Poincar\'e--Hopf theorem (see, e.g., \cite[\S6]{Milnor})
and \eqref{eq:nontrivial-positive-determinant}, \eqref{eq:zero point-determinant},
\eqref{eq:half point index}, we have 
\begin{equation}\label{eq:Poincar hopf}
      0=\chi(\R^2/\Z^2)=\sum_{q:\nabla \Theta_S(q)=0}\ind_{q}\nabla \Theta_S
=1+(3-s)-s+N. 
\end{equation}
Thus $(N,s)=(0,2)$ or $(2,3)$, which means that $p_{\tau}$ has at most one pair of nontrivial
critical points.

When $N=2$, we have $s=3$. Suppose $\pm z$ are the nontrivial critical points.
Since $p_{\tau}$ attains a global minimum on the compact manifold $E_{\tau}$,
only nontrivial critical points can be minimum points.
Since $p_\tau(z,t)=p_\tau(-z,t)$, both points attain the same global minimum.
Moreover, \eqref{eq:nontrivial-positive-determinant}
implies that $\pm z$ are nondegenerate global minimum points.
Conclusion (3) follows from \eqref{eq:nontrivial-positive-determinant},
\eqref{eq:half point index}, and Lemma \ref{lem:global-triangle}.

When \(N=0\), we have \(s=2\), and
Conclusion (2) follows directly from
\eqref{eq:half point index}.
This completes the proof.
\end{proof}

\begin{proof}[\bf Proof of Theorem \ref{thm:smooth-degeneracy}]
For any $q_0\in\tfrac12\Z^2/\Z^2\setminus\{0\}$, 
we first prove that $\mathcal D_{q_0}\ne\varnothing$.

Note that 
\[
\begin{pmatrix}1&0\\1&1\end{pmatrix}
\begin{pmatrix}\frac12\\0\end{pmatrix}
=\begin{pmatrix}\frac12\\\frac12\end{pmatrix},\quad 
\begin{pmatrix}0&1\\1&0\end{pmatrix}
\begin{pmatrix}\frac12\\0\end{pmatrix}
=\begin{pmatrix}0\\\frac12\end{pmatrix}.
\]
Thus, the elements of $\tfrac12\Z^2/\Z^2\setminus\{0\}$ can be 
transformed into one another under the action of $GL(2,\Z)$.
So, there exist
$U_-,U_+\in GL(2,\Z)$ such that
\begin{equation}\label{eq:half period transfer}
       U_-^{-1}q_0=\left(\frac{1}{2},0\right),
\quad U_+^{-1}q_0=\left(\frac{1}{2},\frac{1}{2}\right),~\pmod{\mathbb Z^2}.
\end{equation}

Based on formulas  \eqref{eq:circle-signs} and
 \eqref{eq:half period transfer}, together with $\Theta_{USU^T}(q)=\Theta_S(U^{-1}q)$
and $\Theta_I(x_1,x_2)=\theta(x_1;1)$ $\theta(x_2;1)$, we have
\[
\begin{aligned}
    \det D^2\Theta_{U_-U_-^T}(q_0)=\det D^2\Theta_I(\frac{1}{2},0)<0<\det D^2\Theta_I(\frac{1}{2},\frac{1}{2})=  \det D^2\Theta_{U_+U_+^T}(q_0).   
\end{aligned}
\]
Since $\operatorname{Sym}_2^+(\R)$ is convex and 
$\det D^2 \Theta_{S}(q_0)$ is continuous in $S$, the intermediate value theorem 
implies that there exists $S\in \operatorname{Sym}_2^+(\R)$ such that
 $\det D^2\Theta_{S}(q_0)=0$, i.e., $\mathcal D_{q_0}\ne\emptyset$.

 Next, we prove that $\mathcal D_{q_0}$ is a real analytic embedded hypersurface
 in $\operatorname{Sym}_2^+(\R)$. Define
 \[
 d_{q_0}(S):=\det D^2\Theta_S(q_0) : \operatorname{Sym}_2^+(\R)\longrightarrow \R.
 \]
 For any $S_0\in\mathcal D_{q_0}$, we have $d_{q_0}(S_0)=0$.
 Moreover, by Proposition \ref{prop:reverse-half} and Lemma 
\ref{lem:no-nontrivial-saddle}, similarly to the proof of Theorem \ref{thm:global-classification},
 there exist local perturbations $u(x,y)=\Theta_{S_0}
 (q_0+xv+yw),h(x),f(x)$
 satisfying \eqref{eq:u f inequalities} and 
 \[
 D^2\Theta_{S_0}(q_0)v=0,\quad
 D^2\Theta_{S_0}(q_0)w=\lambda w,
 \quad \lambda>0.
 \]
 Therefore, in the orthonormal basis $(v,w)$, we have 
 \begin{equation}\label{eq:Theta S0 two derivative}
     D^2\Theta_{S_0}(q_0)
 =\begin{pmatrix}0&0\\0&\lambda\end{pmatrix}.
 \end{equation}
Write
\[
D^2\Theta_{S_0+\rho vv^T}(q_0)
=\begin{pmatrix}
a(\rho)&b(\rho)\\
b(\rho)&c(\rho)
\end{pmatrix}.
\]
Then $a(0)=b(0)=0$ and $c(0)=\lambda$. Hence
\[
\left.\frac{d}{d\rho}\right|_{\rho=0}
\det D^2\Theta_{S_0+\rho vv^T}(q_0)
=\left.(a'c+ac'-2bb')\right|_{\rho=0}
=\lambda a'(0).
\]
Note from \eqref{eq:matrix-heat} that
\[
\left.\frac{d}{d\rho}\right|_{\rho=0}
\Theta_{S_0+\rho vv^T}(q)
=\partial_v^2\Theta_{S_0}(q).
\]
Since $a(\rho)=\partial_v^2\Theta_{S_0+\rho vv^T}(q_0)$,
we have
\[
a'(0)
=\partial_v^2
\left(
\left.\frac{d}{d\rho}\right|_{\rho=0}
\Theta_{S_0+\rho vv^T}
\right)(q_0)
=\partial_v^4\Theta_{S_0}(q_0).
\]
Therefore,
\begin{equation}\label{eq:d_q derivative}
\begin{aligned}
Dd_{q_0}(S_0)[vv^T]
=\left.\frac{d}{d\rho}\right|_{\rho=0}
\det D^2\Theta_{S_0+\rho vv^T}(q_0)=\lambda\partial_v^4\Theta_{S_0}(q_0).
\end{aligned}
\end{equation}

 Using the Taylor expansion \eqref{eq:valley-order} of $h(x)$ and formula \eqref{eq:half-orders}, we have 
$h(x)=O(x^3)$ as $x\to 0$. Thus, from the Taylor expansion of $u(x,y)$ at $y=h(x)$, we have 
\[
\begin{aligned}
f(x)-u(x,0)=u(x,h(x))-u(x,0)=O(h(x)^2)=O(x^6),~x\to 0.
\end{aligned}
\]
Therefore, comparing the fourth-order Taylor coefficients of $f(x)$ and $u(x,0)$, we obtain that
\begin{equation}\label{eq:u fourth derivative=f fourth derivative}
 0<f^{(4)}(0)=u_{xxxx}(0,0)=\partial_v^4
 \Theta_{S_0}(q_0).
\end{equation}
Here, $f^{(4)}(0)>0$ follows from Proposition \ref{prop:reverse-half} and Lemma 
\ref{lem:no-nontrivial-saddle}.
Therefore,
\begin{equation}
   \begin{aligned}
Dd_{q_0}(S_0)[vv^T]
=\lambda\partial_v^4
\Theta_{S_0}(q_0)=\lambda f^{(4)}(0)>0.
\end{aligned}    
\end{equation}
Consequently, $Dd_{q_0}(S_0)\ne \vec{0}$ for any $S_0\in \mathcal D_{q_0}=d_{q_0}^{-1}(0)$.
Since $d_{q_0}$ is real analytic,
the implicit function theorem
implies that $\mathcal D_{q_0}=d_{q_0}^{-1}(0)$
is a real analytic embedded hypersurface
in $\operatorname{Sym}_2^+(\R)$.

By Theorem \ref{thm:global-classification} (2),
for any $S_0\in\mathcal D_{q_0}$,
the other two nonzero half-periods are
nondegenerate saddle points.
Hence
\[
\mathcal{D}_{q_0}\cap \mathcal{D}_{\tilde{q}_0}=\emptyset,\quad\forall q_0,\tilde{q}_0\in \tfrac12\Z^2/\Z^2\setminus\{0\},~q_0\neq \tilde{q}_0.
\]
Therefore, for any $\tilde{q}_0\in \tfrac12\Z^2/\Z^2\setminus\{0,q_0\}$, we have
$\det D^2\Theta_{S_0}(\tilde{q}_0)<0$. By continuity,
there exists a neighborhood
$U$ of $S_0$ such that
\[
\det D^2\Theta_S(\tilde{q}_0)<0,\quad
\forall S\in U,~\tilde{q}_0\in \tfrac12\Z^2/\Z^2\setminus\{0,q_0\}.
\]
This means that the other
two nonzero half-periods remain
nondegenerate saddle points when $S\in U$.
Hence \eqref{eq:p tau critical number} follows immediately from Theorem \ref{thm:global-classification}
and \eqref{eq:heat kernel theta function equality}.

Finally, let $S_j\in U_-$ tend to $S_0$ and let $q_j$ be a nontrivial critical point of $\Theta_{S_j}$.
By compactness, every subsequential limit of $q_j$ is a critical point of $\Theta_{S_0}$.
All critical points of $\Theta_{S_0}$  other than $q_0$ are nondegenerate half-integer points.
Near each of them, the implicit function theorem gives a unique critical point for nearby $S$, namely the same fixed half-integer point.
Thus none can be a limit of the nontrivial points $q_j$, which implies $q_j\to q_0$.
Hence the pair splits from $q_0$.
\end{proof}

\subsection*{Acknowledgments}  Z. Chen is supported by National Key R\&D Program of China (No. 2023YFA1010002) and NSFC (No. 12222109). 
S. Zhang is supported by the Postdoctoral Fellowship Program and China Postdoctoral Science Foundation under Grant Numbers BX20250062 and 2026M793381. The AI tools were used
during the preparation of this manuscript. The authors assume full responsibility
for the mathematical validity, accuracy, and integrity of the proofs presented in the
manuscript.

\subsection*{Data availability statement}  No data was used for the research described in the article.

\subsection*{Conflict of interest} There is no conflict of interest.

\end{document}